\RequirePackage{ifluatex}

\documentclass[12pt]{article}
\usepackage[margin=1.15in]{geometry}
\usepackage{graphics}
\usepackage{ifpdf}
\usepackage{amsmath,amssymb}
\usepackage{hyperref}
\usepackage{epsfig}
\usepackage{epstopdf}
\usepackage{upgreek}
\usepackage{authblk}
\usepackage[maxbibnames=10, citestyle=numeric-comp,bibstyle=numeric, giveninits=true, backend=bibtex, sorting=none]{biblatex}
\usepackage[font=small,labelfont=bf]{caption}
\usepackage{titlesec}
\usepackage[titletoc,toc,title]{appendix}
\usepackage{abstract}  
\usepackage{setspace}
\usepackage[]{subfig}
\usepackage{bm}
\usepackage{float}

\DeclareMathAlphabet      {\mathbfit}{OML}{cmm}{b}{it}

\usepackage{titlesec}
\titleformat{\paragraph}
{\normalfont\normalsize\bfseries}{\theparagraph}{1em}{}
\titlespacing*{\paragraph}
{0pt}{3.25ex plus 1ex minus .2ex}{1.5ex plus .2ex}

\makeatletter
\def\@seccntformat#1{\@ifundefined{#1@cntformat}%
	{\csname the#1\endcsname\quad}
	{\csname #1@cntformat\endcsname}
}
\makeatother

\AtEveryBibitem{%
	\clearfield{isbn}
	\clearfield{issn}
	\clearfield{month}
}
\renewbibmacro{volume+number+eid}{%
	\printfield{volume}%
	\setunit{\space}%
	\printfield{number}%
	\setunit{\space}%
	\printfield{eid}
}
\DeclareFieldFormat[article]{number}{\mkbibparens{#1}}
\DeclareFieldFormat[article]{volume}{\mkbibbold{#1}}
\bibliography{Source_File}

\title{\large{\textbf{Varying-order NURBS discretization: An accurate and efficient method for isogeometric analysis of large deformation contact problems}}}
\author{Vishal Agrawal}
\author{Sachin S. Gautam\footnote{Corresponding author email: ssg@iitg.ac.in}}
\affil[]{\small{\textit{{~~~~~~~~~~~~~~Department of Mechanical Engineering, \newline
				Indian Institute of Technology Guwahati, Guwahati, Assam, India, 781039}}}}
\affil{\small {Published\footnote{This pdf is the personal version of an article (with minor corrections) whose final publication is available at  https://www.journals.elsevier.com/computer-methods-in-applied-mechanics-and-engineering} ~in: \textit{Computer Methods in Applied Mechanics and Engineering} \newline 
		\hspace*{1.1cm} DOI:\href{https://www.sciencedirect.com/science/article/pii/S0045782520303108?dgcid=author}{10.1016/j.cma.2020.113125}, \newline Submitted: 27 Nov. 2019, Accepted: 8 May 2020}}
\date{}

\begin{document}
\sloppy	
\maketitle
\vspace{-5em}
	
\renewenvironment{abstract}
{\begin{quote}
	\noindent \rule{\linewidth}{.5pt}\par{\bfseries \abstractname.}}
	{\medskip\noindent \rule{\linewidth}{.5pt}
	\end{quote}
}	
\begin{abstract}	

A novel varying-order based NURBS discretization method is proposed to enhance the performance of isogeometric analysis (IGA) technique within the framework of computational contact mechanics. The method makes use of higher-order NURBS for the contact integral evaluations. The minimum orders of NURBS capable of representing the complex geometries exactly are employed for the bulk computations. The proposed methodology provides a possibility to refine the geometry through controllable order elevation for isogeometric analysis. To achieve this, a higher-order NURBS layer is used as the contact boundary layer of the bodies. The NURBS layer is constructed using different surface refinement strategies such that it is accompanied by a large number of additional degrees of freedom and matches with the bulk parameterization.

The capabilities and benefits of the proposed method are demonstrated using the two-dimensional frictionless and frictional contact problems, considering both small and large deformations. The results with the existing fixed-order based NURBS discretizations are used for comparisons. Numerical examples show that with the proposed method, a much higher accuracy can be achieved even with a coarse mesh as compared to the existing NURBS discretization approach. It exhibits a major gain in the numerical efficiency without the loss of stability, robustness, and the intrinsic features of NURBS-based IGA technique for a similar accuracy level. The simplicity of the proposed method lends itself to be conveniently embedded in an existing isogeometric contact code after only a few minor modifications.
		
			
\end{abstract}

\section{Introduction} \label{sec:Introduction}
In the last two decades, non-uniform rational B-splines (NURBS) based isogeometric analysis~\cite{hughes2005} approach has been established as an advantageous computational technology for various classes of problems, especially to those where the geometric approximation and inaccurate capture of the driving forces at the interface adversely influence the accuracy of the result. For example, in cohesive zone modelling problems~\cite{Verhoosel2010, Dimitri2014, Dimitri2017}, estimation of the fatigue life of interfacial cracks in the functionally graded material plates~\cite{Bhardwaj2015, Bhardwaj2016, Singh2019}, evaluation of through-thickness cracks in homogeneous and isotropic plates~\cite{Singh2018a}, in three-dimensional problems having straight and curved crack fronts~\cite{Singh2018b}, and for the analysis of cracks in the functionally graded piezoelectric materials~\cite{Singh2020}, etc. This is attributed to the distinguished intrinsic features of its underlying basis function, viz. the ability to represent complex geometry exactly even with a coarse mesh, variational diminishing and convex-hull properties, tailorable inter-element continuity, and non-negativeness~\cite{hughes2005, Cottrell2009}. Contact modeling belongs to one of these classes that have particularly been benefited from the aforementioned features of the IGA technology. Higher-order smoothness of the NURBS discretized geometry directly provides the continuous normal vector field across the boundary of the contact elements. Thus, it eliminates the need for additional contact surface smoothening approaches~\cite{Padmanabhan2001, El-Abbasi2001, Wriggers2001, Krstulovi2002, Al-Dojayli2002, Stadler2003} often utilized in the context of the finite element method.

The first investigations on the treatment of contact problems using isogeometric analysis are conducted by Temizer et al.~\cite{Temizer2011, Temizer2012}, Lu~\cite{Lu2011}, and De Lorenzis et al.~\cite{DeLorenzis2011, DeLorenzis2012}. 

Temizer et al.~\cite{Temizer2011} applied the isogeometric analysis to three-dimensional thermomechanical frictionless contact using the mortar contact formulation. For the regularization of the contact constraints, the penalty method is adopted. They demonstrated that NURBS based discretization achieves superior results in terms of quality and robustness over its counterpart Lagrange-polynomial based discretization. However, a very fine mesh is still required to get the result that closely matches with the exact solution. Moreover, the non-mortar contact formulation is found to be over-constrained in nature and leads to inaccurate results. Lu~\cite{Lu2011} combined the isogeometric analysis with the segment-to-segment based contact formulation presented by Papadopoulos and Taylor~\cite{Papadopoulos1992}. He showed that intrinsically smooth NURBS discretized geometry alleviates the non-physical oscillations of the contact forces and is capable of accurately describing the intricate mechanics of smooth materials such as fabrics. De Lorenzis et al.~\cite{DeLorenzis2011} introduced a mortar-based contact formulation in the context of IGA for the investigation of two-dimensional large deformation contact with Coulomb friction. They demonstrated that the magnitude of the non-physical oscillations of the normal and tangential contact reaction forces reduce on increasing the interpolation order of NURBS discretization. They also showed that the oscillation error in the Hertzian contact stress at the boundary of the contact region is reduced if a very fine mesh is used. Moreover, the distributions of the contact responses obtained with the Lagrange-polynomial based discretizations are found to be highly sensitive to the interpolation order. Later, De Lorenzis et al.~\cite{DeLorenzis2012} applied IGA to three-dimensional large deformation frictionless contact with the mortar method. The contact constraints are enforced using the augmented-Lagrangian approach~\cite{ALART1991}. Temizer et al.~\cite{Temizer2012} introduced a three-dimensional mortar-based frictional contact formulation as an extension of works in~\cite{Temizer2011, DeLorenzis2011}. Again, it is shown that on increasing the interpolation order of the NURBS discretization the smoothness of contact forces improves across the contact interface. A point-to-segment contact formulation for two-dimensional frictionless contact is presented by Matzen et al.~\cite{Matzen2013}, and recently extended to a weighted point-based contact approach for three-dimensional contact~\cite{Matzen2016}, both with the IGA. Many researchers have extensively studied the application of mortar method in the context of IGA~\cite{Kim2012, Temzier2013Multiscale, Dittmann2014, BRIVADIS2015, SEITZ2016, Duong2018}, e.g. Brivadis et al.~\cite{BRIVADIS2015} investigated the isogeometric mortar method from the theoretical as well as the numerical point of view. Seitz et al.~\cite{SEITZ2016} developed a dual mortar method for isogeometric contact analysis and analyzed the spatial rate of convergence for the mesh tying and contact problems. Recently, Duong et al.~\cite{Duong2018} introduced a segmentation-free isogeometric mortar contact formulation. A number of works have explored the application of collocation methods to different contact problems using IGA, e.g.~\cite{KRUSE2015, DELORENZIS2015, Weeger2018}. De Lorenzis et al.~\cite{DeLorenzis2014} presented a comprehensive overview on the growth of various isogeometric based treatment procedures and their advantages for contact mechanics as compared to traditional finite element analysis.

From the above-reviewed literature on contact modeling using IGA, it is evident that NURBS-based discretization delivers significantly superior performance in terms of accuracy, stability, and robustness over Lagrange-polynomial based discretization. This is due to the aforementioned distinguished intrinsic features of the NURBS over the Lagrange polynomials. However, the application of the existing NURBS-based discretization approach to isogeometric contact analysis is computationally expensive since due to the rigid tensor nature of the NURBS structures, the underlying mesh can be refined only in a global manner. Moreover, the interpolation order of NURBS functions employed for the discretization of the contact boundary layer and for the remaining bulk domain of the body can only be elevated uniformly. From the analysis point of view, this may not be desirable since the accuracy of the contact solution is primarily governed by the description of the contact interface. According to literature, to enable the local mesh refinement in the context of isogeometric contact, T-splines~\cite{BAZILEVS2010Tspline, DORFEL2010, Scott2011, Dimitri2014Tspline, Dimitri2014, Dimitri2017}, NURBS-based hierarchical refinement~\cite{Temizer20161, Temizer20162}, and locally refined (LR) NURBS~\cite{Sauer2017} approaches have been adopted. Dimitri et al.~\cite{Dimitri2014Tspline, Dimitri2014, Dimitri2017} employed the T-splines based discretizations for analyzing the frictionless contact, cohesive/contact interface modeling, and mixed-mode debonding problems. Temizer and Hesch~\cite{Temizer20162}, and Hesch et al.~\cite{ Temizer20161} used the hierarchical NURBS for the treatment of frictionless and frictional contact, respectively. Zimmermann and Sauer~\cite{Sauer2017} introduced LR NURBS for the analysis of contact computations of solids and membranes. However, the idea to refine the geometry through a controllable order elevation strategy remains unexplored, as already noted by Temizer et al.~\cite{Temizer2012}.	

In the context of finite element contact analysis, a considerable research progress has been made that integrates the intrinsic features of NURBS and higher-order Lagrange polynomials with the standard FE discretization~\cite{Corbett2014, Corbett2015, RASOOL2016182, MALEKIJEBELI2018, Otto2018, DIAS2019}. Corbett and Sauer~\cite{Corbett2014, Corbett2015} introduced a NURBS-enriched contact element formulation that combines the geometric smoothness of NURBS with the efficiency characteristic of the FE discretization. The potential contact layer of each finite element is locally replaced by a NURBS layer, resulting in ``NURBS-enriched contact finite elements." Their work is based on the contact element enrichment strategy of Sauer~\cite{Roger2011_enrichment, Sauer2013}, which used higher-order Lagrange and Hermite functions on the contact surface. Konyukhov and Schweizerhof~\cite{Konyukhov2009, Konyukhov2013, Konyukhov2015} introduced the hierarchical enrichment of the Lagrange functions space, which allows the anisotropic refinement of the contact boundary layer with the higher-order Lagrange polynomial, while the linear order polynomials remain fixed for the construction of the interior domain of a contact body. Later, Maleki-Jebeli et al.~\cite{MALEKIJEBELI2018} proposed a hybrid isogeometric-finite element discretization method where the advantages of NURBS are exploited in two-dimensional cohesive interface contact/debonding. The bulk domain is given by FE discretization. The transition from NURBS to FE discretizations is carried out using the so-called ``transition elements", initially presented in~\cite{Corbett2014}. Otto et al.~\cite{Otto2018, Otto2019} proposed a coupled FE-NURBS discretization approach where an auxiliary NURBS layer is placed between the contact zone of higher-order FE discretized contacting bodies. In order to tie the NURBS layer with FE discretization, the pointwise and mortar mesh tying approaches are used. Recently, Dias et al.~\cite{DIAS2019} presented a higher-order mortar-based contact element, where the contact interface and bulk are discretized using the  $ hp- $FEM. In contrast to the enrichment approach of Corbett and Sauer~\cite{Corbett2014, Corbett2015}, they used higher-order Lagrange functions as a basis for the discretization of the contact interface, domain, and approximation of the solution field. The NURBS are only used to map the position of FE contact layer nodes.

On the other hand, in the context of isogeometric contact analysis, no such effort that directly allows the controllable order elevation of NURBS discretized structures and accordingly complements the performance of IGA while fully retaining its intrinsic key features has been devoted. 

In the present work, we thus propose a novel varying-order (VO) based NURBS discretization method to improve the performance of NURBS-based IGA technique in the framework of computational contact mechanics. In the proposed method, the user-defined higher-order NURBS functions are employed for the contact computations. The minimum-order of NURBS capable of representing the complex geometries exactly are used for the description of overall bulk domain that does not come into contact. To achieve the VO based NURBS discretization, a new higher-order NURBS layer is used as the contact boundary layer of an initially NURBS discretized geometry. This avoids the application of higher-order NURBS in the region away from the contact interface. The layer is constructed using different surface refinement strategies in such a manner that it is accompanied by a large number of additional degrees of freedom across the contact interface and matches with the bulk parameterization. The proposed varying-order based NURBS discretization method is accordingly denoted by N$ _p-$N$_{p_c} $, where N$_p $ is the order of the NURBS functions utilized for the description of the bulk domain and N$_{p_c}, (p_c>p) $ is the order of NURBS used for the contact boundary layer. Due to its simplicity, the proposed approach can be conveniently embedded into an existing IGA contact code after only a few modifications.

In the this work, the proposed methodology is applied to two-dimensional frictionless and frictional contact problems, considering both small and large deformations. A classical full-pass version of the contact algorithm introduced by Sauer and De Lorenzis~\cite{unbiased2013, unbiased2015} is combined with the VO based NURBS discretization method in this work as a first step towards developing a simple yet computationally efficient technique for isogeometric contact analysis. Such a contact algorithm has previously introduced for Lagrange-polynomial based discretizations~\cite{Fischer2005, DIAS2019}, fixed-order based NURBS  discretizations~\cite{Temizer2011, DeLorenzis2011, DeLorenzis2012, Dimitri2014Tspline, Dimitri2014, Dimitri2017} and recently for T-spline based discretizations~\cite{Dimitri2014Tspline, Dimitri2014, Dimitri2017}. The penalty method is employed for the regularization of impenetrability and sticking contact constraints.

The remainder of the paper is structured as follows. In Section~\ref{sec:NURBS_modelling}, the existing geometric modeling procedure using NURBS and the strategies used for the refinement of a NURBS discretized geometry are briefly discussed. Section~\ref{sec:Formulation} describes the mathematical formulation for a large deformation contact between two deformable bodies considering friction. The varying-order based NURBS discretization technique and its implementation into the existing IGA contact code are fully described in Section~\ref{sec:global_strategy}. In Section~\ref{sec:numerical_example}, the performance of the proposed method is demonstrated w.r.t. fixed-order based NURBS discretizations using various numerical examples. Finally, Section~\ref{sec:conclusion} concludes the paper with future directions.

\section{Preliminaries} \label{sec:NURBS_modelling}
In this section, we briefly discuss the existing NURBS based discretization procedure used for the continuum and the contact surface in the context of IGA. Next, an overview of the different strategies used for the refinement of a NURBS geometry is presented. The associated issues are also highlighted. For the detailed description we refer to the monographs by Pigel and Tiller~\cite{nurbsbook} and Cottrell et al.~\cite{Cottrell2009}.

\subsection{Geometric modeling using NURBS} \label{sec:NURBS_discretization}	
The IGA has emerged as a successful integration of computer-aided-design (CAD) and finite element analysis (FEA) technique~\cite{Cottrell2009}. It utilizes the CAD polynomials as a basis for the modeling of complex shaped geometries exactly and approximation of unknown solution fields. In current work, NURBS functions are used to describe the continuum as well as contact surface of a geometry.

It is known that NURBS are built from B-splines~\cite{nurbsbook}. Thus, we begin with a brief introduction to B-splines. For a given knot vector $ \Xi^i $ along the $ \xi^i~(i=1,2) $ parametric direction, a B-spline function of interpolation order\footnote{In this work, the linear, quadratic, cubic, etc. piecewise polynomials are denoted with the functions of order $ p = 1, 2, 3, etc., $ respectively. This is as per the notation introduced by Hughes et al.~\cite{hughes2005}. The term ``order'' is usually referred to as ``degree'' in the literature based on computational geometry.} $ p_i $ is defined using the following Cox-de Boor recursive relation~\cite{nurbsbook}:
\begin{eqnarray}
\textnormal{for $p_i=0$, \quad}
N_{l,0}(\xi^i) &=& \left\{\begin{array}{lr}
1, & \textrm{ if }~ \xi_l^i\leq \xi^i < \xi_{l+1}^i\,, \\
0, & \textrm{otherwise}\,.
\end{array} \right.  \\ \newline
\textnormal{for $p_i>0$, \quad}
N_{l,p_i}(\xi^i) &=& \frac{\xi^i-\xi_{l}^i}{\xi_{l+p_i}^i-\xi_l^i} N_{l,p_i-1} +   \frac{\xi_{l+p_i+1}^i-\xi^i}{\xi_{l+p_i+1}^i- \xi_{l+1}^i} N_{l+1,p_i-1}\,,
\end{eqnarray}
where $ N_{l,p_i}(\xi^i) \geq 0 $ and $ \xi_l^i \in \mathbb{R} $ is the $ l^{\textrm{\scriptsize{th}}} $ knot of the knot vector $ \Xi^i = \{\xi_1^i, \xi_{2}^i, \dots, {\xi_{n_i+p_i+1}^i} \} $. A knot vector is a set of non-decreasing values of parametric coordinates, i.e. $ \xi_l^i \leq \xi_{l+1}^i $. Here, $ n_i $ denotes the total number of control points for $ p_i $ order of B-spline defined along the $ \xi^i $ parametric direction. A $ p^{\textrm{\scriptsize{th}}}_i $ order of B-spline function offers $ C^{p_i-m} $ continuity across each knot $ \xi^i$, where $ m $ is the knot multiplicity. In an \textit{open} knot vector, which are commonly adopted in IGA, the first and last knot entries are repeated by $ p_i+1 $ times that make the NURBS functions interpolatory at the end knots.

NURBS are the projective transformation of the B-splines~\cite{nurbsbook}. Thus, NURBS additionally utilize the weight values in their construction as compared to B-splines. For a given knot vector $ \Xi^i $, a $ p^{\textrm{\scriptsize{th}}}_i $ order of univariate NURBS polynomial along the $ \xi^i $ parametric direction is defined as~\cite{nurbsbook}
\begin{equation}  \label{eq:NtoR}
R_{l}^{p_i}(\xi^i) = \frac{w_{l}^i}{ \sum_{A=1}^{n_i} w_{A}^i N_{A,p_i}(\xi^i)} N_{l,p_i}(\xi^i)\,,	 
\end{equation}
where $ R_{l}^{p_i}(\xi^i) \geq 0~~\forall\, \xi^i\in \Xi^i\,,$ and $ w_l^i > 0 $ is the weight value. For a specified control points vector $ \mathbf{X}^i = \{\mathbf{X}_{l}^i\}_{l=1}^{n_i} \,,$ a $ p_i^{\scriptsize{\textrm{th}}} $ order of NURBS curve defined along the $ \xi^i $ parametric direction can be constructed using the linear combination of univariate NURBS functions and the control points coordinates as
\begin{equation}\label{eq:NURBS_curve}
\mathbf{C}(\xi^i) = \sum_{l=1}^{n_i} R_{l}^{p_i}(\xi^i) \mathbf{X}_{l}^i\,.
\end{equation}
A bivariate continuum patch is constructed by the tensor product of the two univariate NURBS curves defined in the $ \xi^1 $ and $ \xi^2 $ parametric directions as
\begin{equation}\label{eq:NURBS_surface}
\mathbf{S}(\xi^1,\xi^2) = \sum_{l=1}^{n_1} \sum_{m=1}^{n_2} R_{lm}^{p_1,p_2} (\xi^1,\xi^2)\mathbf{X}_{lm}\,,
\end{equation}
where $ \{\mathbf{X}_{lm}\}_{l,m=1}^{n_1,n_2} $ is the coordinate vector of the control points net $ \mathbf{X} $. Moreover, $ R_{lm} (\xi^1,\xi^2) > 0 $ is a bivariate NURBS function that is given by the tensor product of two univariate B-spline functions defined along the $ \xi^1 $ and $ \xi^2 $ parametric directions as
\begin{equation}\label{eq:bivariate_NURBS_basis}
R_{lm}^{p_1,p_2} (\xi^1,\xi^2) =  \frac{w_{lm}}{W(\xi^1,\xi^2) } N_{l,p_1}(\xi^1)N_{m,p_2}(\xi^2)\,,
\end{equation}
where $W(\xi^1,\xi^2)= {\sum_{l=1}^{n_1}\sum_{m=1}^{n_2} w_{lm} N_{l,p_1}(\xi^1)N_{m,p_2}(\xi^2)} $ is a normalized weight function that is defined in terms of weight point $ w_{lm} $ and the B-spline functions. A composition of knot vector with associated control points accompanied by weights constitutes a patch. In the existing NURBS based discretization procedure, the parametrization for the contact surface is directly inherited from the bulk parametrization, but only in any one of the two parametric directions.

\subsection{Refinement strategies}	
In the context of IGA, a NURBS described geometry can be refined by means of knot insertion ($ h- $refinement), order-elevation ($ p- $refinement), and $ k- $refinement based strategies~\cite{hughes2005, Cottrell2009}. In the knot insertion strategy, an additional knot in the knot vector is inserted. If the inserted knot value is unique, an additional knot-span, consequently an element, without changing the inter-element continuity of the NURBS, is introduced. On the other hand, the repetition of knots reduces the smoothness of the NURBS functions across the element boundary.

In the order elevation based strategy, the geometry is refined by means of raising the order of the NURBS interpolations. For the elevation of order from $ p_i $ to $ p_i + t $, each knot value is repeated by $ t $ times. As a result, the continuity of the NURBS remains unchanged during the order-elevation ($ p- $refinement). 

The particular order of application of order elevation and knot insertion strategies yields $ k- $refinement~\cite{hughes2005}. During $ k- $refinement, first, the inter-element continuity of the NURBS interpolations is increased, and then, the additional elements within the given knot vector are introduced. On the other hand, if their application order is reversed, i.e. the knot insertion is performed first before the order elevation strategy, the inter-element continuity of the NURBS remains unchanged and knot values are repeated due to order-elevation. This, as a result, yields a large number of control points as compared to $ k- $refinement based strategy.

Within the existing NURBS-based discretization approach, the interpolation order of the NURBS constructed geometry is uniformly elevated to refine its contact boundary layer. However, this approach is not computationally favorable since the higher-order NURBS used for the evaluation of the contact integrals additionally have to be employed for the computation of the large majority of the region that is away from the contact interface. It, therefore, becomes desirable to develop an improved NURBS-based discretization method that provides a possibility to perform controllable order-elevation for a NURBS described geometry, as also suggested by Temizer et al.~\cite{Temizer2012}. This method is introduced in this paper and is described in detail in Section~\ref{sec:global_strategy}.

\section{Frictional contact formulation} \label{sec:Formulation} 	
In this section, we present the mathematical formulation for two-dimensional large deformation contact between two deformable bodies, considering Coulomb friction. The presented formulation is based on~\cite{laursen2003, wriggers2006, unbiased2013, unbiased2015}.

\subsection{Problem description and weak form} 	
It is assumed that two elastic bodies come into contact and undergo large deformations, see Fig.~\ref{fig:contact_system} for schematic illustration. According to the convention in~\cite{wriggers2006}, one of them is denoted as the slave body $ \mathcal{B}^\textrm{s} $ and another as master body $ \mathcal{B}^\textrm{m} $ occupying a bounded domain $ \Omega^{k} $ in $ \mathbb{R}^2 $. In the following, the superscript $ k = \{\textrm{s}, \textrm{m}\} $ is used to specify the slave and master bodies, respectively. The boundary $ {\Gamma}^{k} $ of a body $ \mathcal{B}^{k} $ can consists following three distinct parts: $ \Gamma^{k}_{\sigma} $ where boundary tractions are prescribed, $ \Gamma^{k}_{u} $ where displacements are prescribed, and $ \Gamma^{k}_{\textrm{c}} $ as the interface where contact will occur. It is assumed that these boundaries satisfy: $ \Gamma^{k}_{\sigma} \, \cap \,  \Gamma^{k}_{u} \, \cap \Gamma^{k}_{\textrm{c}} = \emptyset $.

The deformed (or current) configuration of each body $ \mathcal{B}^{k} $ is given by in terms of their generic point coordinates $ \bm{x}^k $ that are expressed as
$ \bm{x}^k = \bm{X}^k + \bm{u}^k\,, $	where $ \bm{X}^k $ denotes the coordinates of a generic points $ \bm{x}^k $ in the initial (or reference) configuration and $ \bm{u}^k $ is the displacement field of the same points of body $ \mathcal{B}^{k}$. The master contact surface $ \Gamma_\mathrm{c}^\mathrm{m} $ is parametrized using the convective coordinate $ {\xi}^\textrm{m} $ and  the covariant tangent vector $  {{\bm{\uptau}}}_1^{} $ at $ {\xi}^\textrm{m} $ is defined as
\begin{equation}\label{eq:covariant_tangent_vector}
{{\bm{\uptau}}}_1^{}(\xi^{\mathrm{m}}) = \frac{\partial \bm{x}^\textrm{m}}{\partial{\xi^\textrm{m}}} = \bm{x}^\textrm{m}_{\,, \,{\xi^\textrm{m}}}\,.
\end{equation}
The contravariant vector follows from $ {\bm{{{\uptau}}}}^1 :=  m^{11}  \bm{\uptau}_1^{} $, where the inverse metric is given by $ m^{11} :=  1/( \bm{\uptau}_1^{} \cdot  \bm{\uptau}_1^{})$.

We write the variational form for a system involving contact using the principal of the stationary total potential energy as~\cite{wriggers2006}
\begin{equation}\label{eq:total_potential}
\delta \, (\Pi_{\textrm{int}} + \Pi_{\textrm{c}} - \Pi_{\textrm{ext}} ) = 0 ~ \Leftrightarrow ~ \delta \Pi_{\textrm{int}} + \delta\Pi_{\textrm{c}} - \delta\Pi_{\textrm{ext}} = 0\,,~~~~ \forall\, \delta\boldsymbol{u} \in \mathcal{V}\,. 
\end{equation}
where $ \Pi_{\textrm{int}} $, $ \Pi_{\textrm{c}} $ and $ \Pi_{\textrm{ext}} $ are the total internal potential, contact potential, and potential stemming from the external loads, respectively. The function $ \mathcal{V} $ denotes the space of kinematically admissible variation of the displacement field $ \delta\boldsymbol{u} $. For the examples analyzed in this work $ \Pi_{\textrm{ext}}  = 0 $ holds. According to~\cite{wriggers2006, unbiased2015}, the total internal potential and contact potential are defined as 
\begin{eqnarray}
\delta \Pi_{\textrm{int}}  &=&  \sum_{k}^{\textrm{s},\textrm{m}} \left[ \int_{\Omega^k} \textrm{grad}(\delta  \boldsymbol{u}^k) : \boldsymbol{\sigma}^k~\textrm{d}{v} \right] , ~~~\textrm{and } 
\label{eq:int_potential} \\	
\delta \Pi_{\textrm{c}} &=& - \int_{\Gamma^\textrm{s}_{\textrm{c}}} \boldsymbol{t}^\textrm{s} \cdot \delta{\boldsymbol{u}}^\textrm{s}  ~\textrm{d}\Gamma \, + \, \int_{\Gamma^\textrm{s}_{\textrm{c}}} \boldsymbol{t}^\textrm{s} \cdot \delta{\boldsymbol{u}}^\textrm{m}  ~\textrm{d}\Gamma\,, \label{eq:contact_potential} 
\end{eqnarray}
where $ \boldsymbol{\sigma}^k $ denotes the Cauchy stress tensor at a point $ \bm{x}^k \in \Omega^k $ and $ \bm{t}^s $ is the contact traction at a point on the contact interface of the slave body.

In the following sections, we present the contact search procedure for the determination of contact surface and the kinematical relations for the normal and frictional contact that are necessary for the evaluation of contact potential $ \delta \Pi_{\textrm{c}} $.

\subsection{Normal contact} \label{sec:normal_contact}
It is known that the contact surface in which the one body comes into contact to another is not known \textit{a priori}. Thus, for the determination of the contact surface in the current configuration, a distance function $ \hat{d}({\xi}^\textrm{m}) := || \bm{x}^{\textrm{s}} - {\bm{x}}^\textrm{m}(\xi^{\textrm{m}})|| $
is introduced. Here, $ {\bm{x}}^\textrm{m}(\xi^{\textrm{m}}) $ denotes the physical coordinates of an arbitrary point computed at the parametric point $ \xi^{\mathrm{m}} $ on $ \Gamma_{\textrm{c}}^{\mathrm{m}} $ in the current configuration. With this definition, the distance between a given (fixed) slave point $ \bm{x}^\textrm{s} \in \Gamma^\textrm{s}_\textrm{c} $ on the slave contact surface and an arbitrary point {$ {\bm{x}}^\textrm{m}(\xi^{\textrm{m}}) \in \Gamma^\textrm{m}_\text{c} $} on the master surface is calculated. For the determination of closest contact point $ \bar{x}^{\mathrm{m}} $ onto the master surface corresponding to a slave point, the minimum of distance function which can be expressed as~\cite{wriggers2006}
\begin{equation}\label{eq:projection_point}
\frac{d}{d\xi^{\textrm{m}}} \, \hat{d}({\xi}^\textrm{m}) = [\bm{x}^{\textrm{s}} - {\bm{x}}^\textrm{m}(\xi^{\textrm{m}})] \cdot  \bm{\uptau}_1^{}(\xi^{\textrm{m}}) \,.
\end{equation}	
The above derivative vanishes if the line describing the distance between the slave point $ \bm{x}^{\textrm{s}} $ and master point $ {\bm{x}}^\textrm{m}(\xi^{\textrm{m}}) $ is orthogonal to the direction of the covariant tangent vector $ \boldsymbol{\uptau}_1^{} $ computed at $ {\xi}^{\textrm{m}} $, see Eq.~(\ref{eq:covariant_tangent_vector}). The corresponding parametric point $ \xi^{\textrm{m}} $ at which $ \frac{d}{d\xi^{\textrm{m}}} \, \hat{d}({\bar{\xi}}^\textrm{m}) = 0 $ is referred as closest projection point $ \bar{\xi}^{\textrm{m}} $, which is denoted by a bar over \mbox{$ \xi^{\mathrm{m}} $} in the subsequent description. Since the  Eq.~(\ref{eq:projection_point}) is nonlinear, an iterative solution method is required to compute the closest  point $ \bar{\xi}^{\textrm{m}} $. For this, the local Newton-Raphson method, starting from an initial guess ${\xi}^{\textrm{m}} = \tilde{\xi}^{\textrm{m}} $ is used. It is computed by taking the second derivative of the distance function is given by
\begin{equation}\label{eq:initial_guess}
\frac{d^2}{d(\xi^{\textrm{m}})^2} \, \hat{d}({\xi}^\textrm{m}) = [\bm{x}^{\textrm{s}} - {\bm{x}}^\textrm{m}(\xi^{\textrm{m}})] \cdot  \bm{x}^{\textrm{m}}_{,\,\xi^{\textrm{m}}\xi^{\textrm{m}}}(\xi^{\textrm{m}})  - \boldsymbol{\uptau}_1^{} \cdot \boldsymbol{\uptau}_1^{} \,.
\end{equation}
\begin{figure}[!t]
	\centering
	{\includegraphics[width=0.8\linewidth]{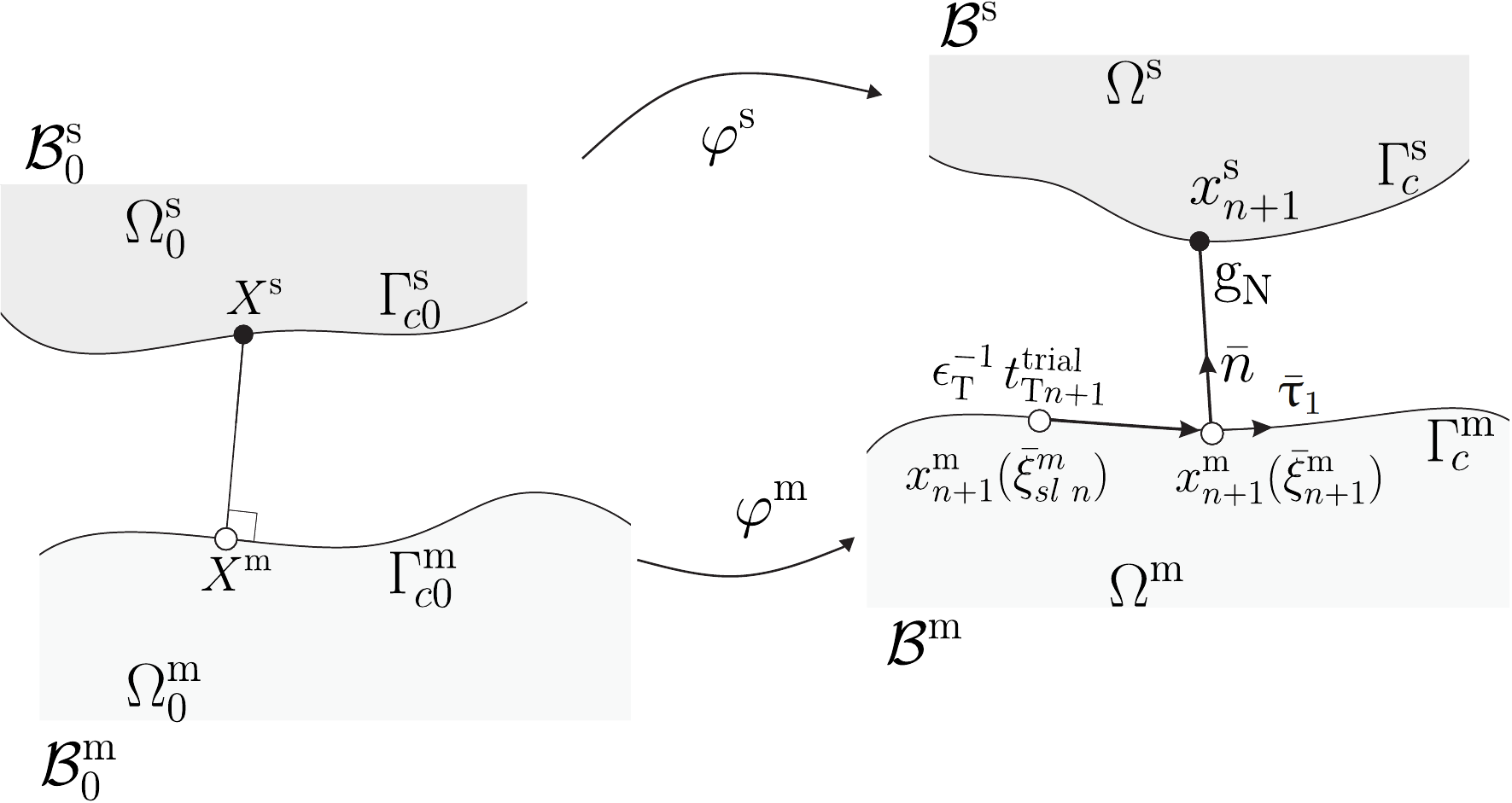}}
	\caption{A schematic illustration of the kinematic contact quantities for a given point $ \bm{x}^{\textrm{s}}_{n+1} $ on contact surface $ \Gamma_{\textrm{c}}^{\textrm{s}} $ and a projection point $ \bm{x}^{\textrm{m}}_{n+1} (\bar{\xi}^{\textrm{m}}_{n+1})$ on $ \Gamma_{\textrm{c}}^{\textrm{m}} $ in the current configuration of a system involving frictional contact.}
	\label{fig:contact_system}
\end{figure}	
The physical coordinates of the closest point $ \bar{\xi}^{\textrm{m}} $ on $ \Gamma^\textrm{m}_\text{c} $ corresponding to the slave point $ \bm{x}^\textrm{s} $ is denoted by $ \bar{\bm{x}}^\textrm{m} = {\bm{x}}^\textrm{m}(\bar{\xi}^{\textrm{m}}) $. Further, the covariant tangent vector computed at projection point $ \bar{\xi}^{\textrm{m}} $ using Eq.~(\ref{eq:covariant_tangent_vector}) is denoted by $ \bar{\bm{\uptau}}_1 $.

Once the physical projection point $ \bar{\bm{x}}^\textrm{m} $ corresponding to a slave point $ \bar{\bm{x}}^\textrm{s} $ is evaluated, the normal gap $ \textrm{g}^{}_{\textrm{N}} $ between these points is determined by
\begin{equation}\label{eq:gN}
\textrm{g}^{}_{\textrm{N}}  = (\bm{x}^\textrm{s} - \bar{\bm{x}}^\textrm{m})\cdot \bar{\bm{n}}\,,
\end{equation} 
where $ \bar{\bm{n}} $ is the outward unit normal vector at $ \bar{\bm{x}}^\textrm{m} $, see Fig.~\ref{fig:contact_system}. It is used to distinguish whether the contact between the two bodies is active or non-active. In case of penetration $ (\textrm{g}^{}_{\textrm{N}}  < 0) $, the normal contact tractions $ \text{t}_\text{N} = \boldsymbol{t} \cdot \bar{\mathbf{n}} $ are activated that avoid the penetration of the contacting regions. In order to enforce the impenetrability in the normal direction, following Karush-Kuhn-Tucker (KKT) conditions are considered
\begin{equation}\label{eq:KKT1}
\textrm{g}^{}_{\textrm{N}}  \geq 0\,,~~~ \text{t}_{\text{N}} \leq 0\,,~~~\textrm{g}^{}_{\textrm{N}} \, \mathrm{t}^{}_{\textrm{N}} = 0\,, ~~~~~~\textrm{on}~ \Gamma_\mathrm{c}^s\,.
\end{equation}

The above-described contact constraints cannot be directly incorporated into the contact variational form as they lead to the non-smooth relationship between the normal gap and the contact pressure, resulting in a non-smooth normal contact constitutive law. In the present work, for the regularization of the contact constraints defined in Eq.~(\ref{eq:KKT1}) the penalty method is adopted. According to~\cite{wriggers2006, unbiased2015}, the penalty regularized normal contact constraint is defined as
\begin{equation}\label{eq:normal_contact}
\bm{t}_{\text{N}} = 
\left\{\begin{array}{lr}
-\epsilon_{\text{N}} \, \textrm{g}_{\textrm{N}}{\boldsymbol{\bar{n}}},  & \textrm{if}~ \textrm{g}_{\textrm{N}} < 0\\
\bm{0},  & \textrm{otherwise}\,,
\end{array}	\right.
\end{equation}
where $ \epsilon_{\text{N}} > 0 $ denotes the normal penalty parameter.

\subsection{Frictional contact} \label{sec:frictional_contact}
For frictional contact, the contact traction vector $ \boldsymbol{t} $ is described in terms of the normal and tangential components as~\cite{unbiased2015}
\begin{equation}\label{eq:contact_traction}
\boldsymbol{t} = \textrm{t}_{\text{N}}\bar{\bm{n}} - \textrm{t}_{\text{T}}\bar{\bm{\uptau}}_1
\end{equation}	
where $ \textrm{t}_{\textrm{N}} $ is defined in Eq.~(\ref{eq:normal_contact}), and  $ \textrm{t}_{\text{T}} = \bm{t}\cdot \bar{\bm{\uptau}}_1 $ is the tangential traction. Also, $ \bar{\bm{n}} $ and $ \bar{\bm{\uptau}}_1 $ are the unit normal and tangent vectors computed at $ \bar{\xi}^{\textrm{m}} $. The tangential traction $ \bm{t}_{\textrm{T}} = \textrm{t}_{\textrm{T}}\bar{\bm{\uptau}}_1 $ is determined based on the stick-slip behavior that occurs during frictional contact. The distinction between the stick and slip states is drawn by the slip function $ \Phi_{} = || \bm{t}_{\text{T}} || -\mu_f \,	\text{t}_{\text{N}} \,, $	
where $ \mu_f $ denotes the friction coefficient. It states sticking for $ \Phi_{} < 0 $ and sliding if $ \Phi_{} = 0 $. The function $ \Phi_{} = 0$ corresponds to a surface in the traction space $ \{ \text{t}_{\text{N}}, \text{t}_{\text{T}} \} $. The tangential traction during the frictional sliding is defined, e.g. for Coulomb's law, as
\begin{equation}\label{eq:friction_law}
\bm{t}_{\textrm{T}} = -\mu_f \, \textrm{t}_{\textrm{N}} \frac{\dot{
		\textrm{g}}_{\textrm{T}}}{|| \dot{
		\textrm{g}}_{\textrm{T}} ||}
\end{equation}
where $ \mu_f $ denotes the Coulomb friction coefficient and $ \dot{\textrm{g}}_{\textrm{T}} $ is the relative tangential sliding velocity.

In the context of computational contact mechanics, the constitutive equations for friction are often formulated in the framework of elasto-plasticity~\cite{wriggers2006}. According to this, the total tangential slip $ \mathbf{g}_{\textrm{T}} $ can be decomposed into an elastic (or stick) part $ \mathbf{g}_{\textrm{T}e} $ and a plastic (or slip) part $ \mathbf{g}_{\textrm{T}sl} $ as
\begin{equation}\label{key}
\mathbf{g}_{\textrm{T}}  = 		\mathbf{g}_{\textrm{T}e} + 		\mathbf{g}_{\textrm{T} sl}\,.
\end{equation}
The tangential stick $ \mathbf{g}_{\textrm{T}e} $ is characterized by no relative tangential motion between the surfaces during the frictional contact. Such a physical behavior enforces a geometrical constraint in the relative motion of the contact surface. For its regularization the penalty method is used herein. The tangential slip $ \mathbf{g}_{\textrm{T}sl} $ is characterized by the sliding of a slave point $ \bm{x}^\textrm{s} $ relative to the master surface $ \Gamma_{\textrm{c}}^{\textrm{m}} $. It is a dissipative process. The constitutive evolution equation for tangential slip $ \mathbf{g}_{\textrm{T}sl} $ can be derived by considering the dissipation during friction. We write the evolution law for the tangential slip as~\cite{wriggers2006}	
\begin{equation}\label{eq:evolution_law}
\dot{\mathbf{g}}_{\textrm{T}sl} = \gamma \frac{\bm{t}_{\text{T}}}{||\bm{t}_{\text{T}}||}\,.
\end{equation}
where $ \gamma $ is the slip parameter that is determined using the following conditions
\begin{equation}\label{eq:KKT_friction}
\gamma \geq 0\,,~~~~\Phi_{} \leq 0\,, ~~~{\gamma}\, \Phi_{} = 0\,.
\end{equation}
The above set of equations is known as KKT conditions for Coulomb friction~\cite{wriggers2006}.

The tangential traction is updated based on a classical corrector-predictor approach. In this, first, the (penalty regularized) elastic trial step is computed, followed by the trial slip function, as in~\cite{unbiased2015}, with
\begin{eqnarray}\label{eq:trial_step}
\bm{t}_{\text{T}n+1}^{\textrm{trial}} &=& \epsilon_{\textrm{T}} \left(\, \boldsymbol{x}^{\textrm{m}}_{n+1}(\bar{\xi}^{\textrm{m}}_{n+1}) \,  - \, \boldsymbol{x}^{\textrm{m}}_{n+1}({\xi}^{\textrm{m}}_{{sl}\,n})\, \right), \\
\Phi_{n+1}^{\textrm{trial}} &=& || \bm{t}_{\text{T}n+1}^{\textrm{trial}} || -\mu_f \,	\text{t}_{\text{N}n+1}\,.
\end{eqnarray}
Here, all the quantities defined in the current configuration at time step $ t = t_{n+1} $ are denoted with subscript $ n+1 $. Similarly, all the quantities at previous time step $ t = t_n $ with subscript $ n $. The tangential penalty parameter is denoted by $ \epsilon_{\textrm{T}}> 0 $. Moreover, $ \boldsymbol{x}^{\textrm{m}}_{n+1}(\bar{\xi}^{\textrm{m}}_{n+1}) $ and $ \boldsymbol{x}^{\textrm{m}}_{n+1}({\xi}^{\textrm{m}}_{{sl}\,n}) $ are the physical coordinate corresponding to the closest projection point $ \bar{\xi}^{\textrm{m}}_{n+1} $ and tangential slip point $ {\xi}^{\textrm{m}}_{{sl}\,n} $, respectively, see Fig.~\ref{fig:contact_system}. According to~\cite{wriggers2006, unbiased2015}, the tangential traction is determined by based on the changes in stick-slip status with
\begin{equation}\label{eq:tangential_traction}
\bm{t}_{\textrm{T}} = 
\left\{ \begin{array}{lr}
\bm{t}_{\text{T}n+1}^{\textrm{trial}}, ~~ &  \textrm{if}~ \Phi_{n+1}^{\textrm{trial}} \leq 0\,,\\
\mu_f\, \textrm{t}_{\textrm{N}n+1} \bm{n}_{\textrm{T}n+1} ,~~ & \textrm{otherwise}\,. 
\end{array}	\right.
\end{equation}
where $ \bm{n}_{\textrm{T}n+1} = {	\mathbf{t}_{\text{T}n+1}^{\textrm{trial}}}/{||	\mathbf{t}_{\text{T}n+1}^{\textrm{trial}} ||} $.

\subsection{NURBS discretized weak form} \label{sec:NURBS_weak_form}
Having discretized the body into a set of elements $ \Omega^k = \cup_{e} \Omega^{ke} $, the same NURBS functions, as employed in Eq.~(\ref{eq:NURBS_surface}), are used for the approximation of the unknown displacement field, its variation and current coordinates within each element as
\begin{equation}\label{eq:disp_field}
\boldsymbol{u}^{ke} = \mathbf{R}^k(\boldsymbol{\xi}^k)\,\mathbf{u}^{ke}\,,~~~\delta \boldsymbol{u}^{ke} = \mathbf{R}^k(\boldsymbol{\xi}^k)\,\delta\mathbf{u}^{ke}\,,~~~\bm{x}^{ke} = \mathbf{R}^k(\boldsymbol{\xi}^k)\,\boldsymbol{x}^{ke} ~~~~~ \forall\, \boldsymbol{\xi}^k \in \Omega^{ke}\,.
\end{equation}
Here, $ \mathbf{R}^k = [ R_1^k \boldsymbol{I}, R_2^k \boldsymbol{I}, \, \,.\,.\,.\, ,R_{n_{cp}^{ke}}^k \boldsymbol{I} ] $ is an array of size $ 2\times 2 n_{cp}^{ke} $ with identity matrix $ \boldsymbol{I} $ in $ {\rm I\!R}^2 $. The $ n_{cp}^{ke} = (p_1^k+1) \times (p_2^k+1) $ denotes the total number of control points whose bivariate NURBS basis functions $ (R^{k}_a(\boldsymbol{\xi}^k)) $ have local support in an element $ \Omega^{ke} $. Also, $ \mathbf{u}^{ke} $, $ \delta \mathbf{u}^{ke} $ and $ \bm{x}^{ke} $ denotes the arrays containing the displacement field, its variation and the current coordinates of $ n_{cp}^{ke} $ number of control points in an element, respectively. The modification in the approximation of $ \boldsymbol{u}^{ke} $, $ \delta \boldsymbol{u}^{ke} $ and $ \bm{x}^{ke} $ for VO based NURBS discretization is explained in Sec.~\ref{sec:global_strategy}.

Using the definition of the displacement field given in Eq.~(\ref{eq:disp_field}), the NURBS discretized, Eq.~(\ref{eq:total_potential}), takes the form
\begin{equation}\label{eq:weak_form}
\delta \mathbf{u}^{\textrm{T}} [ \mathbf{f}_{\textrm{int}} +  \mathbf{f}_{\textrm{c}}] = \mathbf{0}\,,~~~~~\forall \, \delta \mathbf{u} \in \mathcal{V}
\end{equation}
where $ \mathbf{f}_{\textrm{int}} $ and $ \mathbf{f}_{\textrm{c}} $ are the global control point force vectors corresponding to the $ \Pi_{\textrm{int}} $, defined in Eq.~(\ref{eq:int_potential}), and $ \Pi_{\textrm{c}} $, in Eq.~(\ref{eq:contact_potential}), respectively.

Since, the variation of displacement $ \delta \mathbf{u} $ is arbitrary, the Eq.~(\ref{eq:weak_form}) can be written as $ \mathbf{f}_{\textrm{int}} +  \mathbf{f}_{\textrm{c}} = \mathbf{0} $. The global forces $ \mathbf{f}_{\textrm{int}} $ and $ \mathbf{f}_{\textrm{c}} $ are obtained based on the assembly of their elemental contributions $\mathbf{f}_{\textrm{int}}^{ke} $ and $ \mathbf{f}_{\textrm{c}}^{ke} $ using the control point connectivity arrays for each body $ \mathcal{B}^k $, see~\cite{Cottrell2009, Agrawal2018}. 

The elemental internal force vector for a volume element $ \Omega^{ke} $ is given by
\begin{equation}\label{eq:element_int_force_vector}
\mathbf{f}_{\textrm{int}}^{ke}  = \int_{\Omega^{ke}} [\mathbf{B}^{k}]^\textrm{T} \, \boldsymbol{\sigma}^{k}~\textrm{d}v\,,
\end{equation}
where $ \mathbf{B}^{k} $ array contains the $ n_{cp}^{ke} $ number of the derivatives of bivariate NURBS functions having local support in an element $ \Omega^{ke} $. In the examples solved herein, the Neo-Hookean material model is used for the description of bodies. The corresponding constitutive relation is given by~\cite{bonet1997}
\begin{equation}\label{key:Neo-Hookean}
{\bm{\sigma}} = \frac{\lambda}{J}(\textrm{ln} J)\boldsymbol{I} + \frac{\mu}{J} (\boldsymbol{FF}^{{\textrm{T}}} - \boldsymbol{I})\,,
\end{equation}
where $ \boldsymbol{F} $ is the deformation gradient,  $ J $ is the determinant of $ \boldsymbol{F} $, and $ \boldsymbol{I} $ is the identity tensor. The Lam\'{e} constants, i.e. the bulk modulus and shear modulus, are defined via $ \lambda = 2\mu\nu / (1-2\nu) $ and $ \mu = E/ 2(1+\nu) $, respectively.

\subsection{Contact algorithm} \label{sec:contact_formulation}
In this work, we adopt the classical full-pass version of the contact algorithm introduced by Sauer and De Lorenzis~\cite{unbiased2015} to integrating the proposed VO based NURBS discretization methodology with the isogeometric contact analysis. The full-pass contact algorithm was first introduced by Fischer and Wriggers~\cite{Fischer2005} for the Lagrange-polynomial based FE analysis of the contact problem. It was later extended by Temizer et al.~\cite{Temizer2011}, De Lorenzis et al.~\cite{DeLorenzis2011, DeLorenzis2012} and Dimitri et al.~\cite{Dimitri2014Tspline, Dimitri2014, Dimitri2017} for the NURBS based isogeometric contact analysis. Recently, Dimitri et al.~\cite{Dimitri2014Tspline, Dimitri2014, Dimitri2017} further extended this algorithm for T-spline based isometric analysis of contact, cohesive zone modeling, and mixed-mode debonding problems. In these literature, a number of terminologies have been adopted to denote the classical full-pass approach such as  ``mortar" in~\cite{Fischer2005}, ``Knot-to-Surface" in~\cite{Temizer2011}, ``non-mortar" in~\cite{DeLorenzis2011,DeLorenzis2012}, and ``Gauss-point-to-surface (GPTS)" in~\cite{Dimitri2014Tspline, Dimitri2014, Dimitri2017}. In the present work, as in~\cite{Dimitri2014Tspline, Dimitri2014, Dimitri2017}, we also denote this approach as GPTS, since within this the contact constraints, i.e. impenetrability in Eq.~({\ref{eq:KKT1}}) and sticking in Eq.~(\ref{eq:KKT_friction}), are enforced independently at each quadrature point that enables the evaluation of the contact integrals. This formulation passes the contact patch test within the quadrature error, which ensures the convergence of the solution upon mesh refinement~\cite{Fischer2005, Zavarise2009}. The only drawback is that the GPTS based contact algorithm leads to over-constrained formulation, as demonstrated by~\cite{Temizer2011, DeLorenzis2011, DeLorenzis2012}. The over-constraining stems from directly enforcing the contact constraints at too many numbers of locations that, as a result, leads to numerical instability, see~\cite{Kikuchi1988} for the contact cases. Further, with GPTS approach, oscillatory responses of contact forces are obtained if a large value of penalty parameter is employed. However, in our examples, we use the value of penalty parameter, for which the loss of the numerical stability issues becomes significant with the GPTS approach, lay beyond those required to obtain a satisfactory quality solution from the engineering perspective. The values lies within a range where the ill-conditioning of the global system matrix related issues still arises. The post-processing smoothening scheme of Sauer~\cite{Sauer2013} can also be used as in~\cite{Dimitri2014Tspline, Dimitri2014} to alleviate the oscillatory responses of contact forces effectively. The incorporation of the mathematically more sophisticated mortar-based isogeometric contact algorithm by Temizer et al.~\cite{Temizer2011, Temizer2012}, De Lorenzis et al.~\cite{DeLorenzis2011, DeLorenzis2012}, Seitz et al.~\cite{SEITZ2016}, Popp et al.~\cite{Popp2010, Popp2013, Popp2018} and Duong et al.~\cite{Duong2018} is subject of future research work.

The elemental contact force vector $ \mathbf{f}_c^{ke} $ on the contact surfaces $ \Gamma_c^{ke} $ of a contact element $ \Omega^{ke}_c $ for a body $ \mathcal{B}^k $ is given as
\begin{equation}\label{eq:contact_force_element}
\mathbf{f}_{\textrm{c}}^{\textrm{s}e}  = -\int_{\Gamma^{\textrm{s}e}_c} [\mathbf{N}^{\textrm{s}}(\xi^{\textrm{s}})]^\textrm{T}\, \boldsymbol{t}^\textrm{s}~\textrm{d}\Gamma\,,~~~~~~~ \mathbf{f}_{\textrm{c}}^{\textrm{m}e}  = \int_{\Gamma^{\textrm{s}e}_c} [\mathbf{N}^{\textrm{m}}(\bar{\xi}^{\textrm{m}})]^\textrm{T} \, \boldsymbol{t}^\textrm{s}~\textrm{d}\Gamma\,,
\end{equation}
where $ \mathbf{N}^{k} $ array consists $ n_{cps}^{ke} = p_1^k+1 $ number of univariate NURBS basis functions having support on the contact surface $ \Gamma_c^{ke} $. In this, all the quantities are computed within the current configuration of the bodies. The contact traction $ \boldsymbol{t}^\textrm{s} $ is computed at each quadrature point using  Eqs.~(\ref{eq:normal_contact}) and~(\ref{eq:tangential_traction}) for frictionless and frictional contact, respectively. An active set strategy that incorporates the contribution of the quadrature points in case of penetration is utilized for the evaluation of contact integral. The changes in the evaluation of $ \mathbf{f}_{\textrm{c}}^{ke} $ with VO based NURBS discretization is detailed in the next section.

The linearization of the internal force vector $ \mathbf{f}_{\textrm{int}}^{ke} $ and contact force vector $ \mathbf{f}_{\textrm{c}}^{ke} $ is necessary for the solution of resulting nonlinear system of equations using the Newton-Raphson method. The linearization of the contact force vector yields the contact tangent matrix, see~\cite{unbiased2013, unbiased2015}, and their corresponding expressions are provided in Appendix~\ref{appendix:A} for the sake of completeness.

\section{Varying-order based NURBS discretization} \label{sec:global_strategy}
In this section, we introduce the theory for the varying-order based NURBS discretization method, followed by the description on its integration into an existing isogeometric contact code.

\subsection{Varying-order NURBS} \label{subsec:proposed_technique}
\begin{figure}[!t]
	\centering
	\subfloat[]{\includegraphics[width=0.8\linewidth]{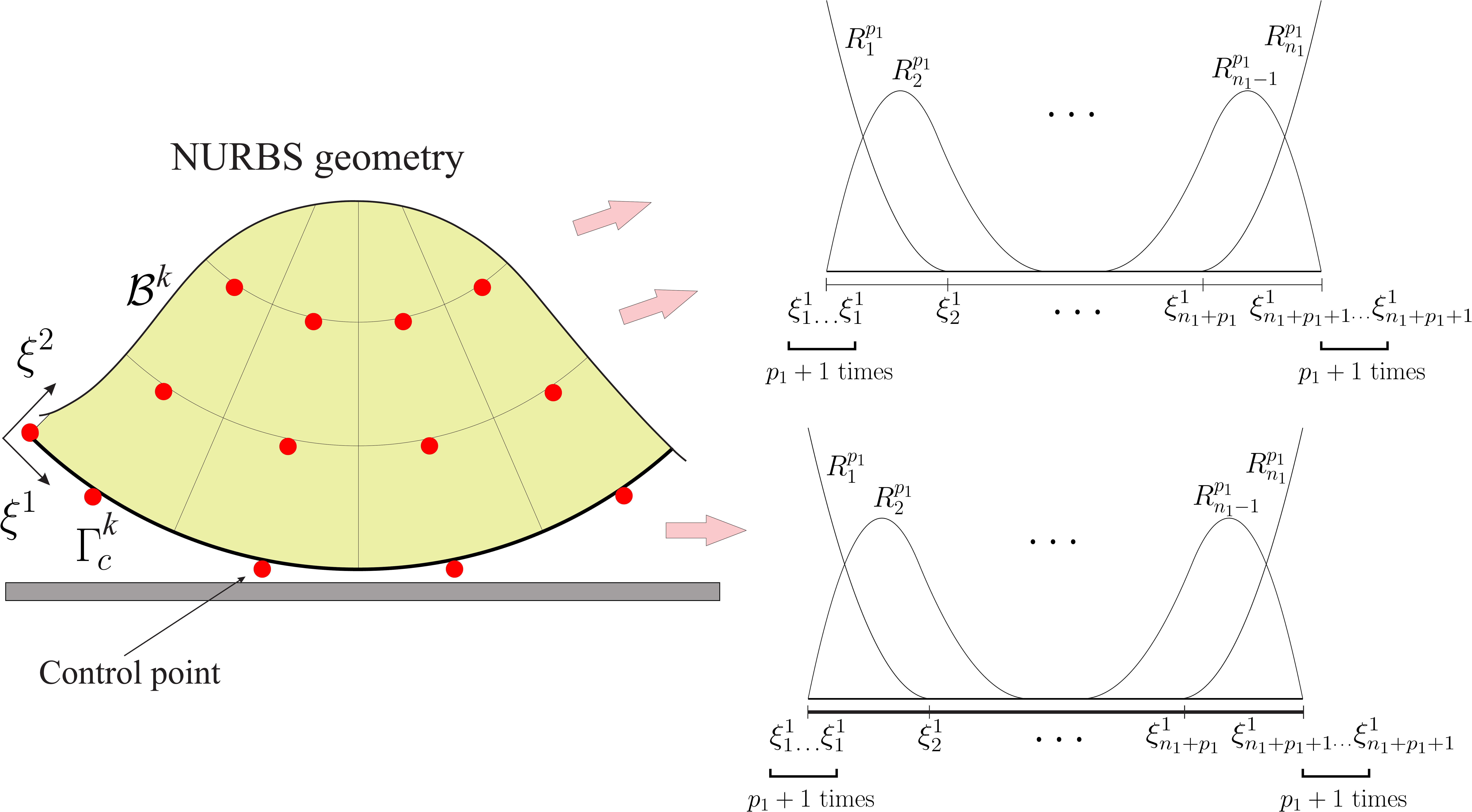}\label{fig:Enrichmenta}} \\
	\subfloat[]{\includegraphics[width=0.8\linewidth]{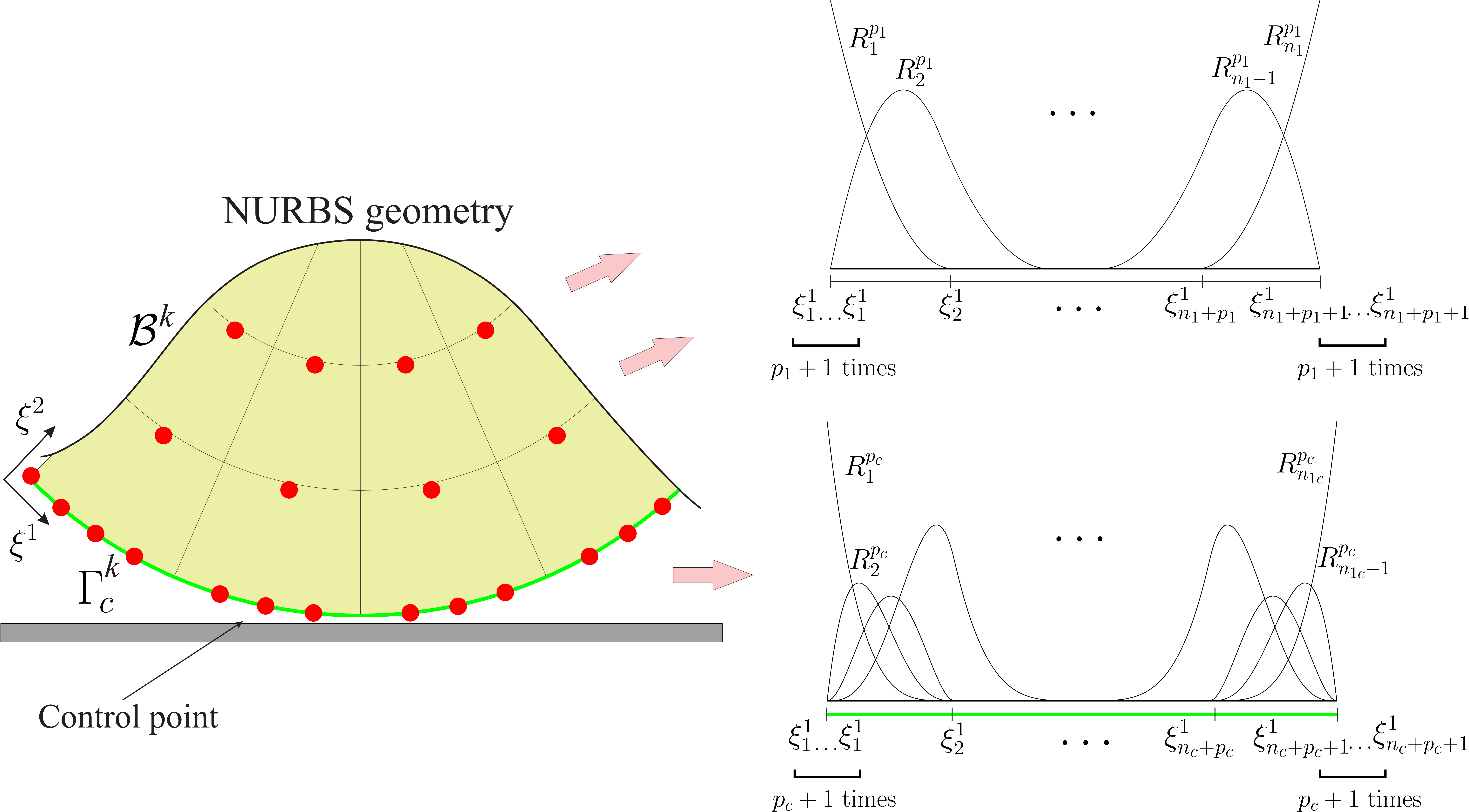}\label{fig:Enrichmentb}}
	\caption{A schematic illustration of the VO based NURBS discretization method for a given geometry. (a) Exact representation of the geometry with minimum $ p_1 $ and $ p_2 $ order of NURBS along the $ \xi^1 $ and $ \xi^2 $ parametric directions with a very coarse mesh. (b) Representation of the VO NURBS discretized geometry where higher-order NURBS (i.e. $ p_c>p_1 $) are used for the contact boundary layer, and minimum order NURBS interpolations are used for the remaining bulk domain. The accompanying control points are shown with red dots and the new contact boundary layer with a bold green line. The corresponding basis functions for the contact surface and bulk domain are also shown.}
	\label{fig:Enrichment}
\end{figure}

The basic concept of the proposed varying-order (VO) based NURBS discretization for isogeometric contact analysis is illustrated in Fig.~\ref{fig:Enrichment}. Consider a body $ \mathcal{B}^k $ having contact boundary layer $ \Gamma_\textrm{c}^k $. Let $ p_1 $ and $ p_2 $ be the minimum orders of NURBS functions that are capable of representing the given geometry in an exact manner. The coarse mesh for the geometry is given by the product of open knot vectors $ \Xi^1 ~\times ~\Xi^2 $ defined along the $ \xi^1 $ and $ \xi^2 $ parametric directions, as shown in  Fig.~\ref{fig:Enrichmenta}. Next, in order to make use of higher-order NURBS functions for contact computations, the originally $ p_1 $ discretized NURBS contact boundary $ \Gamma_\textrm{c}^k $ is replaced with a higher-order $ p_c > p_1 $ of NURBS layer as shown in Fig.~\ref{fig:Enrichmentb}. The resultant discretization is accordingly denoted by N$ _{p}-$N$ _{p_c} $, where N$ _{p} $ $ (p = ~\mathrm{max}(p_1,p_2)) $ is the order of NURBS used for the description of bulk domain and N$ _{p_c} $ $ (p_c>p) $ for the contact boundary layer.

\begin{figure}[!t]
	\centering
	\subfloat[]{\includegraphics[width=0.65\linewidth]{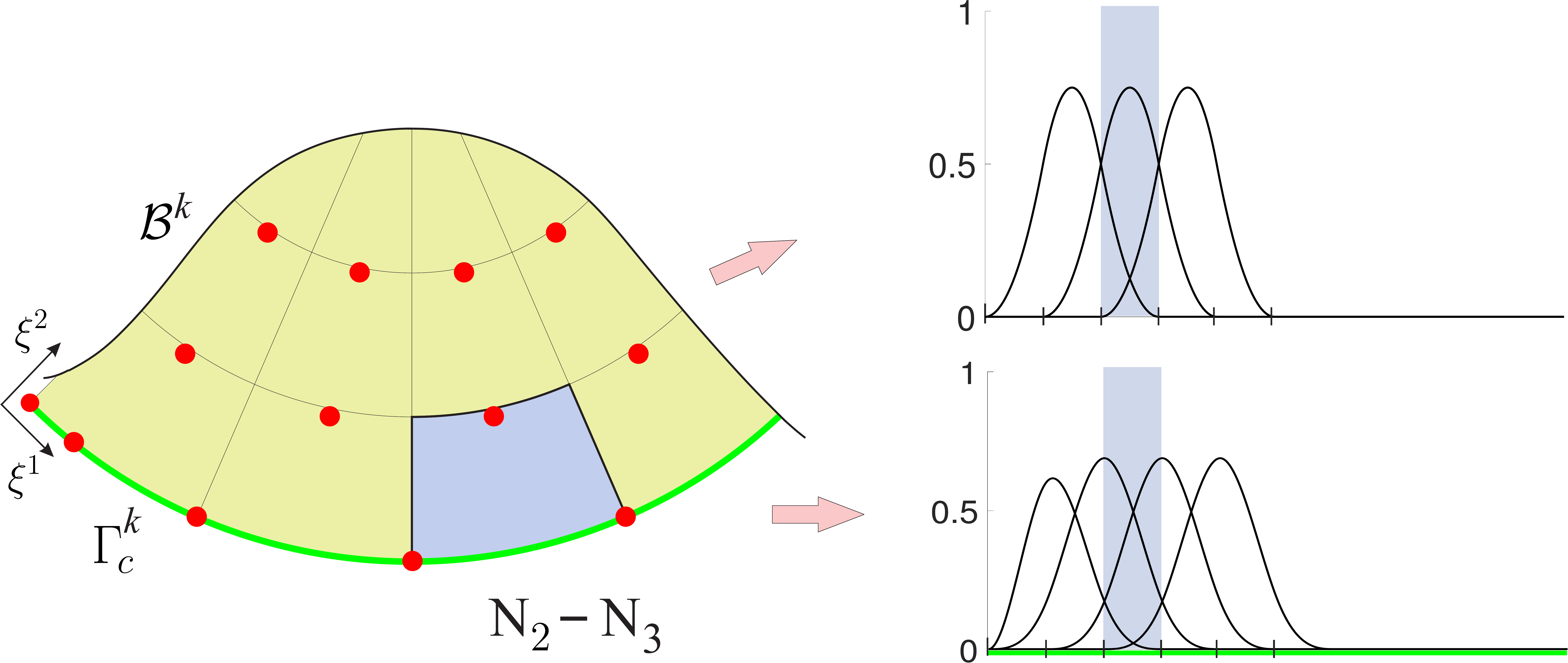}\label{fig:krefine_1}} \\	\subfloat[]{\includegraphics[width=0.65\linewidth]{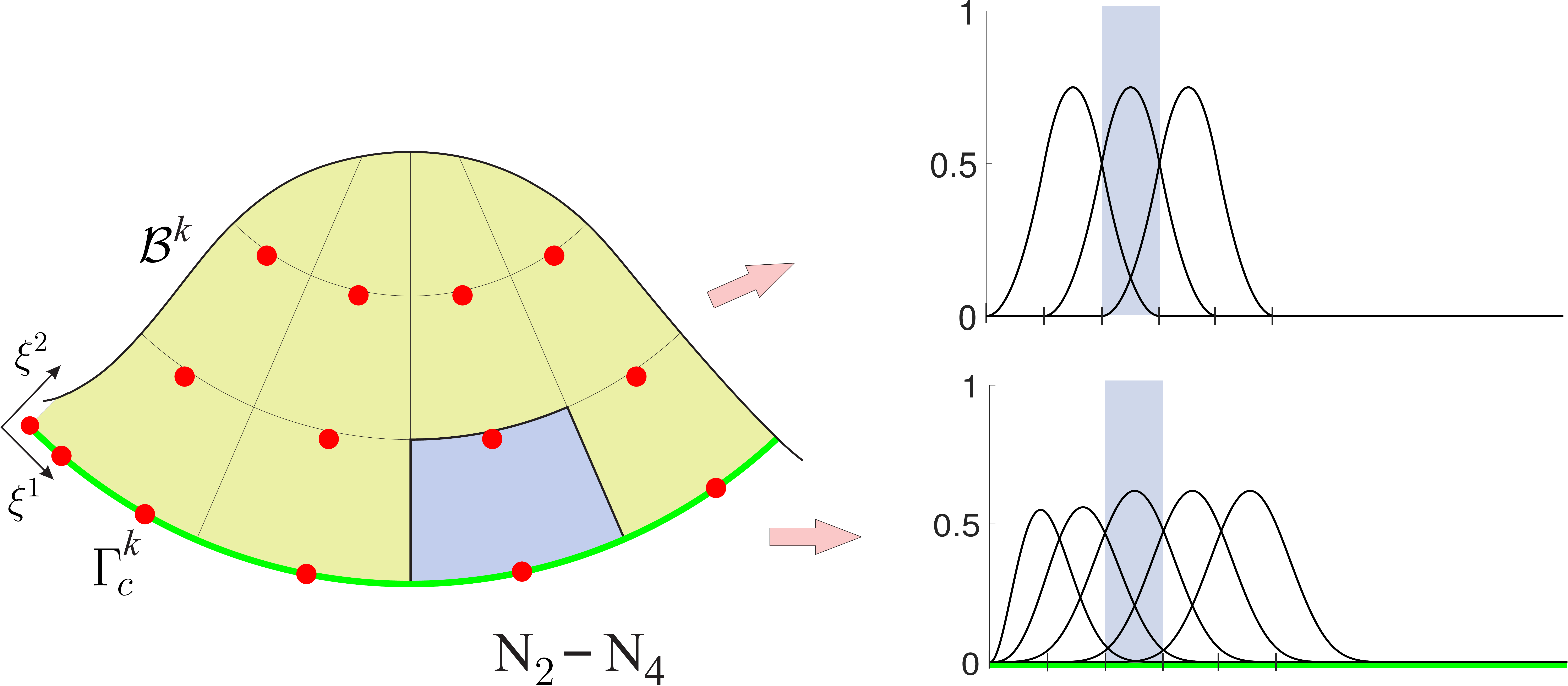}\label{fig:krefine_2}} 
	\caption{VO NURBS based N$ _2- $N$ _{p_c} $ $ (p_c = 3, ~\mathrm{and}~4) $ discretization of an original N$ _2 $ discretized geometry. In this arrangement, the $ k- $refined contact boundary layer is constructed using the (a) $ C^2-$ (with N$ _3 $), and (b) $ C^3-$continuous (with N$ _4 $) NURBS functions. The basis functions corresponding to the contact boundary and bulk part of a (highlighted) contact element for each N$ _2- $N$ _{p_c} $ arrangement are shown in the right column. }\label{fig:VO_krefine}
\end{figure}

\begin{figure}[!t]
	\centering
	\subfloat{\includegraphics[width=0.65\linewidth]{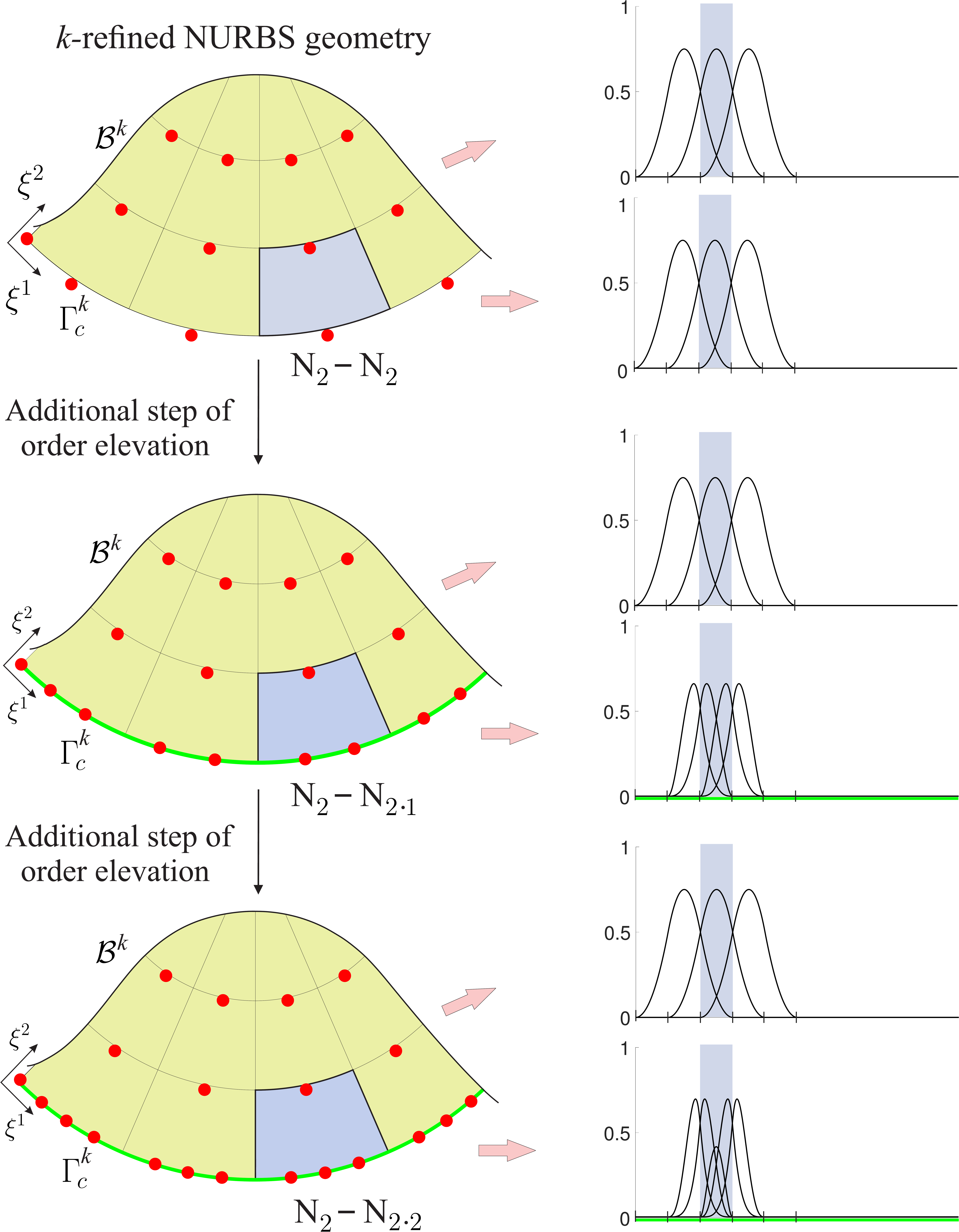}}
	\caption{An illustrative application procedure of the one (middle-row) and two (bottom-row) additional steps of order elevation strategy to the contact boundary layer of a $ k- $refined N$ _2- $N$ _2 $ NURBS geometry (first-row). The basis functions corresponding to a N$ _2- $N$ _{2\cdot p_s} $ $ (p_s = 1, ~\mathrm{and}~ 2) $ discretized contact element are shown in right.}\label{fig:VO_kp_refine}
\end{figure}

The new NURBS layer is constructed either using the $ k- $refinement or through a combination of $ k- $refinement and order-elevation strategies in such a manner that it matches the bulk parametrization, as shown in Fig.~\ref{fig:Enrichmentb}. The application of the $ k- $refinement strategy to the NURBS contact layer increases the order as well as the inter-element continuity of the NURBS functions, see Fig.~\ref{fig:VO_krefine}. On the other hand, application of one or two additional steps of order-elevation to $ k- $refined NURBS layer introduces a large number of additional control points along the contact boundary layer while the inter-element continuity, i.e. $ C^{p_c-1} $, remains unchanged, see Fig. \ref{fig:VO_kp_refine}. With this approach, the resultant VO based NURBS discretization is denoted by N$ _{p}- $N$ _{p_c\cdot p_s} $, where $ p_s $ is the step number of order elevation strategy that is additionally applied to N$ _{p_c} $ NURBS layer.

In Figs.~\ref{fig:VO_krefine} and \ref{fig:VO_kp_refine}, one VO NURBS discretized contact element is highlighted to illustrate the difference between the above two strategies. As shown in Fig.~\ref{fig:VO_krefine}, for employing the higher-continuous NURBS for contact computations, the contact boundary layer of an original N$ _{2} $ discretized example geometry is replaced with a user-defined $ C^{p_c -1} $ continuous N$ _{p_c} $ NURBS layer, where $ (p_c = 3, \mathrm{and}~4) $. For a N$ _{p_c} $ discretization, $ (p_c+1) $ number of basis functions have local support in a contact element, see right column in Fig.~\ref{fig:VO_krefine}. However, as shown in Sec.~\ref{sec:numerical_example}, very smooth basis functions are incapable of accurately capture the sharp changes in the distribution of contact pressure result. This is related to the multi knot-span support of the smooth functions, which widens by an additional knot-span on linearly increasing their smoothness, e.g. from $ C^2 $ to $ C^3 $ as shown in Fig.~\ref{fig:VO_krefine}. Thus, the much higher-continuous NURBS does not considerably improve the accuracy of the contact responses results.

On the contrary, as the basis functions plots show in Fig.~\ref{fig:VO_kp_refine}, the application of additional one step of order elevation to the sufficiently $ C^1 $ smooth contact layer leads to the shrinking of the multi-knot span support of the basis functions. Further, the knot-span support of the basis functions remains unvaried on applying an additional or second step of order elevation to the contact layer. Such a sufficiently smooth and higher-order NURBS functions with a moderate knot-span support can adequately capture the local changes in the contact pressure results as shown in Sec.~\ref{sec:numerical_example}. It is also highlighted that although each type of N$ _2- $N$ _{p_c}$ $ (p_c = 3,~\mathrm{and}~ 4) $ discretization, and its equivalent version N$ _2- $N$ _{p_c \cdot p_s} $, $ (p_c \cdot p_s = 2\cdot1, ~\mathrm{and}~ 2\cdot 2) $, respectively, has same $ (p_c+1) $ or $ (p_c+p_s+1) $ number of non-zero basis functions across the contact boundary of an element. But, the moderate knot-span support of a N$ _2- $N$ _{p_c\cdot p_s} $ is accompanied by additional large number of overall control points across the contact interface as compared to N$ _2-$N$ _{p_c} $ discretization, respectively, see Figs.~\ref{fig:VO_krefine} and \ref{fig:VO_kp_refine}.

This way with the VO based NURBS discretization of a contact geometry: $ (i) $ the higher-order NURBS basis functions are utilized only for the evaluation of contact contributions in a fully NURBS discretized geometry, $ (ii) $ a large number of additional degrees of freedom are introduced across the contact surface without changing the mesh size, and $ (iii) $ the minimum order of NURBS are employed for the description of the bulk domain that does not come into contact. The resulting VO based NURBS discretized contact element is characterized by $ n_{cp}^{e}  = (p_c +1) + (p_1 +1)\times p_2 $ number of control points, where $ p_c + 1$ are present on the contact layer and $ (p_1 +1) \times p_2 $ in its remaining part. The bivariate NURBS basis functions for such an element are defined as
\begin{equation}
\label{eq:Enriched_NURBS_element}
\begin{aligned}	
R_{1}^{p_c,p_2} (\xi^1, \xi^2) &:=  \frac{w_{11}}{W(\xi^1,\xi^2)}  N_{1,p_c}(\xi^1)N_{1,p_2}(\xi^2)\,,  \\
\vdots~~~~~~~   & ~~~~~~~~~~~\vdots \\
R_{p_c+1}^{p_c,p_2} (\xi^1, \xi^2) &:=  \frac{w_{(p_c+1)1}}{W(\xi^1,\xi^2)} N_{p_c+1, p_c}(\xi^1)N_{1,p_2}(\xi^2)\,, \\	
R_{p_c+2}^{p_1,p_2} (\xi^1, \xi^2) &:= \frac{w_{12}}{W(\xi^1,\xi^2)} N_{1, p_1}(\xi^1)N_{2,p_2}(\xi^2)\,,   \\	
\vdots~~~~~~~   & ~~~~~~~~~~~\vdots \\
R_{n_{cp}^{e}}^{p_1,p_2}(\xi^1,\xi^2) &:= \frac{w_{(p_1+1)2}}{W(\xi^1,\xi^2)} N_{p_1+1,p_1}(\xi^1) N_{2,p_2}(\xi^2)\,, 
\end{aligned}	    
\end{equation} 
where the normalizing weight function is given by
\begin{equation}\label{key}
W(\xi^1,\xi^2) = {\sum_{i=1}^{(p_c+1)}  w_{i1} N_{i,p_c}(\xi^1) N_{1,p_2}(\xi^2)  + \sum_{i=1}^{(p_1+1)} w_{i2} N_{i,p_1}(\xi^1)N_{2,p_2}(\xi^2)}\,.
\end{equation}
The basis functions defined in Eq.~(\ref{eq:Enriched_NURBS_element}) exhibit the non-negativity property
\begin{equation}\label{eq:non_negative}
R_a^{p_c,p_2}(\boldsymbol{\xi}) \geq 0 ~~~~\forall \,\boldsymbol{\xi} \in \Omega^e_c ~~~~~~\textrm{where}~ a = 1,2,\dots,n_{cp}^e\,,
\end{equation}
and they also satisfy the partition of unity property: 
\begin{equation}\label{eq:partition_of_unity}
\sum_{a=1}^{n_{cp}^e} R_a^{p_c,p_2} (\boldsymbol{\xi}) = 1\,, ~~~~~~~ \forall \,\boldsymbol{\xi} \in \Omega^e_c\,.
\end{equation}
Using the isoparametric concept, the polynomials used for the VO based NURBS discretization, as defined in Eq.~(\ref{eq:Enriched_NURBS_element}), are employed for the approximation of the unknown displacement field $ \mathbf{u}^e $, its variation $ \delta \mathbf{u}^e $, and the current coordinates $ \bm{x}^e $ within each contact element as 
\begin{equation}\label{eq:contact_surface}
\boldsymbol{u}^{e} = \mathbf{R}(\boldsymbol{\xi})\,\mathbf{u}^{e}\,,~~~\delta \boldsymbol{u}^{e} = \mathbf{R}(\boldsymbol{\xi})\,\delta\mathbf{u}^{e}\,,~~~\bm{x}^{e} = \mathbf{R}(\boldsymbol{\xi})\,\boldsymbol{x}^{e} ~~~~~ \forall \, \boldsymbol{\xi} \in \Omega^e_c\,.
\end{equation}
Here, the new basis function matrix $ \mathbf{R}(\boldsymbol{\xi}) $ contains a total $ n_{cp}^e = (p_c +1) + (p_1 +1)\times p_2 $ number of bivariate NURBS functions $ R^{p_c, p_2}_a(\boldsymbol{\xi}) $ having local support in a VO discretized contact element $ \Omega^e_c $, see Eq.~(\ref{eq:Enriched_NURBS_element}). The array $ \mathbf{B}^{ke} $ in Eq.~(\ref{eq:element_int_force_vector}) containing the derivatives of NURBS functions is also modified according to Eq.~(\ref{eq:contact_surface}). Moreover, the matrix $ \mathbf{N}^{ke} $ in Eq.~(\ref{eq:contact_force_element}) makes use of higher-order of univariate NURBS functions $ R^{p_c}(\xi_1) $ for the evaluation of contact integrals.

The proposed VO based NURBS discretization approach bears some similarities with the NURBS-enriched contact finite element strategy by Corbett and Sauer~\cite{Corbett2014} and hybrid isogeometric-finite element based discretization technique by Maleki-Jebeli et al.~\cite{MALEKIJEBELI2018}. The common focus of these strategies is to improve the performance of the Lagrange-polynomial based finite element method in the context of contact mechanics. For this, the intrinsic features of the NURBS polynomials are coupled with the FE based discretizations. The B\'{e}zier extraction operator~\cite{Bezier_extraction_NURBS} is required to enable the incorporation of NURBS into the FE structure. In contrast, the present work aims to improve the performance of the NURBS-based IGA technique in the context of contact mechanics while fully retaining its original key purpose, i.e. the unified treatment of design and analysis processes by employing the same functions for the analysis that are used for the construction of a geometry exactly. For this purpose, the VO based NURBS discretization method is introduced. This unexplored idea, as initially noted by Temizer et al.~\cite{Temizer2012} presents a possibility to perform the controllable order-elevation based refinements of the NURBS discretized structures. It is noted that the proposed discretization method can also be converted to the enrichment strategy of Corbett and Sauer~\cite{Corbett2014} if the linear order of NURBS\footnote{As they are equivalent to linear order of Lagrange polynomials~\cite{Cottrell2009}} are chosen for the bulk description. The proposed method, however, unlike the enrichment strategy of~\cite{Corbett2014}, does not require the B\'{e}zier extraction operator.

\subsection{Implementation into existing code}
For integrating the VO NURBS discretization strategy into the existing isogeometric contact code, only a few minor modifications are required. First of all, for a given mesh resolution, a $ p_{c} > p_1 $ order of NURBS curve representing the contact boundary layer of the initial NURBS described geometry is constructed. After that, the parametrization for an originally $ p_1 $ NURBS discretized contact layer is replaced with that of the newly constructed $ p_{c} $ order of NURBS curve. For this, a certain number of conditions are need to be fulfilled. The total number of control points defining the geometry must be updated in such a manner that it allows the incorporation of the newly constructed $ p_{c} $ order of contact layer. This means that the connectivity array for contact elements must be adapted in a way that it contains the underlying control points of the VO based discretized geometry. The derived contact element connectivity arrays can have a different length than the bulk element connectivity arrays. The bivariate NURBS basis functions defined in Eq.~(\ref{eq:Enriched_NURBS_element}) are used for the evaluation of elemental quantities for VO based NURBS described contact elements. The univariate $ p_c $ order of NURBS functions are utilized for the evaluation of contact integrals. With the exception of these modifications, no other changes are need to be made in an existing isogeometric contact code. The local quantities, e.g. elemental stiffness matrices and force vector, are assembled to their global part in the same way as with the standard procedure. The reader is referred to Agrawal and Gautam~\cite{Agrawal2018} for a detailed description on the implementation of IGA in a simplified manner. In this work, a default $(p_c+1)\times(p_2+1)$ number of Gauss-Legendre quadrature points are employed for the evaluation of the contact integrals unless stated otherwise. Optimal quadrature rules~\cite{Hughes2010, AURICCHIO2012, FAHRENDORF2018}, which are well-suited for IGA, can also be opted for the reduced numerical evaluation.

\section{Numerical results} \label{sec:numerical_example}
In this section, we show the performance of the proposed VO NURBS based discretization method using three, two-dimensional small and large deformations contact problems. In the first example, a frictional ironing problem is analyzed to demonstrate the superior performance of the proposed method in terms of accuracy, robustness, and efficiency over the standard fixed-order based discretizations. The performance is assessed through by varying the inter-element continuity as well as the interpolation order of the NURBS discretizations using the $ k- $ and $ p-$refinements. In the second and third examples, Hertzian contact and two elastic rings problems are simulated, respectively. In these examples, we demonstrate the coarse-mesh accuracy and cost-efficacy of the proposed method for the evaluation of contact pressure distributions in comparison to standard NURBS discretizations. As compared to the second example, the last one incorporates large deformations and friction.

\subsection{Frictional ironing problem} \label{sec:ironing_problem}
\subsubsection{Problem setup} \label{sec:ironing_setup}
\begin{figure}[!t]
	\centering	
	\subfloat{\includegraphics[scale=0.65]{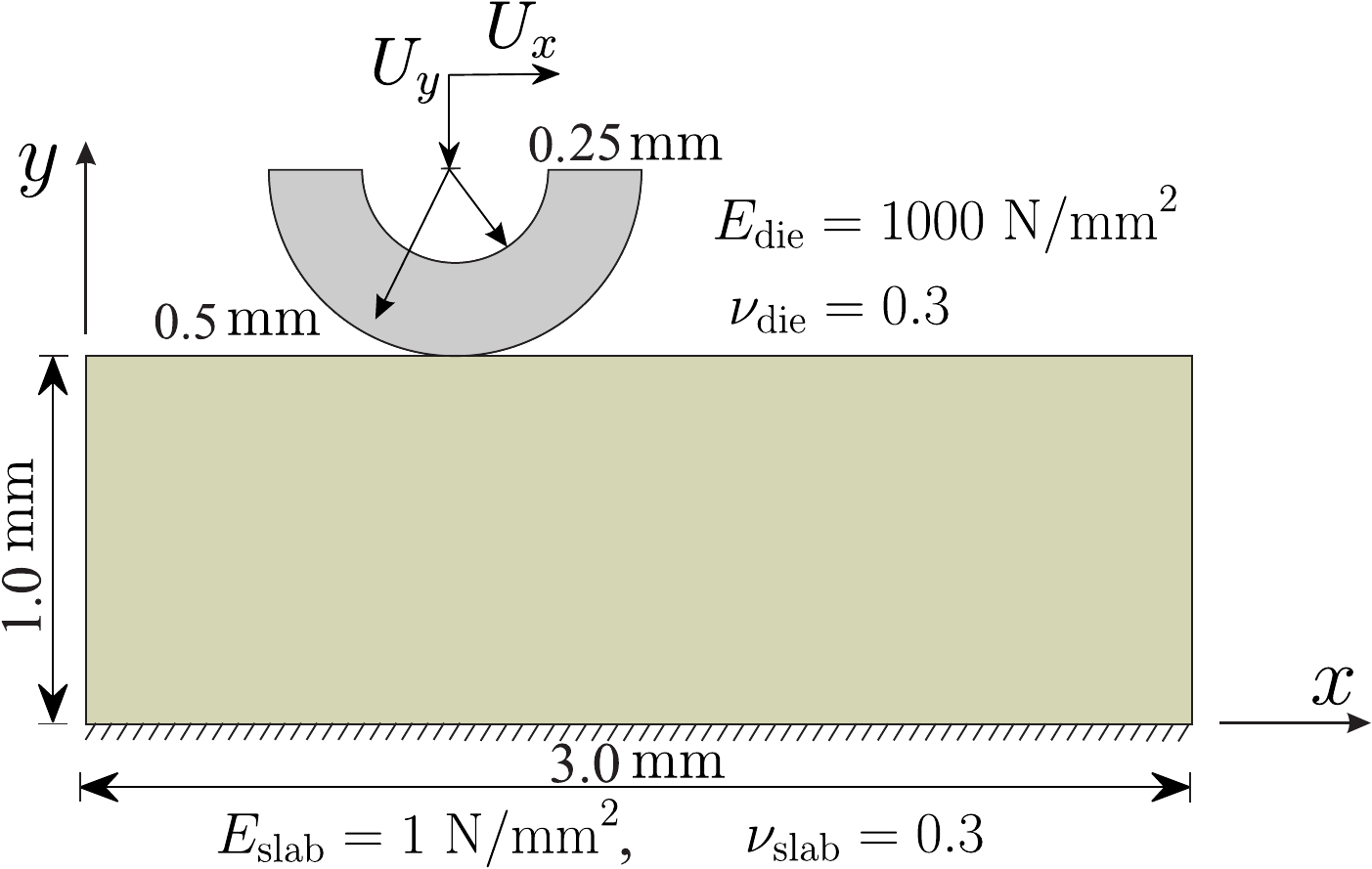}}  
	\caption{The setup of the ironing problem along with the geometric and material details, and the boundary conditions.} \label{fig:ironing_setup}
\end{figure}

In the first example, we consider the frictional contact between two deformable bodies using the setup similar to the one used by De Lorenzis et  al.~\cite{DeLorenzis2011}. The geometric model along with the material details and the boundary conditions are shown in Fig.~\ref{fig:ironing_setup}. This example is used to carefully analyze the performance of the proposed VO based NURBS discretization method over the standard fixed-order based NURBS discretization for large deformation and frictional sliding contact problem. The performance is measured in terms of the accuracy of the results, robustness, and the computational efficiency of the proposed approach.

In the considered example, the die is first pressed onto an elastic slab and then moved horizontally relatively, as shown in Fig.~\ref{fig:ironing_setup}.  A total vertical displacement $ U_y = -0.23 $ mm is applied to the top surface of the die in a uniform manner in $ 46 $ load steps. Thereafter the die is dragged across the slab by applying $ U_x = 1.5 $ mm in $ 250 $ load steps while keeping the vertical displacement fixed. An isotropic Neo-Hookean hyperelastic material model is considered to describe the material behaviour of the die and slab. The material parameters are shown in Fig.~\ref{fig:ironing_setup}. The penalty parameters are taken as $ \epsilon_{\textrm{N}} = \epsilon_{\textrm{T}} = 100  $ N/mm. The Coulomb's friction coefficient $ \mu_{\textrm{f}} = 0.2 $ is considered.
\begin{figure}[!h]
	\centering	
	\subfloat{\includegraphics[scale=0.28]{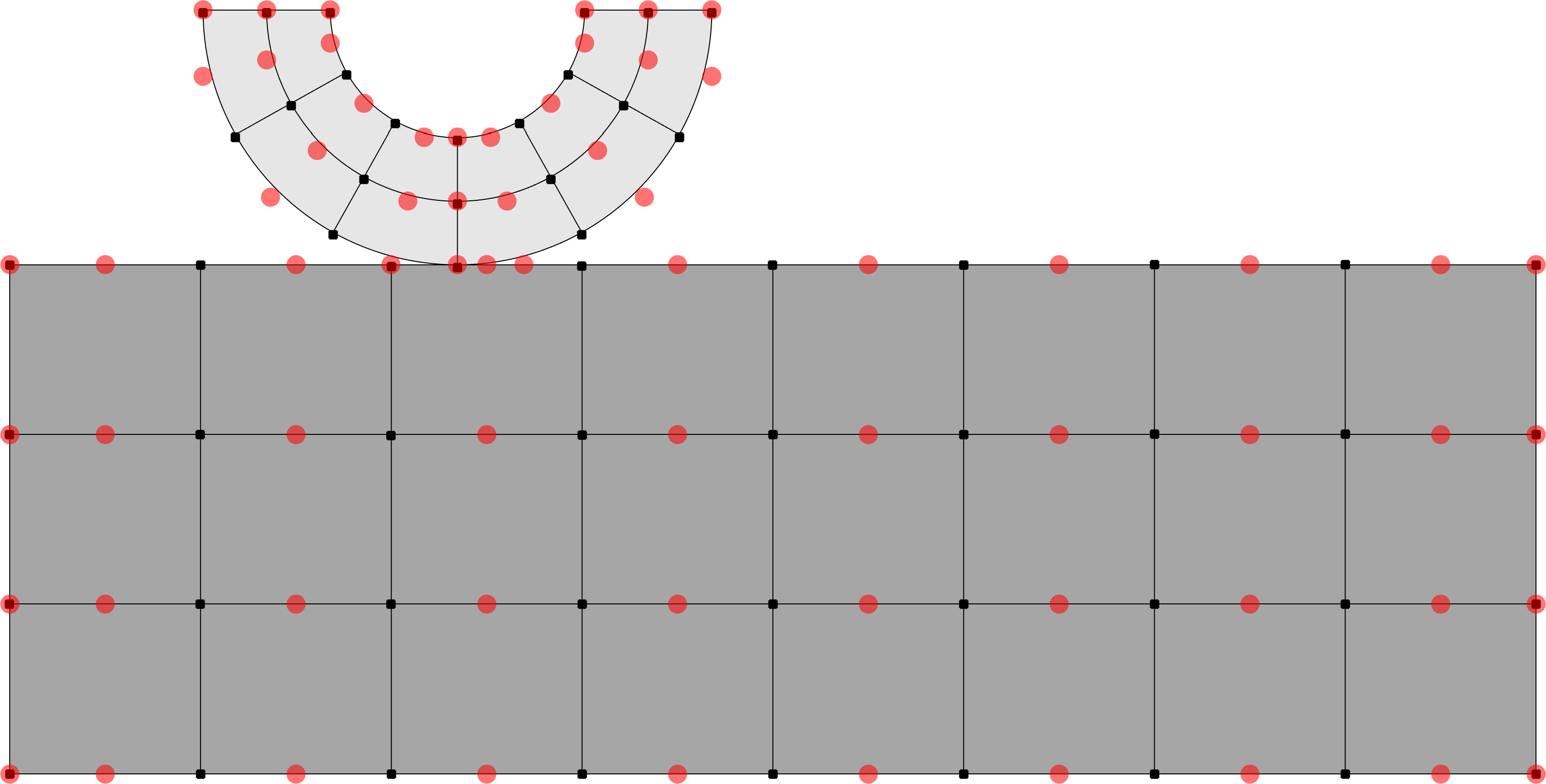}}  
	\caption{The coarsest mesh $ \mathrm{m}_1 $ used for the ironing problem. The control points associated with N$ _2 $ discretization are indicated with the red dots and the unique knot entries on the physical mesh of each body with the black squares.} \label{fig:mesh1_description}
\end{figure}
\begin{table}[!h]
	\begin{center}
		\begin{tabular}{c c c c c c c c}
			\hline
			\textbf{Mesh} &   &   &  &
			\multicolumn{4}{c}{\textbf{Elements}} \\ [0.4ex]
			\cline{5-8}		
			&  &  &  & Die &  &  & Slab \\  
			\hline
			m$ _1 $	&  &  &  & $ 6 \times 2 $  &  &  & $ 8 \times 3 $  \\
			m$ _2 $	&  &  &  & $ 12 \times 4 $ &  &  & $ 16 \times 6 $  \\
			m$ _3 $ &  &  &  & $ 24 \times 8 $ &  &  & $ 32 \times 12 $ \\
			m$ _4 $ &  &  &  & $ 48 \times 16 $ &  &  & $ 64 \times 24 $ \\
			\hline
		\end{tabular} \caption{Number of elements in different meshes used for ironing problem.} \label{table:mesh_size}
	\end{center}
\end{table}

First, a very coarse mesh is chosen to amplify the possible difference between the results with VO and standard NURBS discretizations, as discussed in the following. The considered coarsest mesh is denoted by $ \textrm{m}_1 $ and is illustrated in Fig.~\ref{fig:mesh1_description}. Three other nested meshes that are obtained through the uniform knot insertion and are used for the convergence study are listed in Table~\ref{table:mesh_size}. For a N$ _p $ order of NURBS contact layer, $ 3\times (p+1) $ number of Gauss-quadrature points are used for the evaluation of contact integrals per contact element.
\begin{figure}[!t]
	\centering
	\subfloat{\includegraphics[scale=0.22]{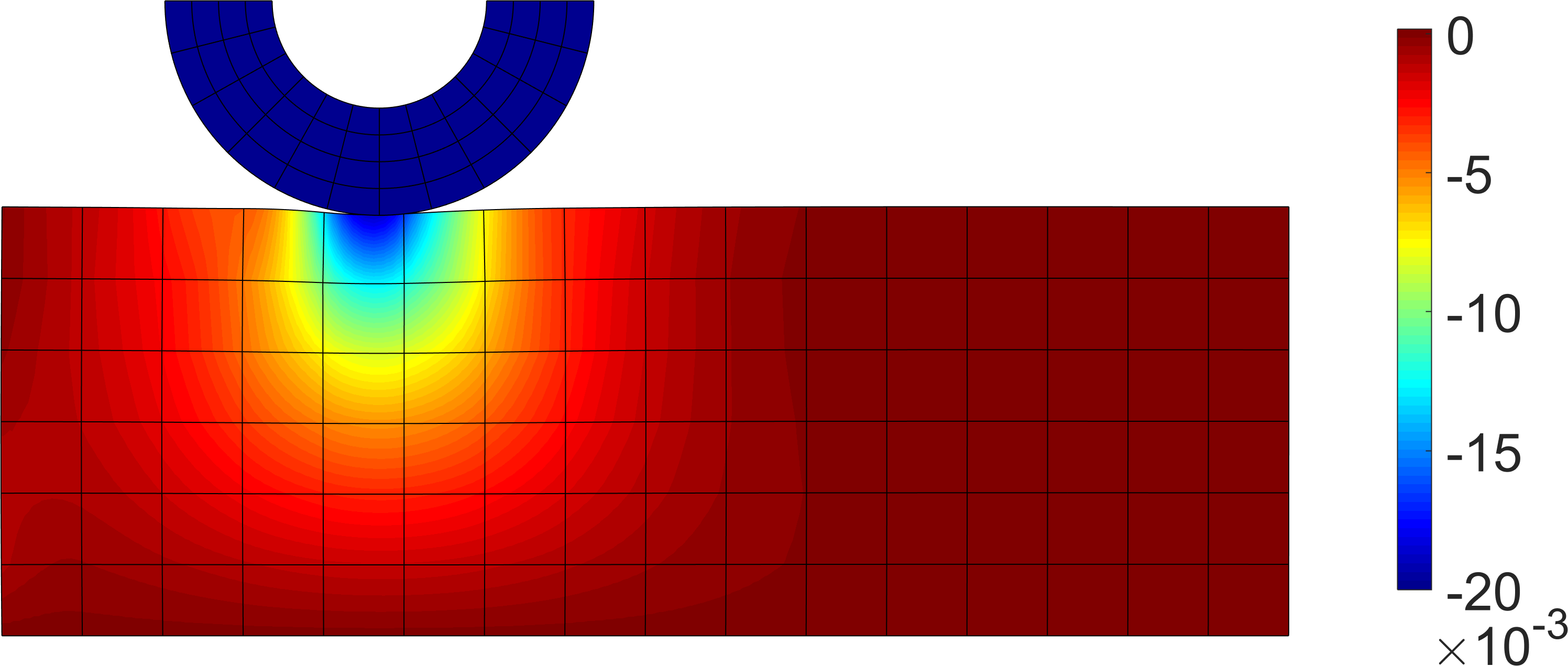}\label{fig:ironing_disp1}} ~~	
	\subfloat{\includegraphics[scale=0.22]{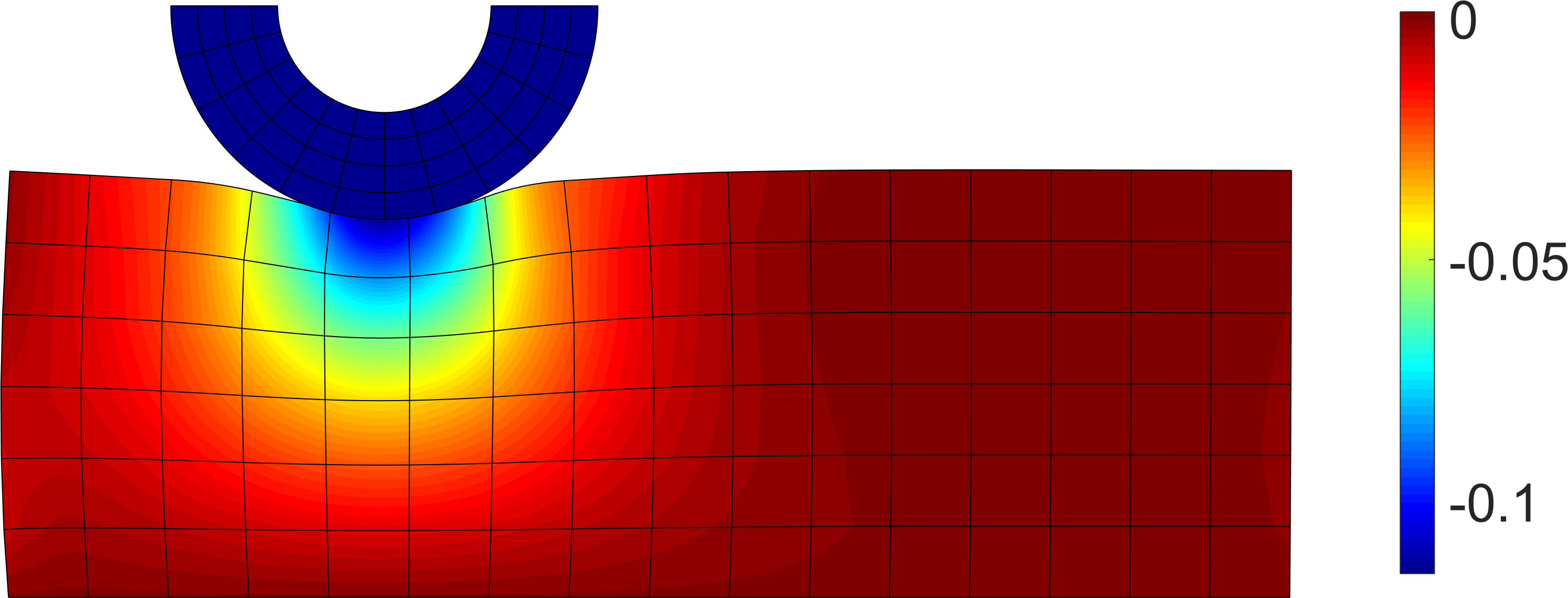}\label{fig:ironing_disp2}} \\
	\subfloat{\includegraphics[scale=0.22]{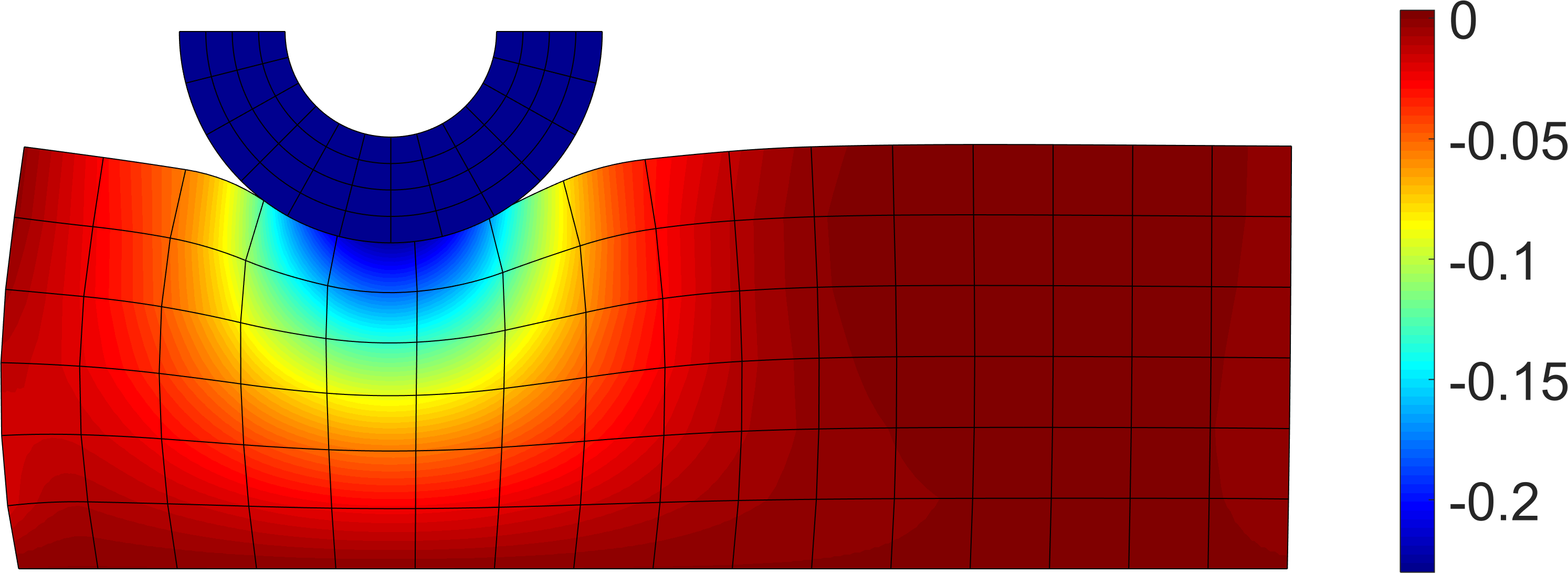}\label{fig:ironing_disp3}} ~~
	\subfloat{\includegraphics[scale=0.22]{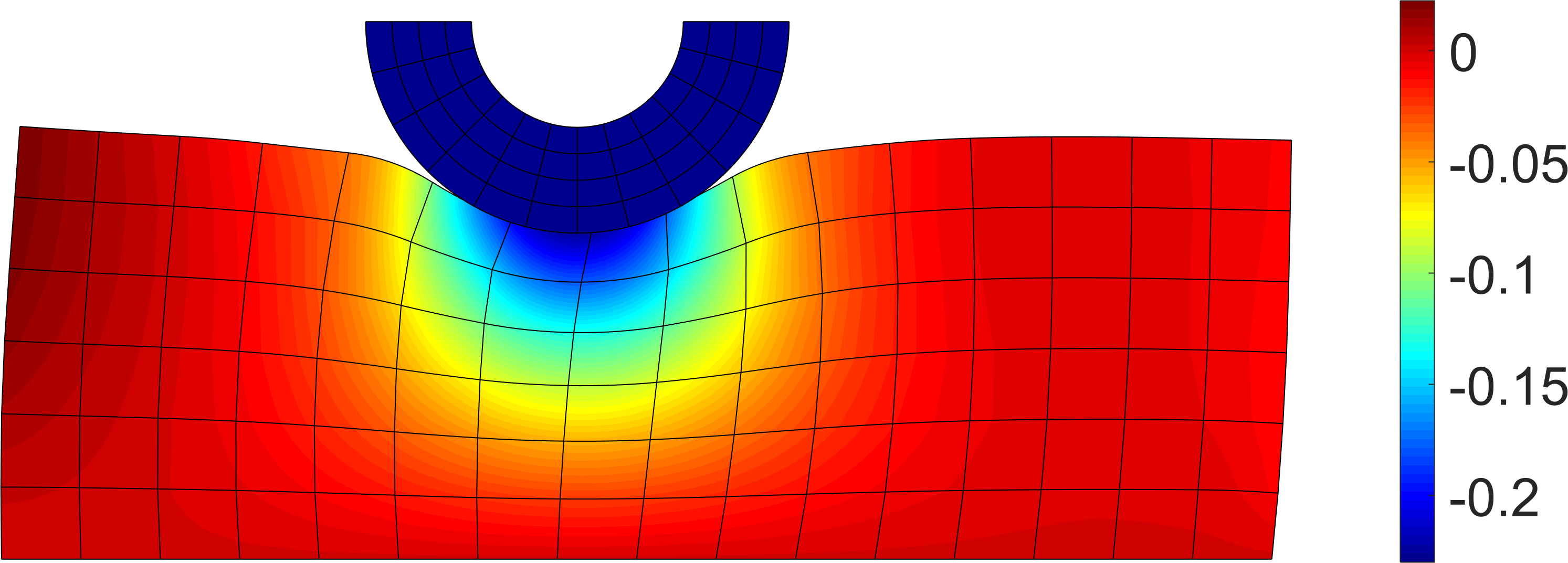}\label{fig:ironing_disp4}} \\
	\subfloat{\includegraphics[scale=0.22]{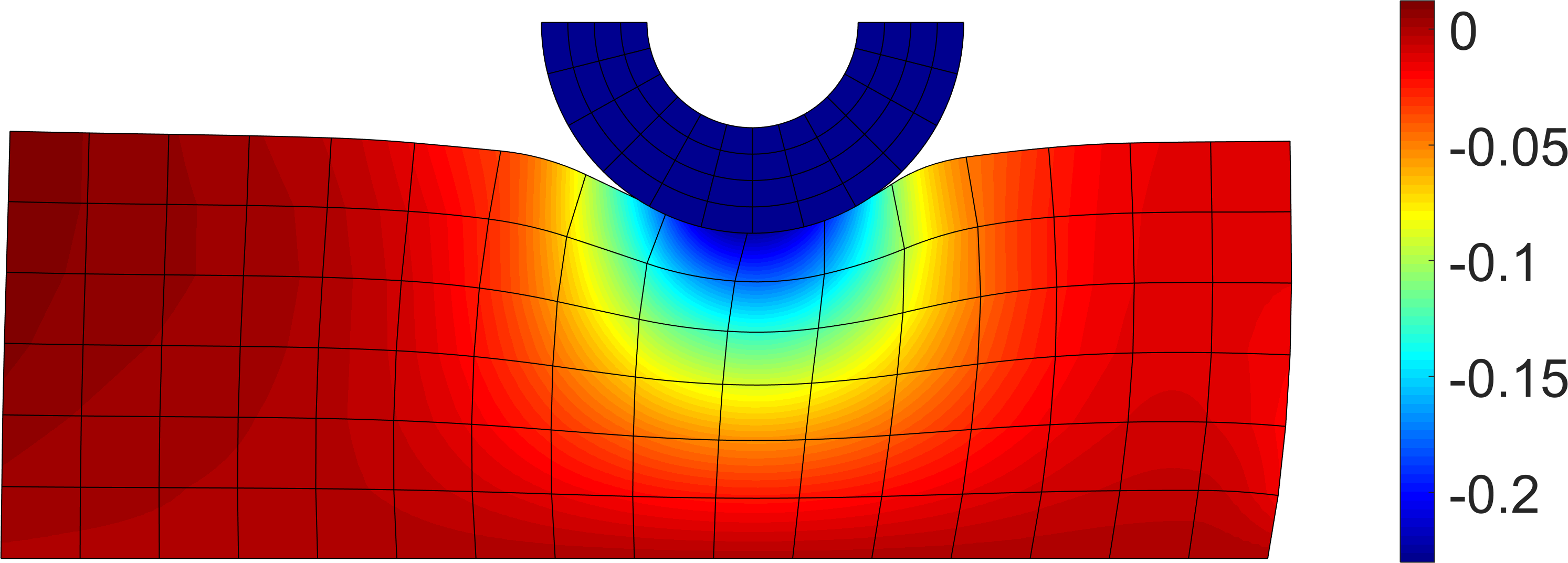}\label{fig:ironing_disp5}} ~~
	\subfloat{\includegraphics[scale=0.22]{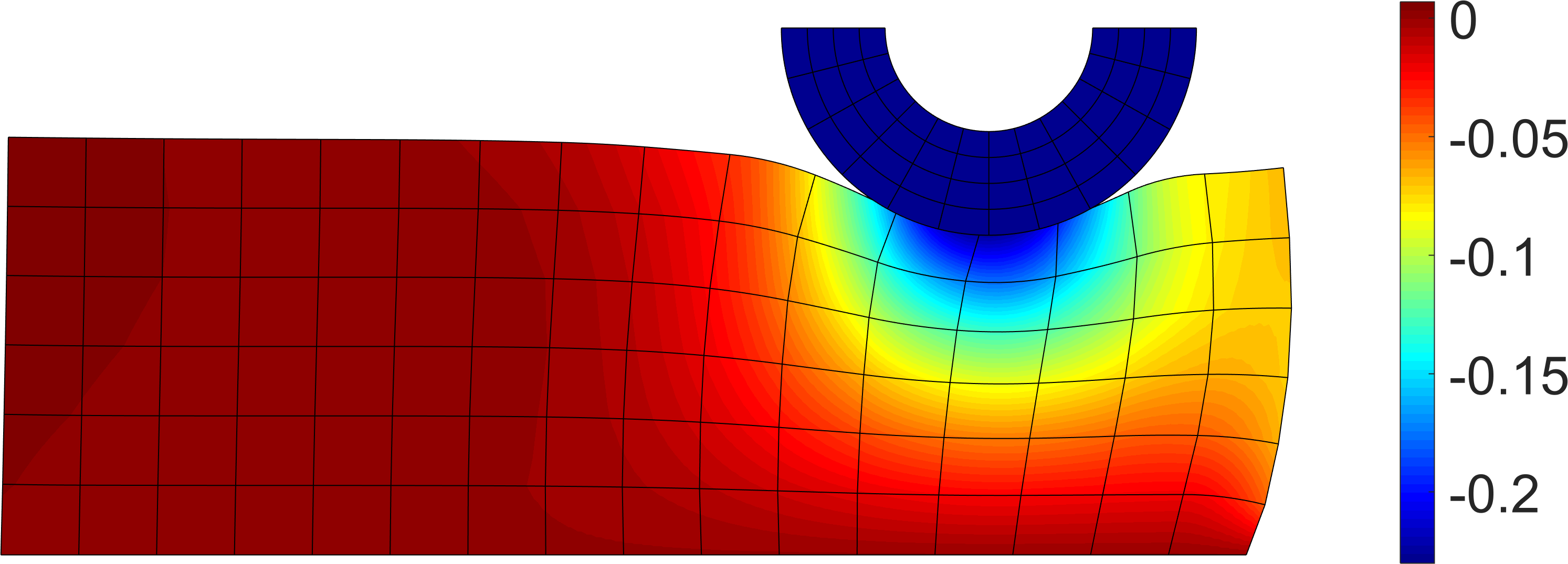}\label{fig:ironing_disp6}}
	\caption{Distribution of displacement field $ u_y $ in the deformed configuration during the compression process at step $ t = {4, 23, \textrm{ and } 46} $ (first three frames) and during the sliding at step $ t = 80, 160, \textrm{ and } 250$ (last three frames) with N$ _2 $ at mesh m$ _2 $.} \label{fig:ironing_disp}
\end{figure}

\subsubsection{Prediction of vertical and horizontal contact forces} \label{sec:contact_reaction_forces}
In this section, we analyze the improvement in the accuracy of the results on increasing the inter-element continuity and interpolation order of the NURBS with the proposed discretization method. The results with standard NURBS based discretizations are used for comparisons.

\paragraph{On increasing the inter-element continuity} \label{sec:ironing_21}
First, the quality of the results with different VO based discretizations: N$_2- $N$ _{p_c} $ ($ p_c = 4,$ and $ 6 $) is compared with that of corresponding standard N$ _p ~(p = 2, 4$ and $ 6) $ order of NURBS discretizations on increasing the inter-element continuity of the NURBS at the coarsest mesh $ \textrm{m}_1 $. In case of VO, the quadratic order of NURBS, which are sufficient to describe the setup exactly, are kept fixed for the discretization of the bulk region, while the interpolation order of the bottom surface of the die and upper surface of the slab are elevated to $ 4 $ and $ 6 $. The results with N$ _3 $ (and N$ _5 $) are omitted from this investigation, as they are similar to those obtained with N$ _4 $ (and N$ _6 $) cases.

\begin{figure}[!t]
	\centering
	\subfloat[]{\includegraphics[width=.55\linewidth]{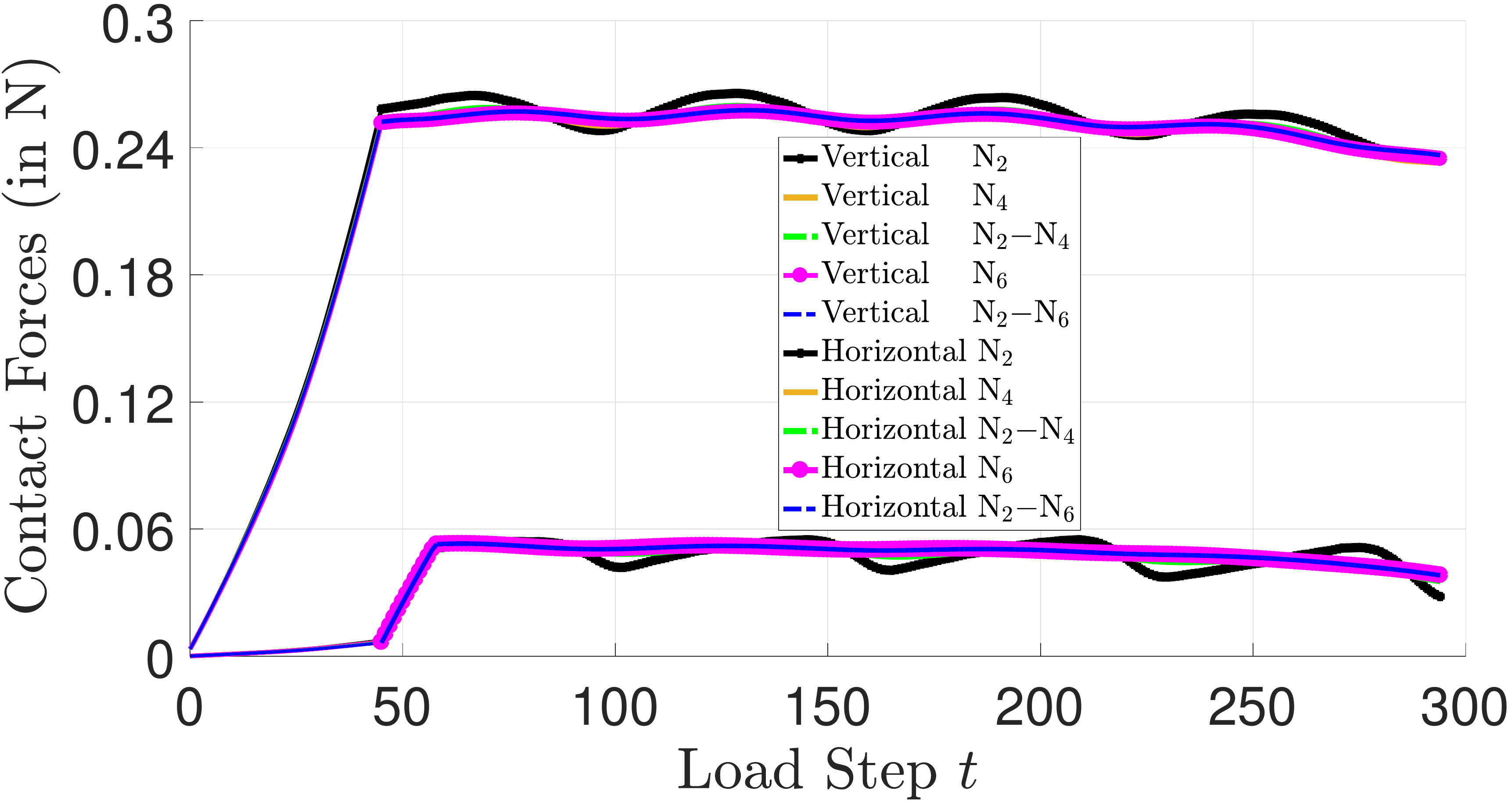} \label{fig:Full}}  \\
	\subfloat[]{\includegraphics[width=.49\linewidth]{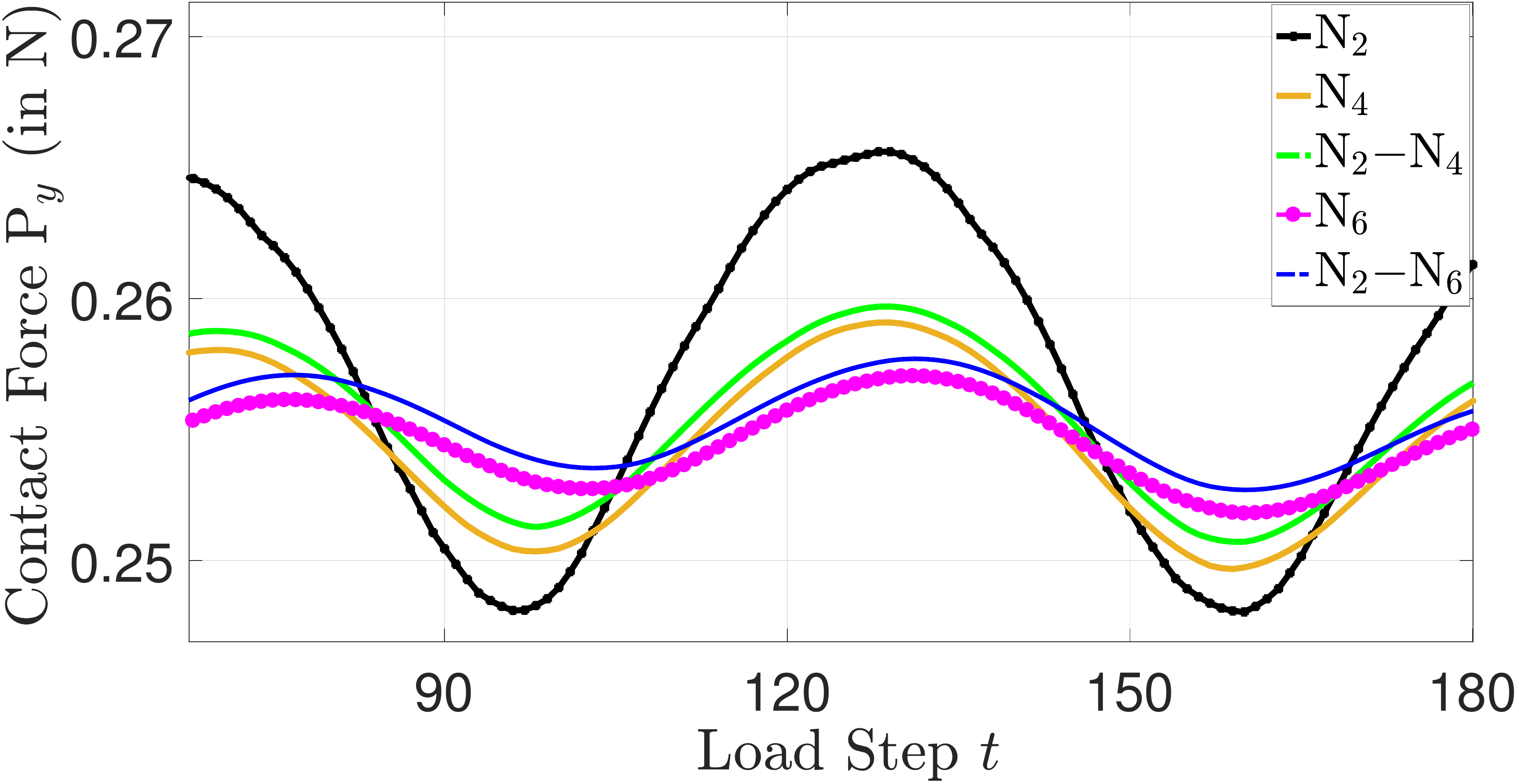}\label{fig:Full_Vertical}}\\
	\subfloat[]{\includegraphics[width=.48\linewidth]{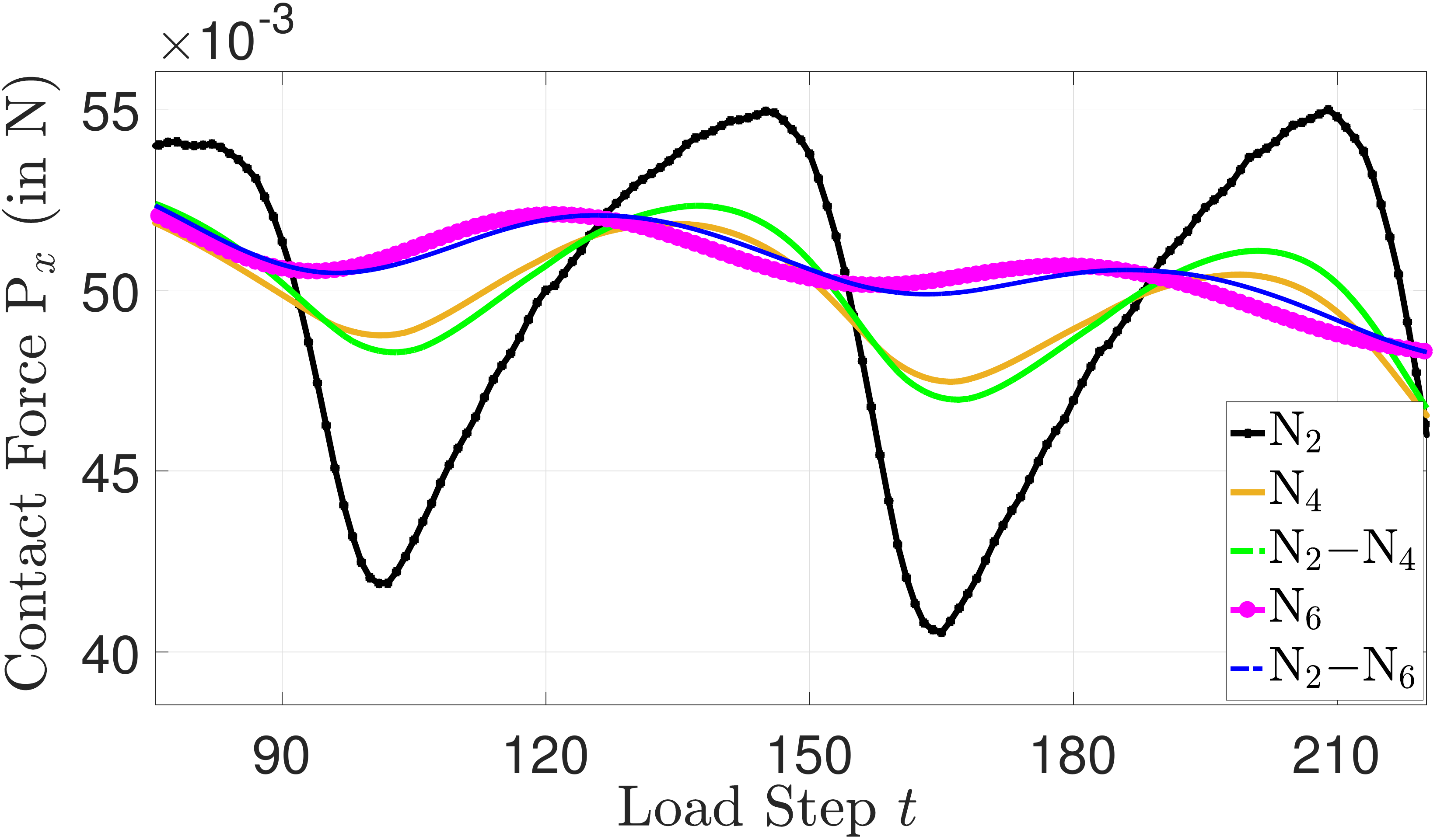}\label{fig:Full_Horizontal}}
	\caption{Ironing problem: (a) Computed total horizontal and vertical contact forces with load step $ t $ for different discretizations at the coarsest mesh m$_1 $. Enlarged view: (b) vertical contact force P$ _y $, and (c) horizontal contact reaction force P$ _x $ for N$ _2- $N$_{p_c}~(p_c = 4,\textrm{ and }6) $ and standard N$ _p~(p=2, 4, \textrm{ and }6) $ discretizations.}  \label{fig:Standard_Plots}
\end{figure}

The deformed configurations of the setup during downward and sliding motions are shown in Fig.~\ref{fig:ironing_disp} with N$ _2 $ using mesh $ \textrm{m}_2 $. During the indentation phase, the sum of contact forces on the top surface of the die should increase gradually, and during the sliding phase, forces should be approximately constant. Figure~\ref{fig:Full} shows the evolution of net vertical contact force P$ _y $ and the horizontal contact force P$ _x $, computed at the top surface of die, as a function of load step $ t $ for both N$ _p $ and N$_2- $N$ _{p_c} $ based discretizations. The enlarged views for P$ _y $ and P$ _x $ are shown in Figs.~\ref{fig:Full_Vertical} and~\ref{fig:Full_Horizontal}, respectively. Two major observations are made. First, the results verify the validity of the proposed VO based NURBS discretization method, as the contact force curves of $ P_y $ and $ P_x $ are nearly indistinguishable for VO based N$ _2-$N$ _{p_c}~(p_c= 4, \textrm{ and } 6) $ and standard N$ _p~ (p = 4, \textrm{ and } 6)$ discretizations. This is due to employing the same $ C^{p-1} -$continuous NURBS for the evaluation of contact integrals and identical number of DOFs present across the contact boundary in both the cases of discretizations. %
The DOF density data for N$ _2-$N$ _{p_c} $ and N$ _p$ is provided in Table~\ref{table:DOF_density_Irnoning} for mesh m$ _1 $. Second, during sliding, large-amplitude periodic oscillations of vertical and horizontal contact forces are present, especially with standard N$ _2 $ based discretization. The source of such large-amplitude non-physical oscillations is the insufficient discretization of the NURBS contact layers with the coarse mesh $ \textrm{m}_1 $. This, as a result, limits the conforming ability of the master body contact layer to the finite deformations introduced by the slave body. For a fixed mesh, conforming ability of contact layers increases on increasing the inter-element continuity with NURBS discretization. Consequently, oscillation error in P$_y  $ and P$_x $ reduces as shown in Figs.~\ref{fig:Full_Vertical} and~\ref{fig:Full_Horizontal}, respectively. With N$ _2-$N$_4 $ and N$ _2 -$N$_6 $, the oscillation amplitude of the horizontal contact force $ \Delta \textrm{P}_x $ reduce to approximately $ 27 \% $ and $ 12 \% $, respectively of that observed with N$_2 $ discretization. A quantitative analysis of the reduction in the oscillation error for N$ _2 -$N$ _{p_c}~(p_c=4,\textrm{ and }6) $ and N$ _p~(p=2,~3,\textrm{ and }4) $ for both the force components is provided in Table~\ref{table:amplitude_error}. The oscillation amplitude of the reaction forces is computed using  $ \Delta P_j := \text{max}(P_j) - \text{min}(P_j),~\text{where }j = x,\,y $. Here, the amplitude observed with N$_2 $ is used as a reference. It is noted that the obtained results and the reduction in the oscillation amplitude on increasing the smoothness of NURBS are in the agreement to those reported in the literature~\cite{DeLorenzis2011, DeLorenzis2012, Temizer2012, Temizer2014, Matzen2016, Temizer20162, Dimitri2017}. Moreover, the convergence behaviour and the reduction in the total analysis time with N$ _2 -$N$ _{p_c}$ as compared to N$ _p $ discretizations are reported in Sec.~\ref{sec:conv_beviour}.
\begin{table}[!hb]
	\begin{center}
		\begin{tabular}{|c| c |c| c| c| c| c | c |}
			\hline
			\textbf{Element} & 
			\multicolumn{3}{|c|}{\textbf{DOF for die}} & \multicolumn{3}{|c|}{\textbf{DOF for slab}} & \textbf{Total} \\[1ex]
			\cline{2-7} 		
			& \textbf{Interface} & \textbf{Bulk} & \textbf{Total} & \textbf{Interface} & \textbf{Bulk} & \textbf{Total} & \textbf{DOF}\\  
			\hline
			N$ _2 $& 18 & 36 & 54 & 20 & 60 & 80 & 134  \\
			N$ _4 $& 26 & 52 & 78 & 24 & 72 & 96 & 174 \\
			N$ _6 $& 34 & 68 & 102 & 28 & 84 & 112 & 214 \\
			N$ _2- $N$_{4} $& 26 & 36 & 62 & 24 & 60 & 84 & 146 \\
			N$ _2- $N$_{6} $& 24 & 36 & 70 & 28 & 60 & 88 & 158 \\
			N$ _2- $N$_{2\cdot1} $& 30 & 36 & 66 & 36 & 60 & 96 & 162 \\
			N$ _2- $N$_{2\cdot 2} $& 42 & 36 & 78 & 52 & 60 &  112 & 190 \\
			N$ _2- $N$_{2\cdot 3} $& 54 & 36 & 90 & 68 & 60 & 128 & 218 \\
			\hline
		\end{tabular} \caption{Degrees of freedom (DOF) density data for the die and slab for different VO and standard NURBS based discretizations at mesh $ \textrm{m}_1 $.} \label{table:DOF_density_Irnoning}
	\end{center}
\end{table}
\begin{table}[!h]
	\begin{center}
		\begin{tabular}{p{5cm} l p{3cm} l l}
			\hline
			\textbf{Element} &  & $ \Delta \textrm{P}_y\, (\%)$ &  & $ \Delta \textrm{P}_x \,(\%)$ \\ \hline
			N$ _2 $            		& & 100   & & 100   \\
			N$ _4 $  		  		& & 49.88 & & 23.31 \\
			N$ _6 $  		   		& & 24.57 & & 11.81 \\
			N$ _2- $N$ _{4} $  		& & 48.11 & & 27.02 \\			
			N$ _2- $N$ _{6} $  		& & 23.83 & & 12.09 \\
			N$ _2- $N$ _{2\cdot1} $ & & 5.98  & & 30.84 \\
			N$ _2- $N$ _{2\cdot 2} $& & 4.31  & & 7.51  \\
			N$ _2- $N$ _{2\cdot3} $ & & 1.89  & & 4.77  \\
			\hline
		\end{tabular} \caption{Reduction in the oscillation amplitude of vertical and horizontal contact forces for different standard and VO based NURBS discretizations at mesh  $ \textrm{m}_{1} $. The oscillation amplitude for N$ _2 $ is used as a reference.} \label{table:amplitude_error}
	\end{center}
\end{table}

\paragraph{On elevating the interpolation order}
Within the context of standard IGA, apart from increasing the inter-element continuity of NURBS discretization, the degree of conformity of the contact layer can be increased by increasing the mesh resolution. However, with the usage of a fine mesh, the  computational cost increases considerably, which is not desirable. This dismisses one of the key advantages of IGA which is its ability to represent the geometry to high-accuracy even with a coarse mesh. Therefore, the focus of the proposed VO NURBS discretization method is to improve the solution quality while keeping the mesh $ \textrm{m}_1 $ fixed. For this purpose, we carry out the additional order-elevation based refinement to N$ _2 $ discretized contact layers of the die and slab at mesh m$ _1 $. With this, although the inter-element continuity remains constant, a large number of additional DOFs are introduced across the contact layers at a fixed mesh attributed to the repetitions of knots with the order elevation process. Performing one and two steps of additional order-elevation to N$ _2 -$N$ _2 $ results in N$ _2 -$N$ _{2\cdot 1} $ and N$ _2 -$N$ _{2\cdot 2} $ discretizations, respectively. The total number of DOFs present in the contact interface and bulk domain corresponding to N$ _2 -$N$ _{2\cdot p_s}~(p_s = 1 \textrm{ and } 2) $ are listed in Table~\ref{table:DOF_density_Irnoning}.
\begin{figure}[!t]
	\centering	
	\subfloat[]{\includegraphics[width=.5\linewidth]{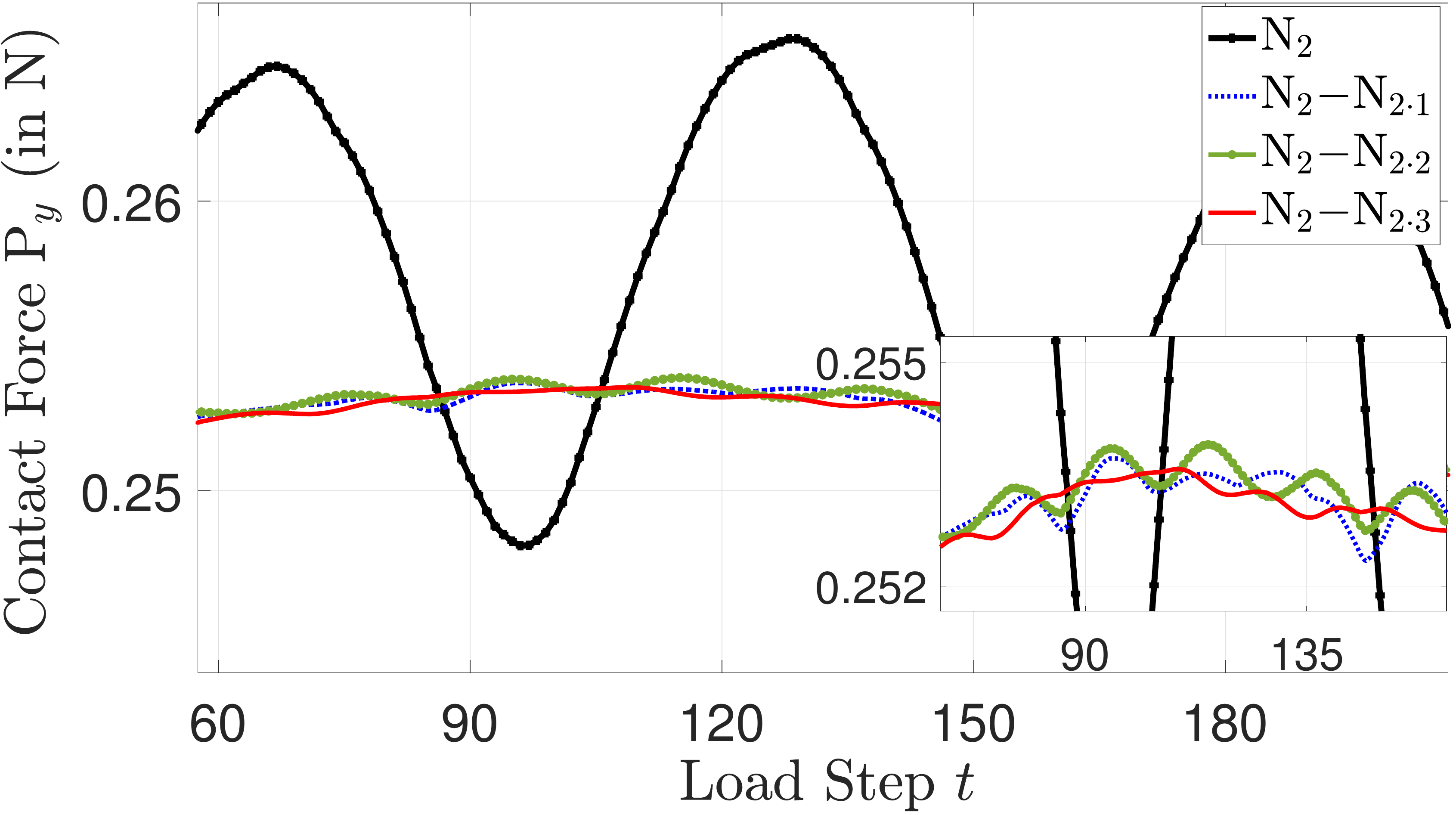} \label{fig:Vertical_N2p}} 	 \\
	\subfloat[]{\includegraphics[width=.5\linewidth]{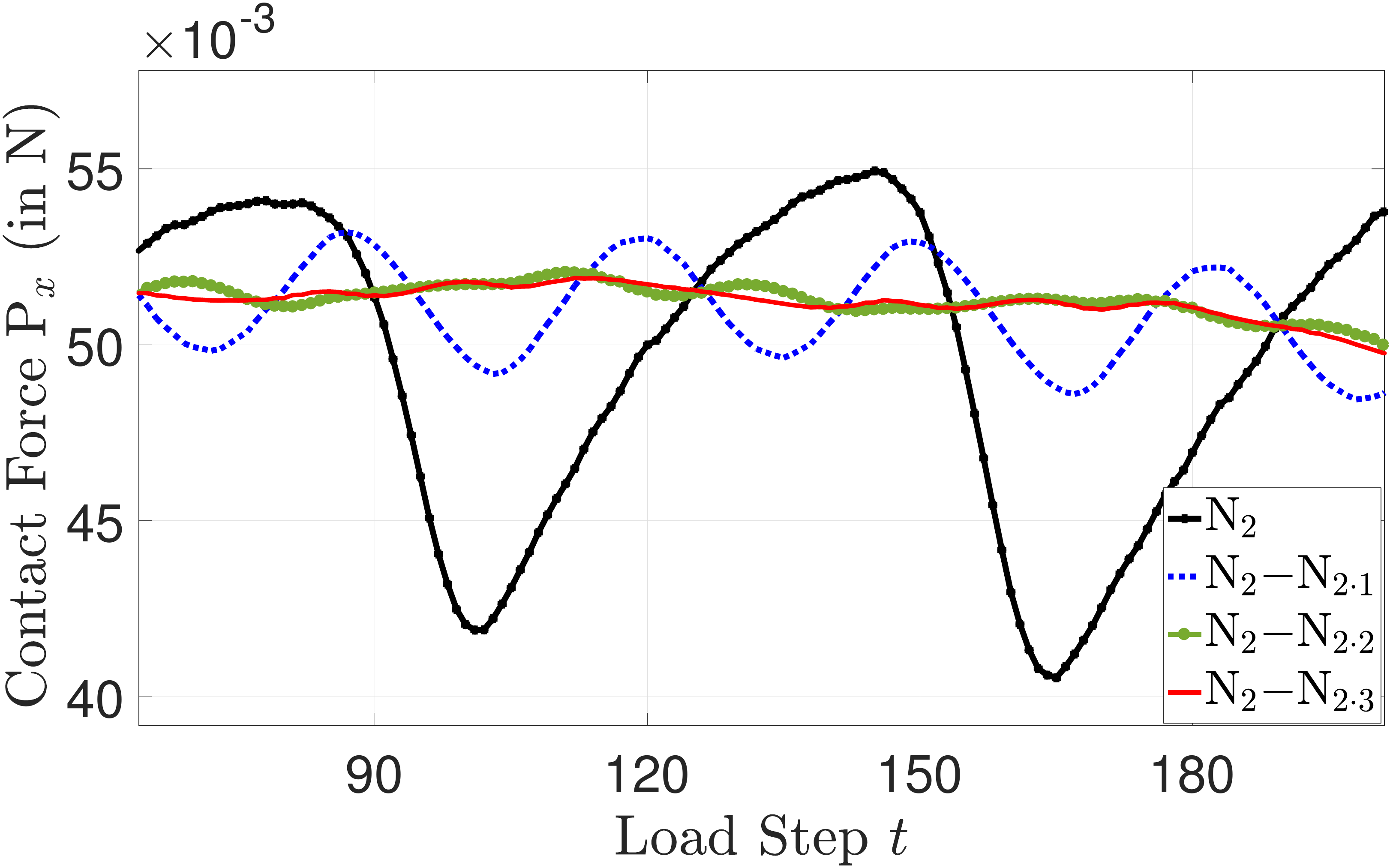}\label{fig:Horizontal_N2p} } 	
	\caption{Comparison of (a) vertical contact force P$ _y $ oscillation, and (b) horizontal contact  force P$ _x $ oscillation with different discretizations at mesh m$_1 $.} \label{fig:N2p}
\end{figure}
\begin{figure}[!t]
	\centering	
	\subfloat{\includegraphics[width=.55\linewidth]{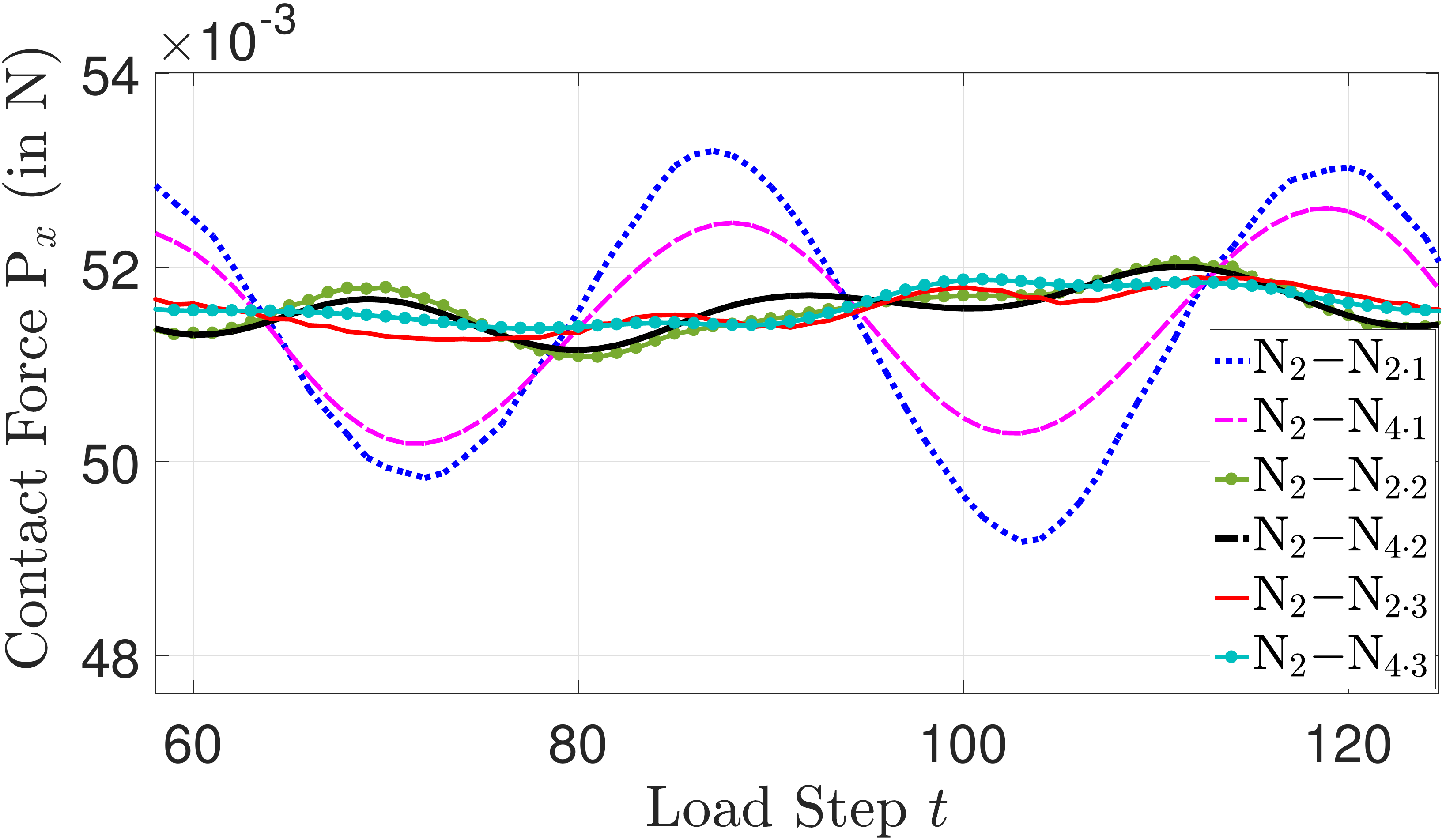}}
	\caption{Comparison of the horizontal contact force oscillation on increasing the inter-element continuity of the NURBS with different VO based discretizations at mesh m$ _1 $.} \label{fig:N2p_N4p_horizontal}
\end{figure}

Figures~\ref{fig:Vertical_N2p}, and~\ref{fig:Horizontal_N2p} show the contact force curve of P$ _y $ and P$ _x $, respectively, with VO based N$_2- $N$ _{2\cdot p_s}~(p_s = 1, \textrm{ and } 2) $ and standard N$ _2 $ based discretizations. As can be seen, the solution quality improves significantly for the new discretizations, especially with N$ _2- $N$ _{2\cdot 2} $. This is due to the moderate multi knot-span of basis functions as discussed in Sec.~\ref{subsec:proposed_technique} and a large number of DOFs present across the contact layer with N$_2- $N$ _{2\cdot p_s}~(p_s = 1, \textrm{ and } 2) $ as compared to N$ _2 $, see Table~\ref{table:DOF_density_Irnoning} for DOFs details. Moreover, the other higher order-elevation based discretization, i.e. N$ _2 -$N$ _{2\cdot 3} $, also reduces the oscillation error, as shown in Figs.~\ref{fig:Vertical_N2p} and~\ref{fig:Horizontal_N2p}. However, only a slight improvement is attained as compared to N$ _2 -$N$ _{2\cdot 2} $, see Table~\ref{table:DOF_density_Irnoning}. This shows that the number of DOFs more than with the N$ _2 -$N$ _{2\cdot 2} $ across the contact layer improves the results only slightly. The quantitative details on the reduction of oscillation amplitude for each discretizations for both the force components are shown in Table~\ref{table:amplitude_error}. 

Finally, we investigate the influence of order-elevation based refinement to the higher continuous NURBS contact layer, i.e. N$ _4 $, again while keeping the original mesh m$ _1 $ fixed. The contact force curve P$ _x $ corresponding to N$ _2- $N$ _{4\cdot p_s}~(p_s = 1,~2 \textrm{ and } 3) $ discretizations are shown in  Fig.~\ref{fig:N2p_N4p_horizontal}. The results with N$ _2- $N$ _{2\cdot p_s}~(p_s = 1,~2 \textrm{ and } 3) $ are used for the  purpose of comparison. The major observation is that on 
increasing the inter-element continuity of the N$ _{2\cdot p_s}$ discretized contact layer to N$ _2- $N$ _{4\cdot p_s}$, the oscillation error reduces marginally. This is expected as N$ _2- $N$ _{4\cdot p_s}$ has total only 12  additional DOFs across the contact interface as compared to N$ _2- $N$ _{2\cdot p_s}$, see Table~\ref{table:DOF_density_Irnoning} for DOFs. Furthermore, it can be observed that more smoother response of contact forces are obtained on increasing the continuity of NURBS discretized contact layer.

With the obtained results, it is evident that additional order-elevation based VO NURBS discretizations, particularly N$ _2 -$N$ _{2 \cdot 2} $, is advantageous over the standard N$ _p ~( p= 2, ~4, \textrm{ and }6)$ as well as N$ _2- $N$ _p $ based VO NURBS discretizations, see Table~\ref{table:amplitude_error}. The N$ _2 -$N$ _{2 \cdot 2} $ reduces the oscillation amplitude of vertical and horizontal contact forces to $ 4.31 \% $ and $ 7.51 \% $, respectively, as compared to N$ _2 $ at a fixed mesh $ m_1 $.

\subsubsection{Convergence behaviour and analysis time} \label{sec:conv_beviour}	
Next, we investigate the convergence of the vertical contact forces oscillation amplitude $ \Delta \textrm{P}_y $ and horizontal contact force oscillation amplitude $ \Delta \textrm{P}_x $ with both the discretization methods upon mesh refinement. Four nested meshes, as described in Table~\ref{table:mesh_size}, are used. In case of VO, the higher-continuous and additional order-elevation based discretizations are considered. Moreover, the total analysis time taken by VO based discretizations is compared to that with the standard NURBS discretizations for each mesh level.

\begin{figure}[!t]
	\centering		
	\subfloat[]{\includegraphics[width=.46\linewidth]{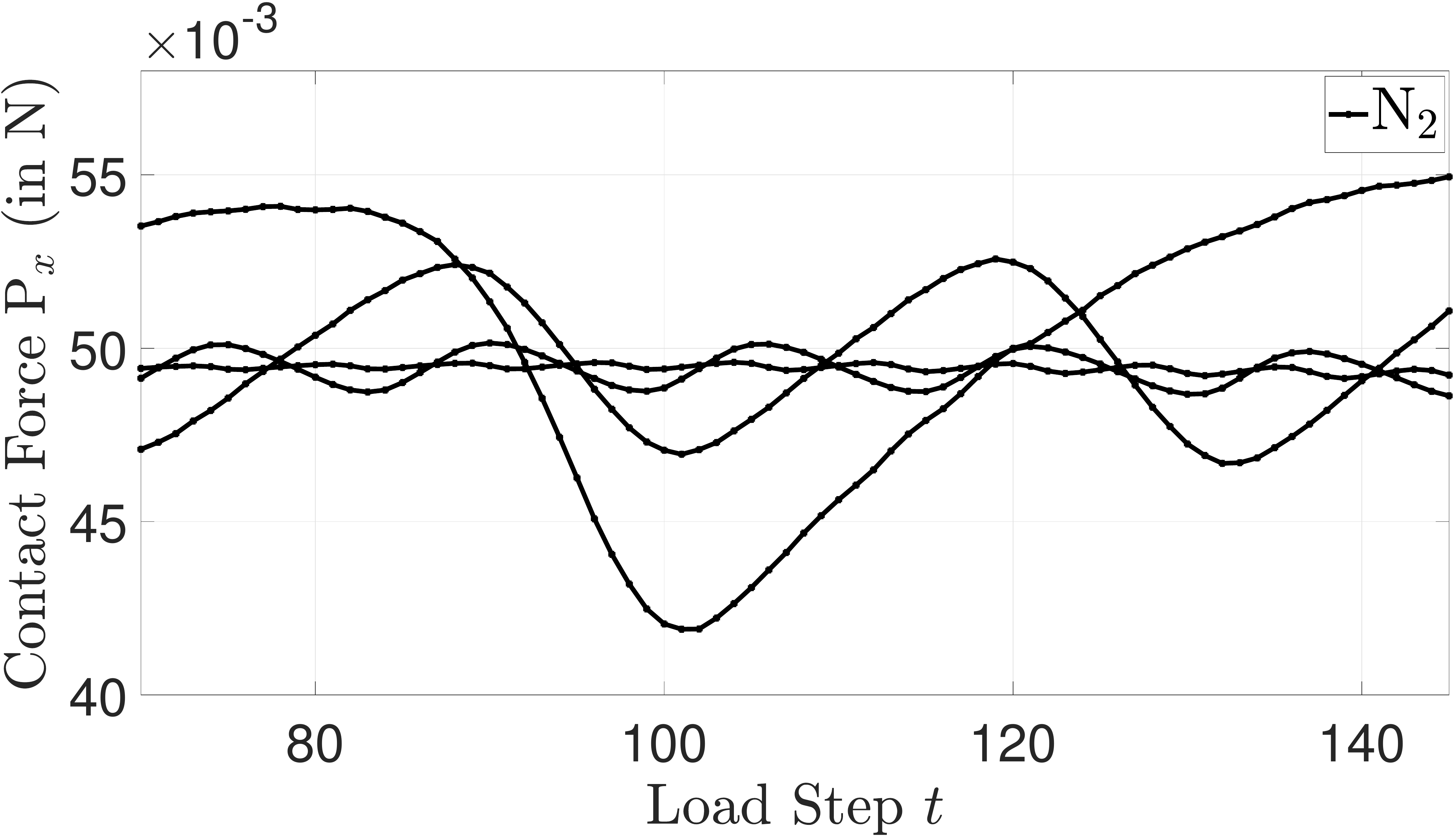} \label{fig:N2_conv}}  ~~ \subfloat[]{\includegraphics[width=.462\linewidth]{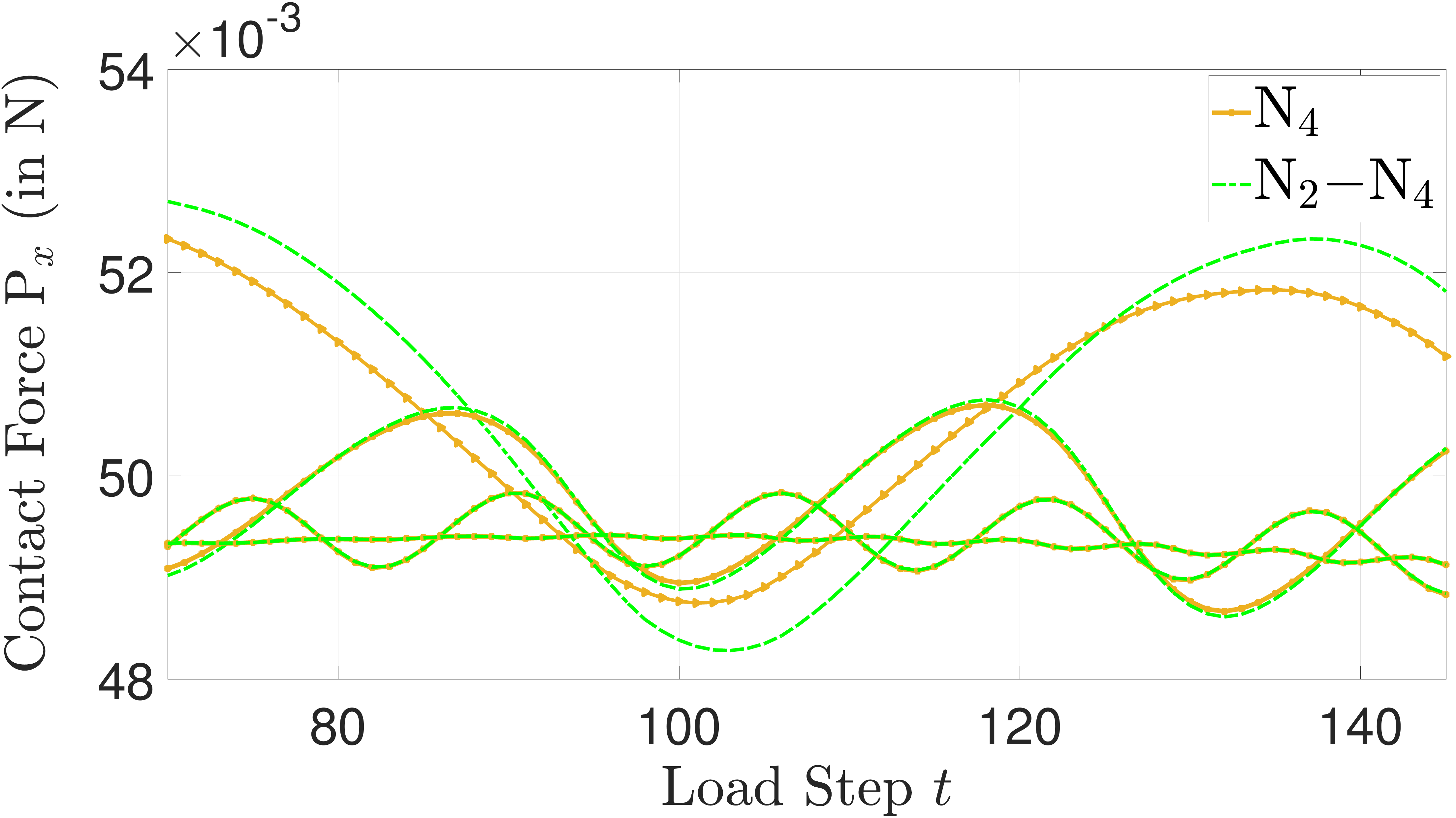} \label{fig:N4_conv}} \\
	\subfloat[]{\includegraphics[width=.46\linewidth]{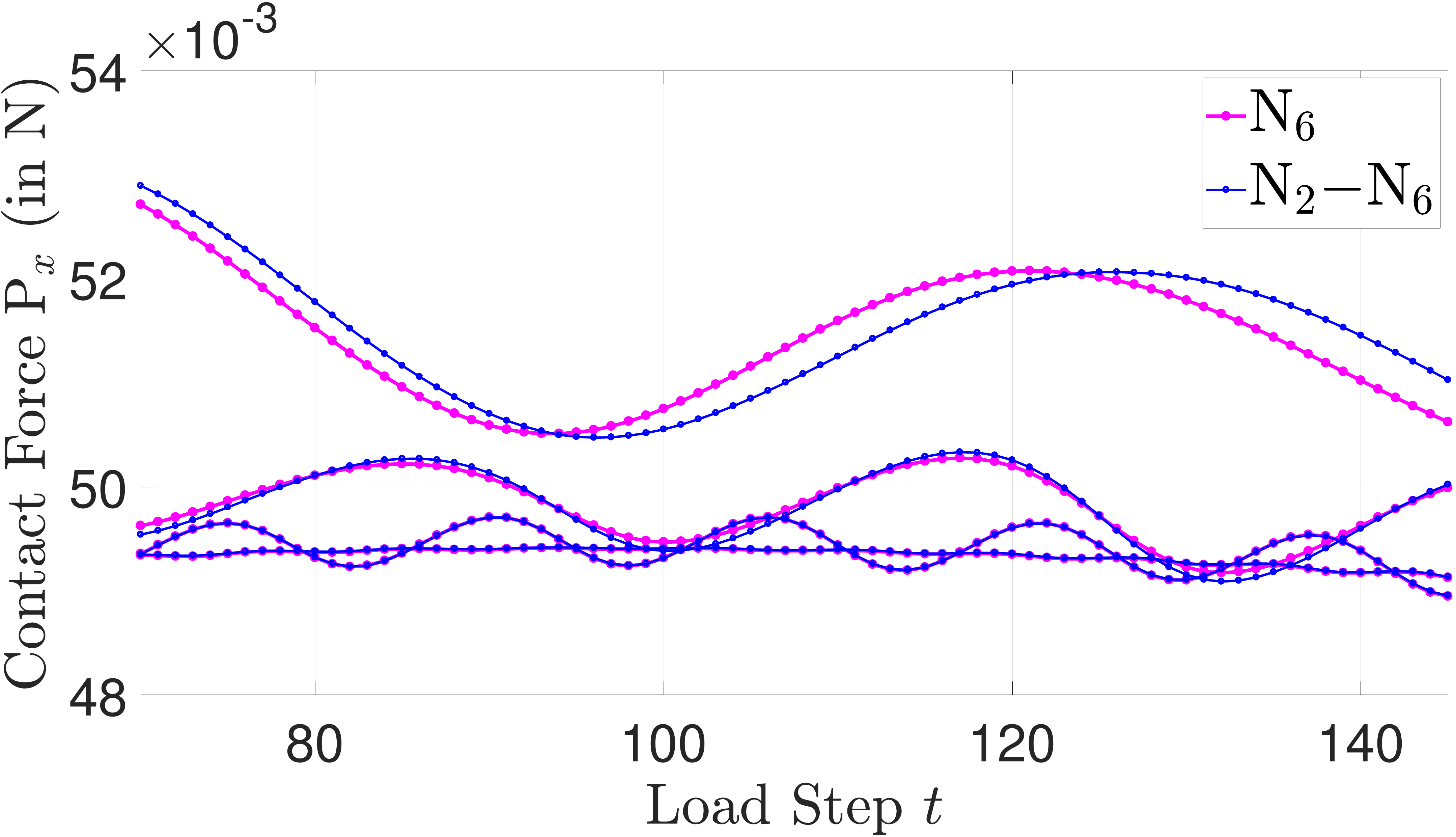} \label{fig:N6_conv}} ~~
	\subfloat[]{\includegraphics[width=.46\linewidth]{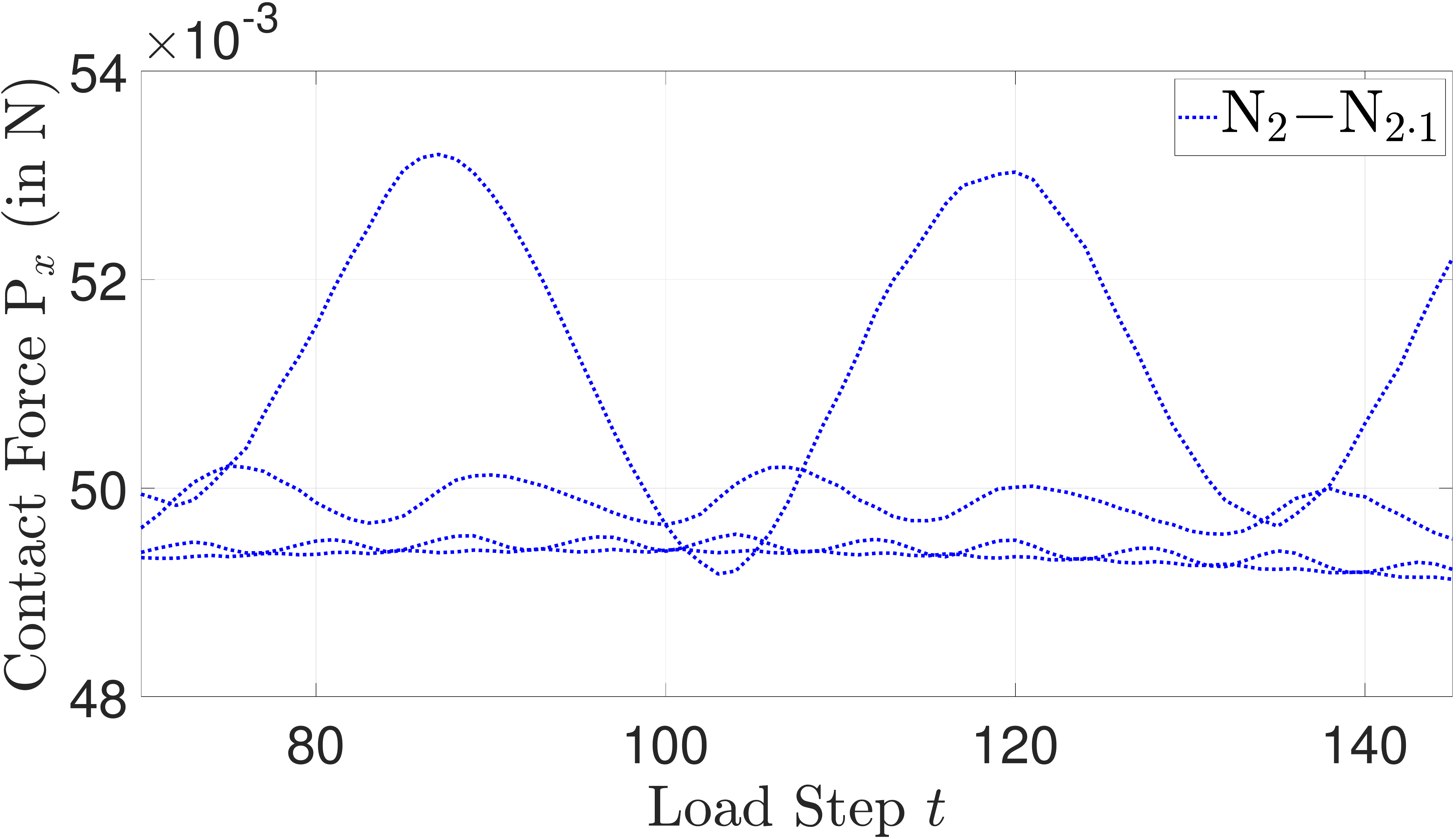} \label{fig:N21_conv}} \\
	\subfloat[]{\includegraphics[width=.46\linewidth]{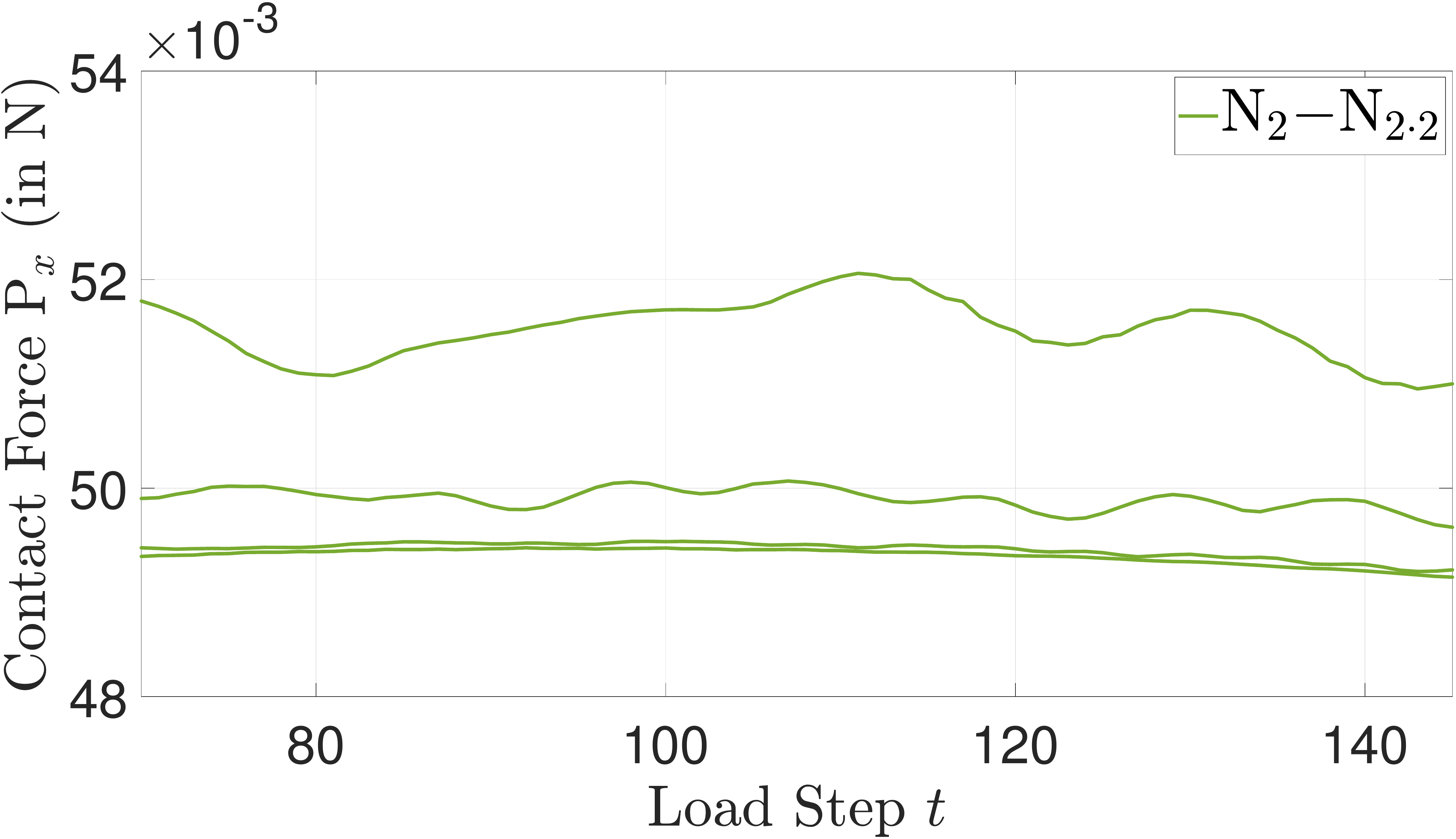} \label{fig:N22_conv}} ~~
	\subfloat[]{\includegraphics[width=.46\linewidth]{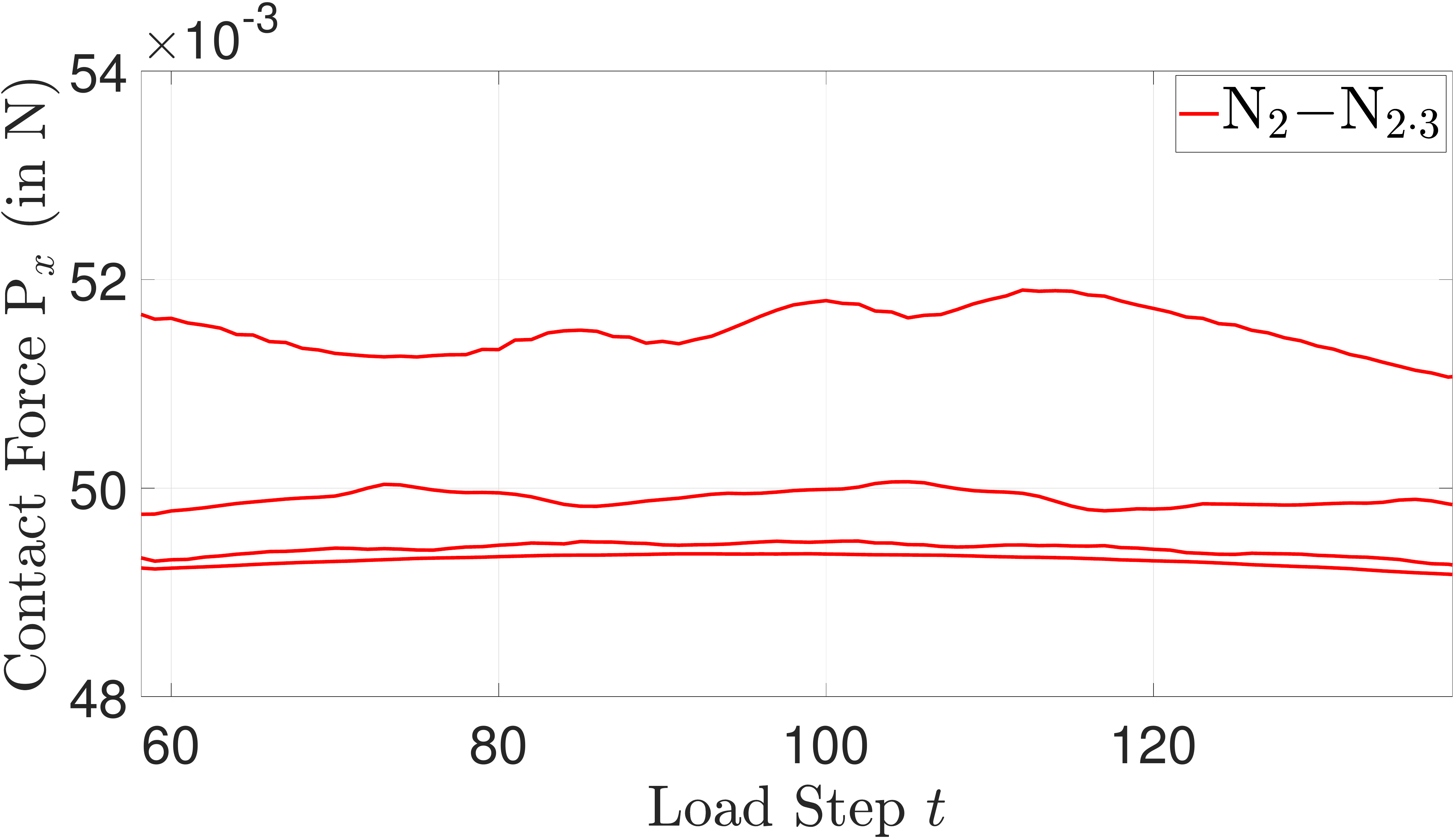} \label{fig:N23_conv}} 
	\caption{Reduction of the horizontal contact force oscillation on increasing the mesh resolution with both the standard and VO based NURBS discretizations.}  \label{fig:conv_behaviour}
\end{figure}

\begin{figure}[!t]
	\centering
	\subfloat[]{\includegraphics[width=.52\linewidth]{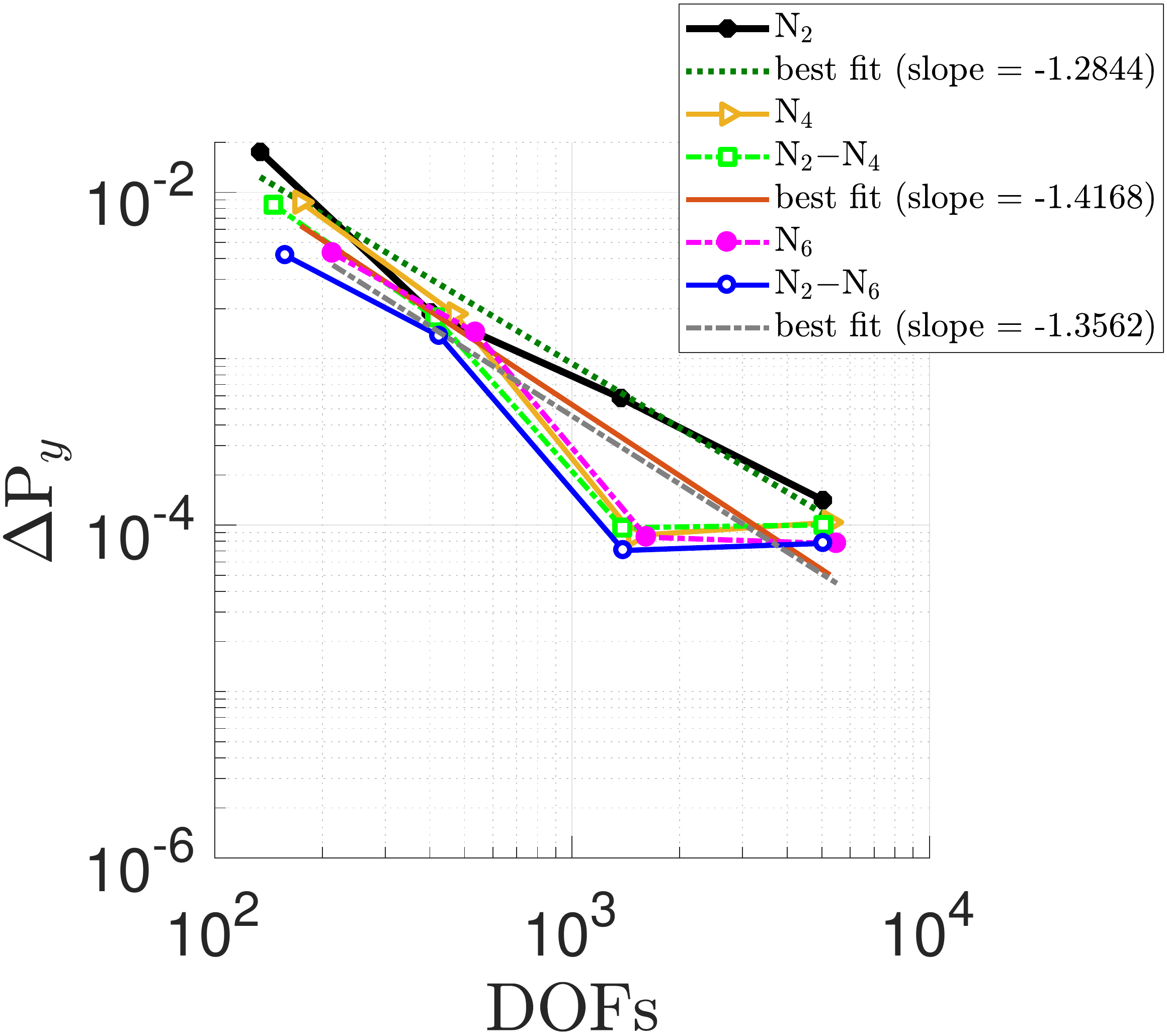}\label{fig:dPy_continuity}} 
	\subfloat[]{\includegraphics[width=.52\linewidth]{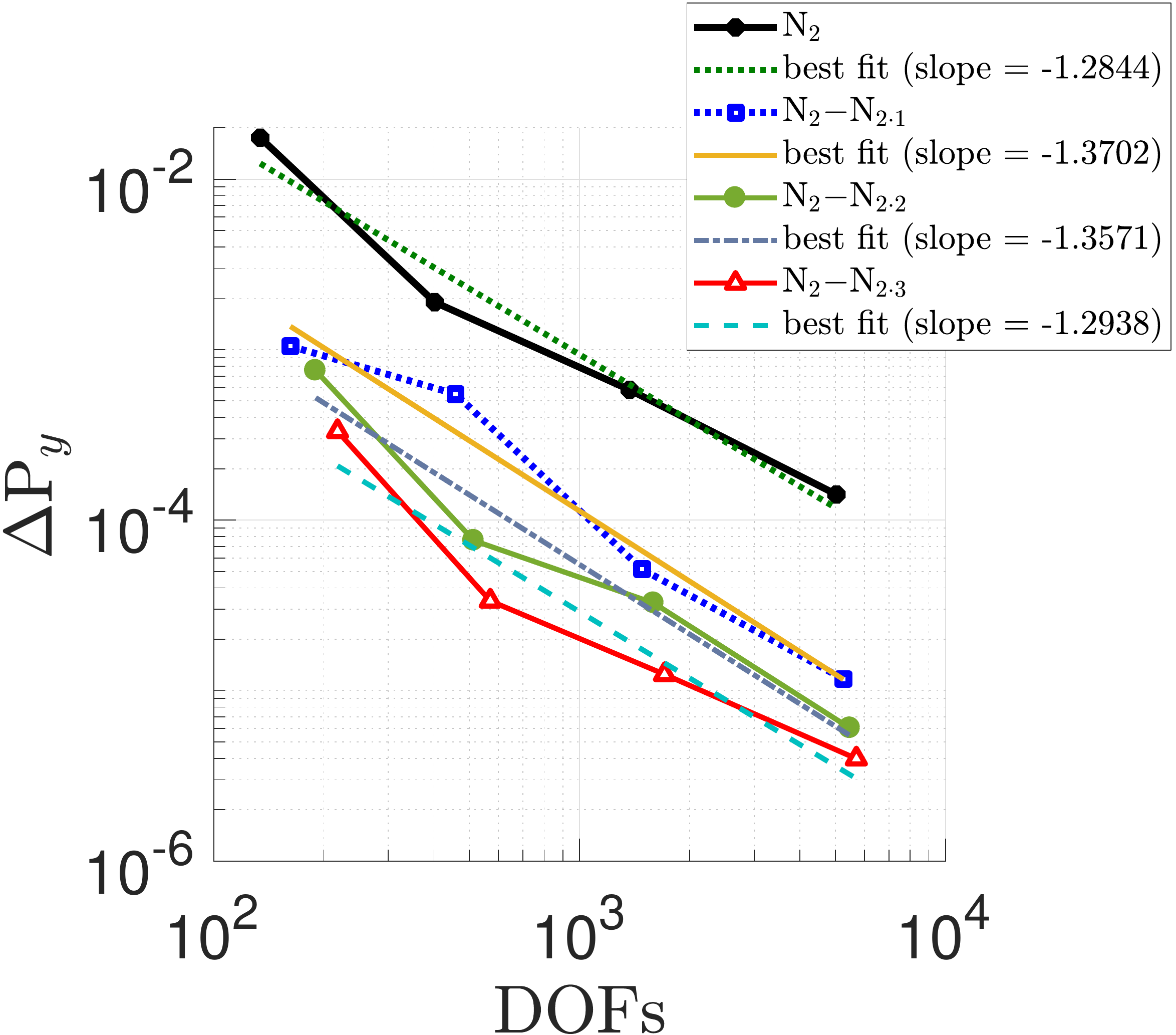} \label{fig:dPy_order}} \\
	\subfloat[]{\includegraphics[width=.52\linewidth]{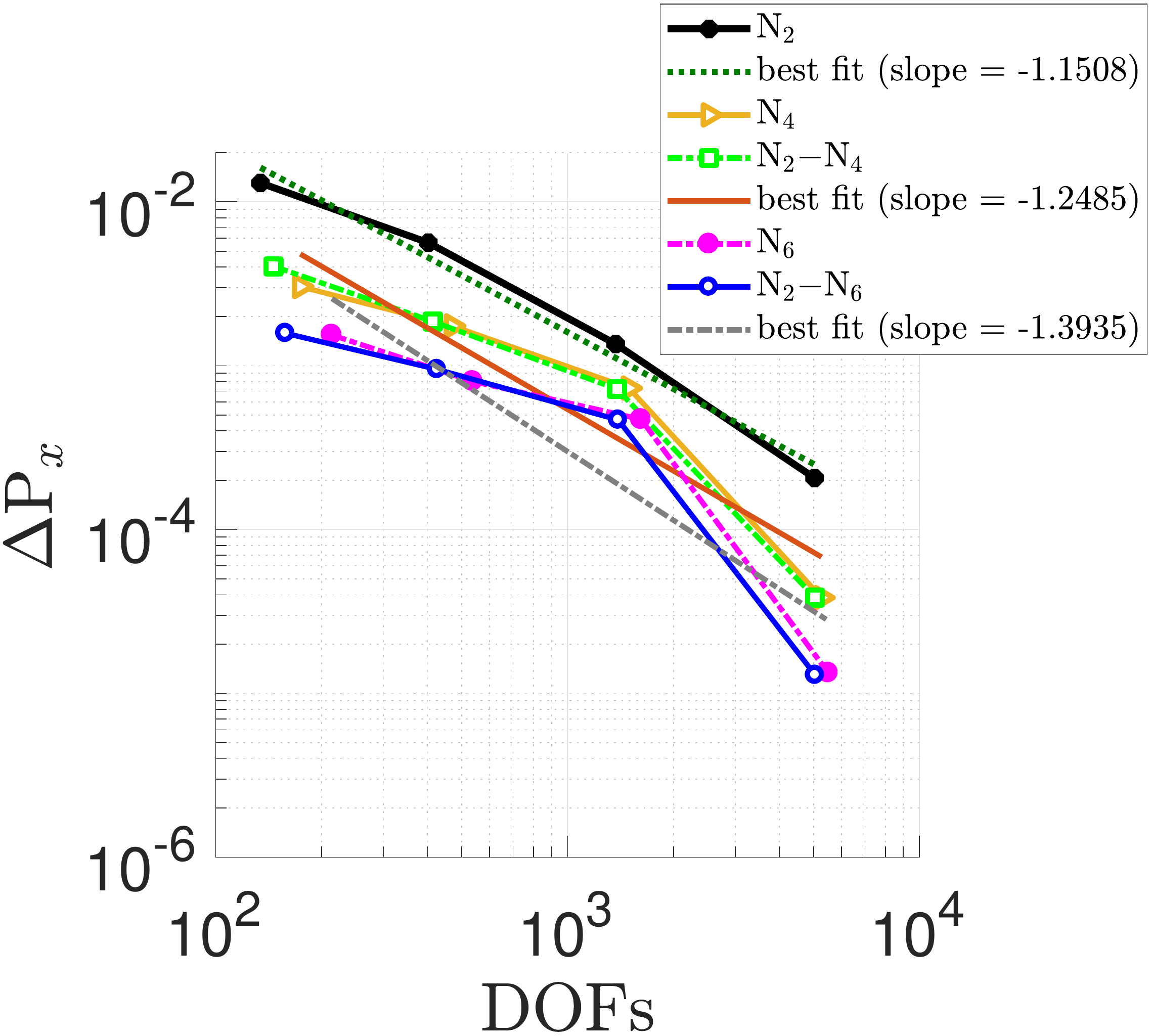}\label{fig:dPx_continuity}} 
	\subfloat[]{\includegraphics[width=.52\linewidth]{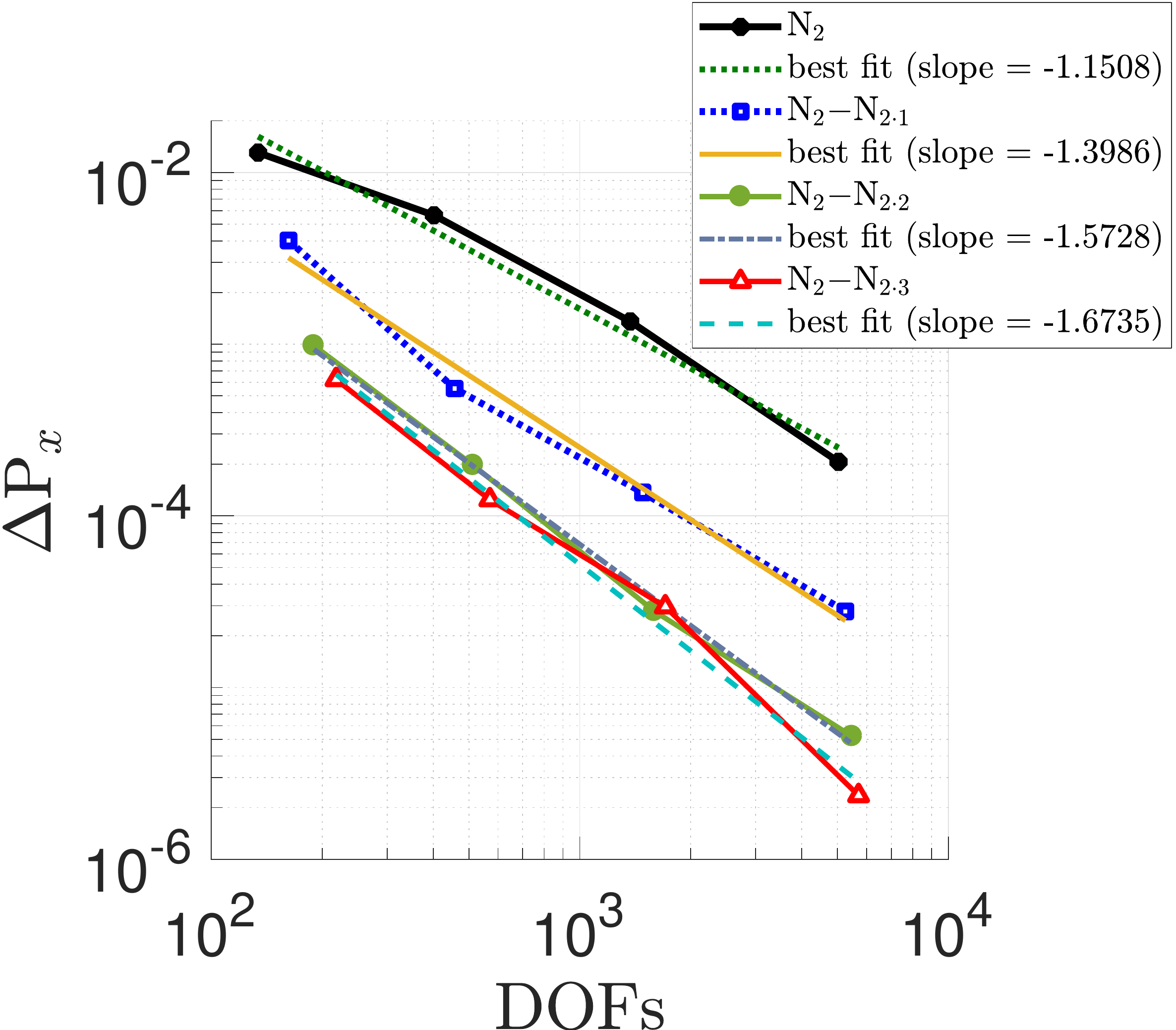} \label{fig:dPx_order}}
	\caption{The oscillation amplitude of the vertical contact reaction force $ \Delta $P$ _y $ (top) and the horizontal contact force $ \Delta $P$ _x $ (bottom) for different VO and standard NURBS based discretizations.}\label{fig:conv_plots}
\end{figure}

\begin{figure}[!ht]
	\centering
	\subfloat[]{\includegraphics[width=.42\linewidth]{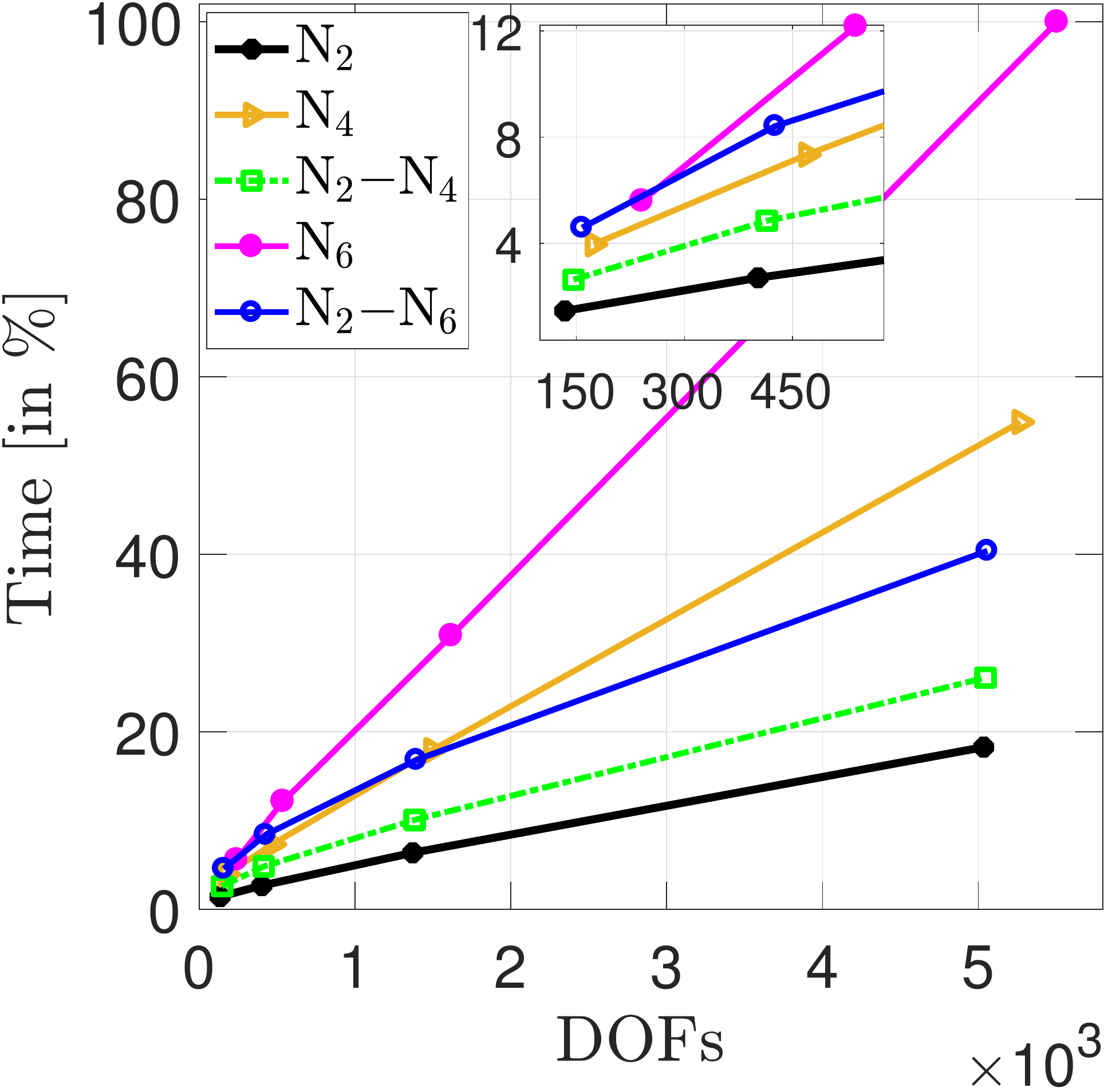} 
		\label{fig:Time_DOFs_Np}}  ~~~
	\subfloat[]{\includegraphics[width=.42\linewidth]{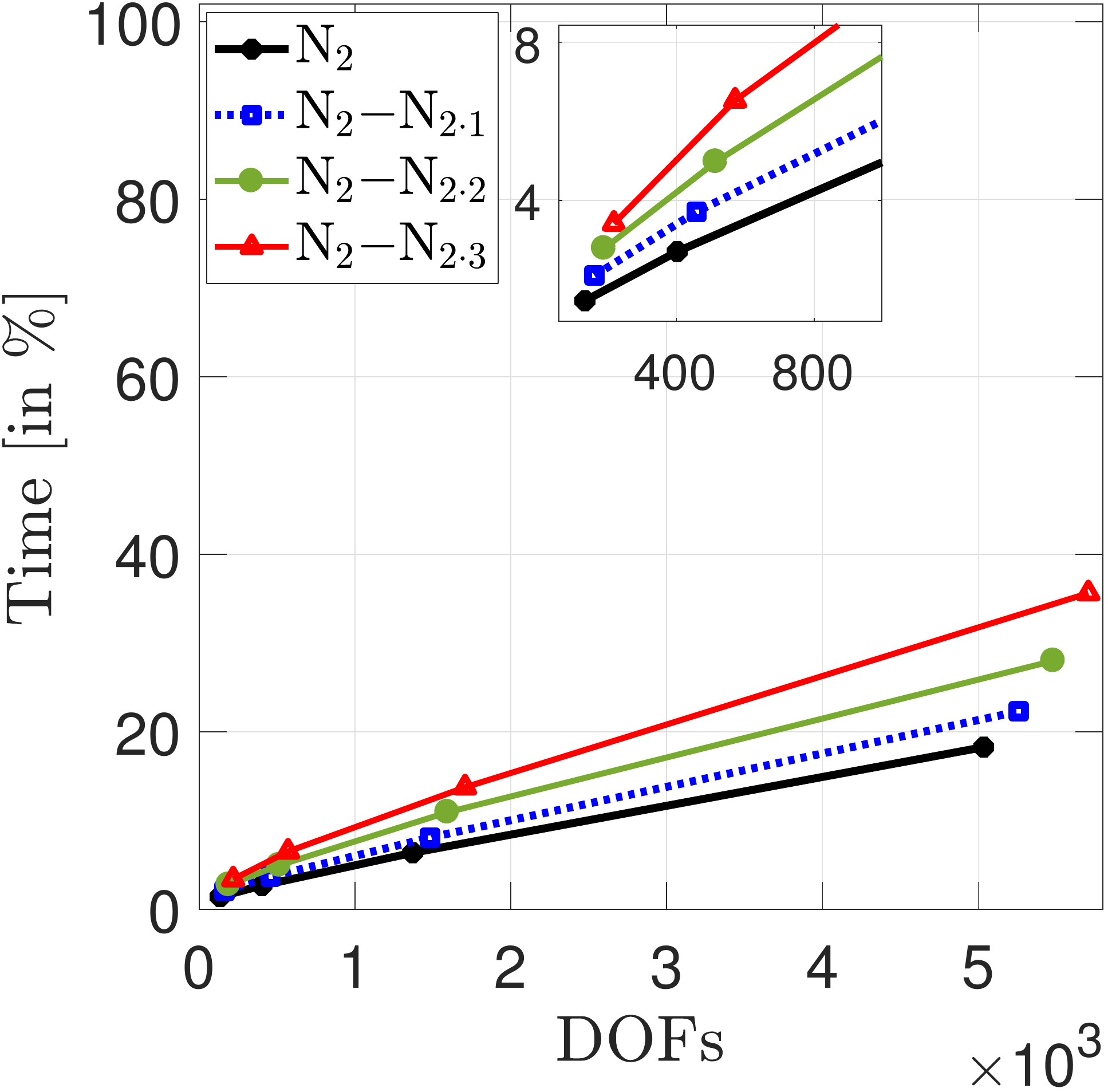} \label{fig:Time_DOFs_N2p}} 
	\caption{Total analysis time (in $ \% $) with different standard and VO based NURBS discretizations with meshes m$ _1 $ to m$ _4 $. The analysis time with N$ _6 $ using mesh m$ _4 $ is used for normalization.}  \label{fig:Time_DOFs}
\end{figure}

Figure~\ref{fig:conv_behaviour} shows the reduction in the oscillation amplitude of horizontal contact force curve P$ _x $ for different discretizations. The convergence plots for both $ \Delta \textrm{P}_y $ and $ \Delta \textrm{P}_x $ with N$ _2- $N$ _{p_c}~(p_c= 4,\,6,\, 2\cdot1,\, 2\cdot 2, \textrm{ and } 2\cdot 3) $ and standard N$ _p~(p=2,\, 4,\, \textrm{6}) $ based discretizations for all the four meshes are shown in Fig.~\ref{fig:conv_plots}. The total processing time taken by different discretizations over the DOFs are shown in Fig.~\ref{fig:Time_DOFs}. For each discretization, the time in $ \% $ is calculated using the following expression
\begin{equation}\label{eq:analysis_time}
\text{Time percentage} = \frac{\text{total analysis time}}{\text{maximum total analysis time}} \times 100\,,
\end{equation}
where the maximum total analysis time is with the standard N$_6 $ order of NURBS using the finest mesh m$ _4 $ (which is the most expensive).

From Figs~\ref{fig:dPy_continuity}, \ref{fig:dPx_continuity}, and~\ref{fig:Time_DOFs_Np} it can be observed that to attain the accuracy equivalent to N$ _4 $, N$ _2-$N$ _4 $ takes approximately $ 33.74 \% $ lower computational cost even at the coarsest mesh m$ _1 $. The reduction in the computational cost with N$ _2-$N$ _4 $ is due to fact that it employs the lower order of NURBS, i.e. N$ _2 $, for the bulk computations as compared to N$ _4 $ based uniform discretization. The gain in the computational efficiency further improves on increasing the mesh resolution, e.g. N$ _2- $N$ _4 $ takes approximately $ 52.44 \% $ lower cost as compared to N$ _4 $ at the finest mesh m$ _4 $. A similar observation is made on comparing the results with more higher-continuous NURBS discretizations, e.g. to deliver the result equivalent to N$ _6 $, N$ _2- $N$ _6 $ takes approximately $ 59.58\% $ lower analysis time at mesh m$ _4 $.

Further, from Figs.~\ref{fig:conv_plots} and~\ref{fig:Time_DOFs} it can be observed that for a fixed mesh, a significant improvement in the accuracy is achieved with additional order-elevation based VO discretizations, i.e. with N$ _2- $N$ _{2\cdot 1} $ and N$ _2- $N$ _{2\cdot 2} $, among all tested cases. The N$ _2- $N$ _{2\cdot 2} $ delivers the most accurate results at a cost slightly more than with the N$ _2 $. Taking a closer look at Figs.~\ref{fig:dPy_order} and~\ref{fig:dPx_order} reveals that the accuracy attained with N$ _2- $N$ _{2\cdot 2} $ at mesh m$ _1 $ and m$ _2 $ is comparable to N$ _2 $ using mesh m$_3 $ and m$ _4 $, respectively. The corresponding total analysis time taken by N$ _2- $N$ _{2\cdot 2} $ is approximately $ 56.62 \% $ and $ 72.78 \% $ lower than with the N$ _2 $ for the similar accuracy level, see Fig.~\ref{fig:Time_DOFs_N2p}. With other higher order based VO discretizations, i.e. with N$ _2- $N$ _{2\cdot3} $, a slight improvement in the accuracy is achieved as compared to N$ _2- $N$ _{2\cdot2} $ for the same mesh level.

In summary, this example demonstrates the advantages of the proposed VO based discretizations in terms of the accuracy of results and computational efficiency over the standard fixed-order based NURBS discretizations. To attain the accuracy equivalent to N$ _p $ based discretization, N$ _2- $N$ _p $ takes much lower computational cost. Further, with the additional order-elevation based VO discretizations, a significant improvement in the accuracy is achieved even at the coarse mesh. It attains a major gain in the numerical efficiency to provide the accuracy as similar to N$ _2 $ based discretization.

\subsubsection{Performance at large indentation depth}
\begin{figure}[!t]
	\centering
	\subfloat{\includegraphics[scale=0.31]{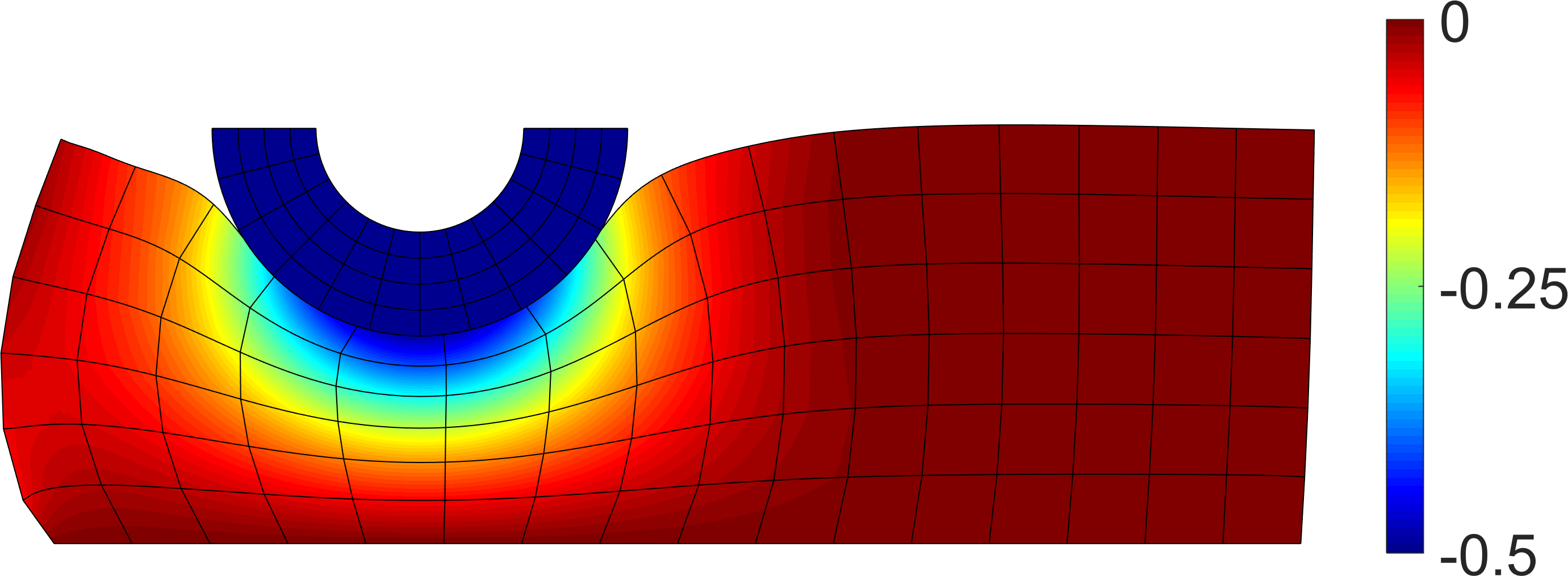}\label{fig:ironingL_disp1}} ~~~~
	\subfloat{\includegraphics[scale=0.31]{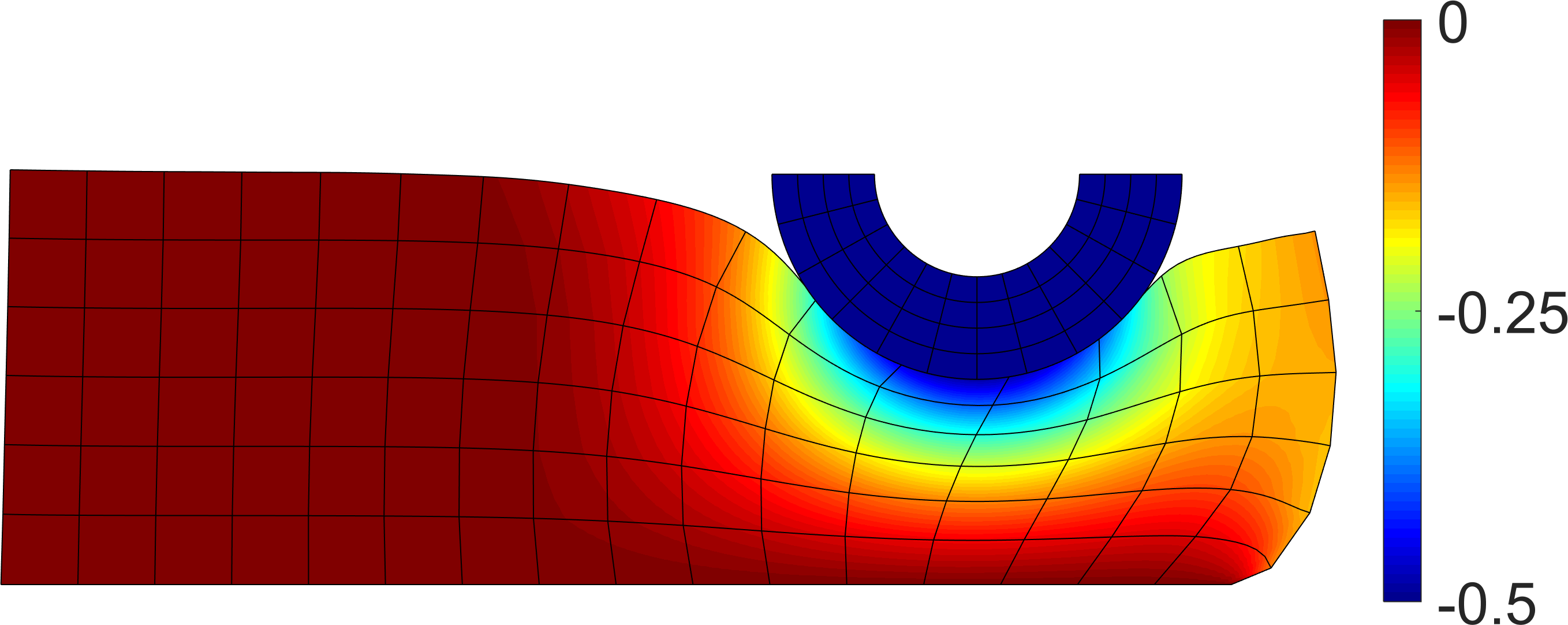}\label{fig:ironingL_disp2}} \\	
	\caption{Ironing problem with large indentation depth: Deformed configuration at the end of compression process (left) and at the end of sliding (right) with N$ _2 $ at mesh m$ _2 $. The color shows the magnitude of the vertical displacement $ u_y $. Here, $ U_y =-0.5 $ and $ E_{\textrm{die}} = 10^4 $.} \label{fig:ironingL_disp}
\end{figure}

Finally, in this section, we examine the influence of the interpolation order of the underlying bulk elements on the solution quality for the ironing problem with large indentation depth. It is known that numerical simulation of this problem can be computationally challenging if the stiffness ratio of the bodies is high and the harder body is dragged relative to the softer one after the significant indentation~\cite{PUSO2004, Temizer2014}. In such a case, the strong interaction between the bulk and contact surface discretizations takes place. Therefore, it becomes important to assess the accuracy of the solution with VO discretization to that of corresponding standard NURBS discretization for large deformation contact. The problem setup used for this analysis is same as described in Sec~\ref{sec:ironing_setup} except the modification that higher stiffness for the die $ E_{\textrm{die}} = 10^4 $ N/mm$ ^2 $ and vertical displacement $ U_y = -0.5 $ mm are considered. The Coulomb's friction coefficient is $ \mu_{\textrm{f}} = 0.1 $. Four nested meshes that are listed in Table~\ref{table:mesh_size} are used for the analysis. The deformed configurations of the setup at the end of compression and sliding processes are shown in Fig.~\ref{fig:ironingL_disp} with N$ _2 $ using mesh m$ _2 $.	
\begin{figure}[!b]
	\centering
	\subfloat{\includegraphics[width=.45\linewidth]{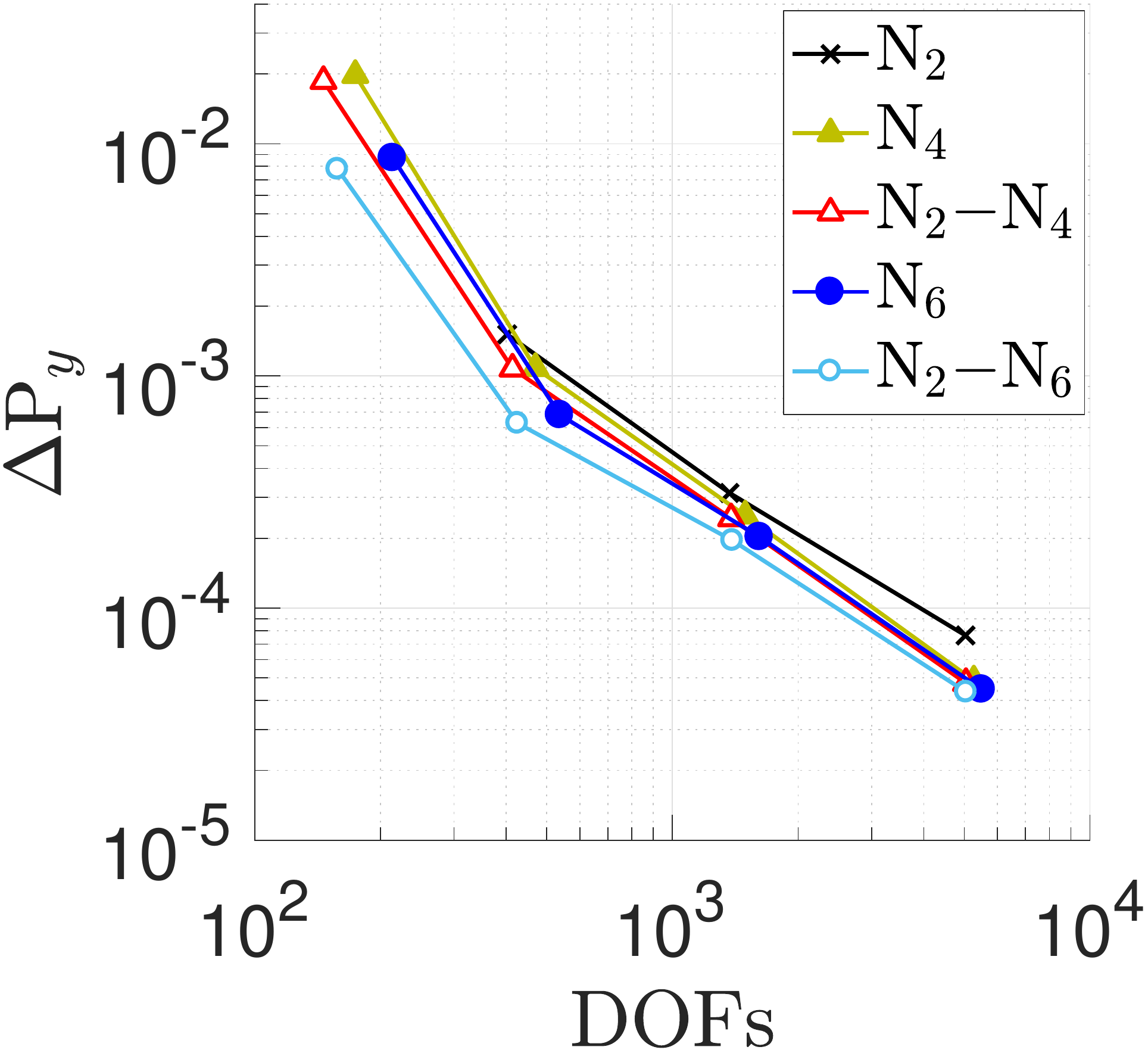}\label{fig:dPy_NpL}} ~~~~
	\subfloat{\includegraphics[width=.45\linewidth]{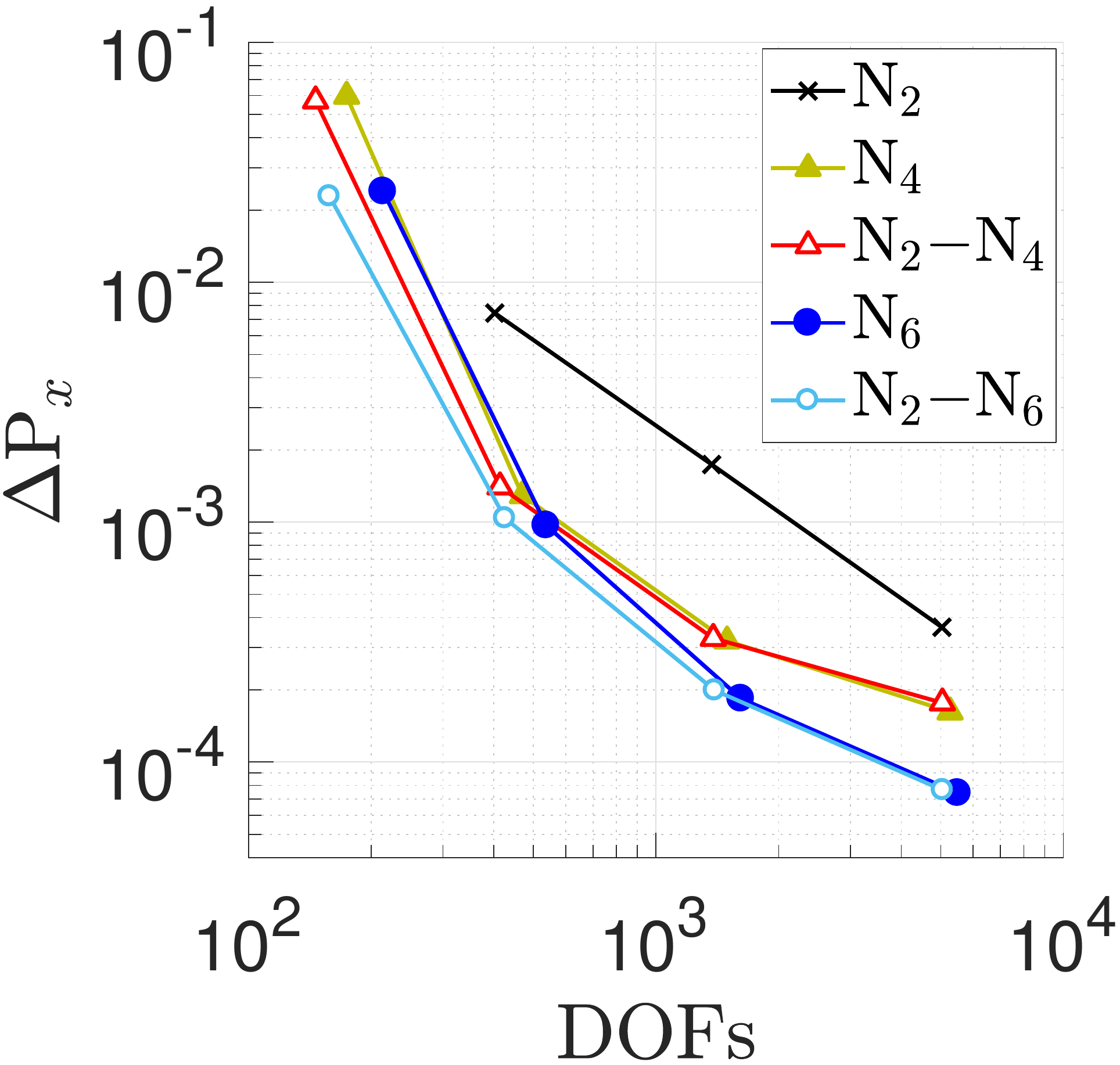} \label{fig:dPx_NpL}} \\
	\caption{Ironing with large indentation depth: Comparison of the oscillation amplitude of contact forces with different standard and VO based discretizations for meshes m$ _1 $ to m$ _4 $.}  \label{fig:conv_plotsL}	
\end{figure}

Figure~\ref{fig:conv_plotsL} shows the comparison of the oscillation amplitude of the vertical contact force $ \Delta \textrm{P}_y $ and horizontal contact force $ \Delta $P$ _x $ with both the discretization methods for meshes m$ _1 $ to m$ _4 $. Two major observations are made. First, although N$ _2 $ based discretization is  $ C^1 -$continuous across the contact layer, it fails to converge at the coarsest mesh m$ _1 $. Clearly, N$ _2 $ with the mesh m$ _1 $ is not sufficient to analyze this problem. On the other hand, the higher-order based NURBS discretization is capable of providing the solution even at the coarsest mesh m$ _1 $. Second, the results obtained with the VO based discretization verify the application of the proposed methodology to large deformation contact case. For the same mesh level, the accuracy of the result achieved with VO based N$ _2- $N$ _{p_c}~(p_c = 4 \textrm{ and } 6) $ discretizations is equivalent to that of fixed-order based N$ _p~(p=4 \textrm{ and } 6) $ discretizations, respectively.

\subsection{Hertzian contact problem}
In the context of isogeometric contact analysis, Hertzian-type contact problem has been widely used to asses the contact pressure distribution across the contact interface on increasing the mesh resolution and interpolation order of the NURBS with the different type of contact formulations, see e.g.~\cite{Temizer2011, DeLorenzis2011, Temizer2012, DeLorenzis2012, Kim2012, Matzen2013, Matzen2016, Dimitri2014Tspline, SEITZ2016}. Therefore, as a second example, we consider it to demonstrate the coarse mesh accuracy of the VO NURBS discretization method as compared to standard NURBS discretizations.
\begin{figure}[!h]   
	\begin{center}	
		\subfloat{\includegraphics[scale=0.3]{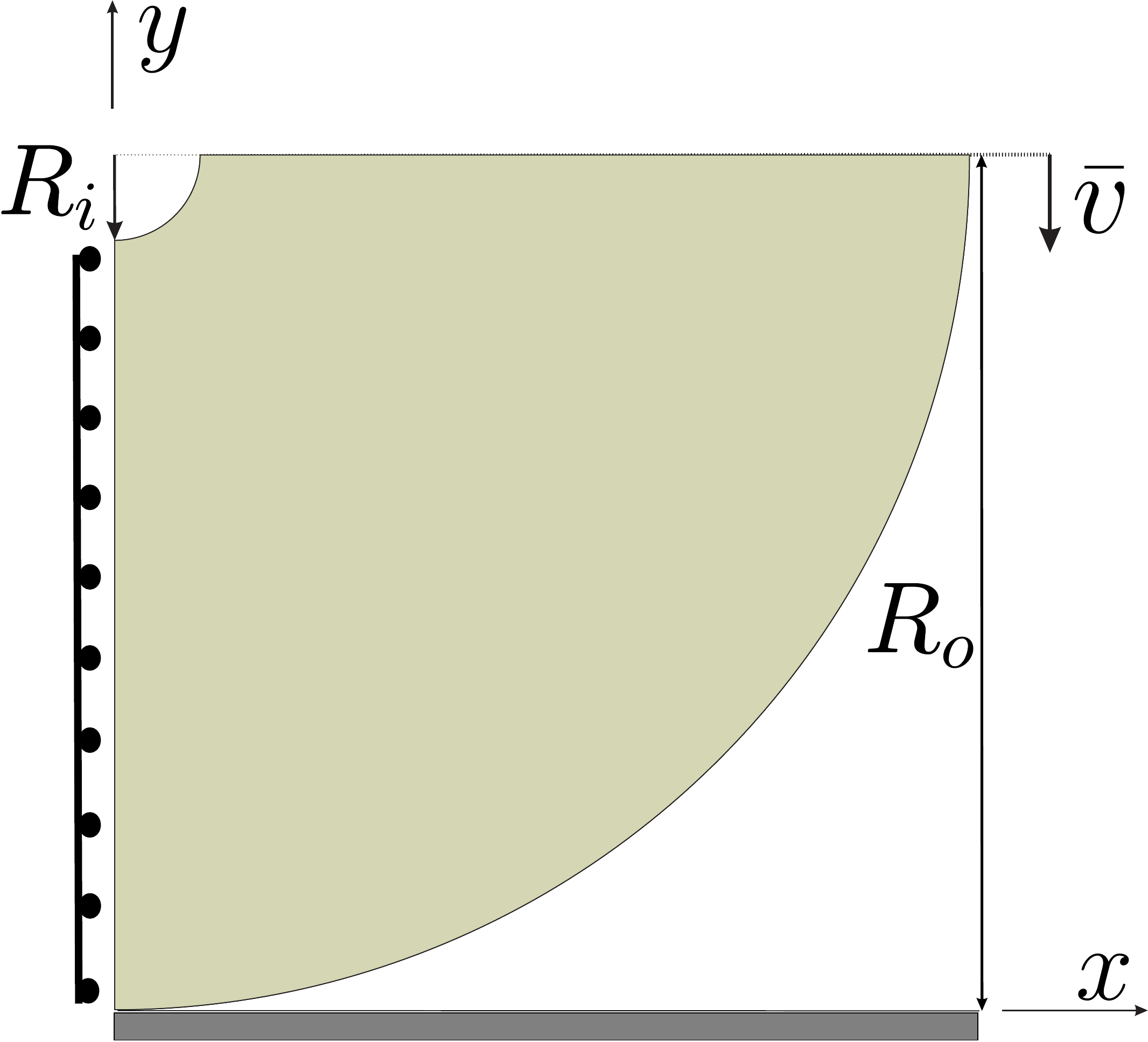}}
		\caption {Setup of the Hertzian problem.}
		\label{fig:Hertz_setupa}
	\end{center}
\end{figure}

\begin{figure}[!t]
	\centering	
	\subfloat[]{\includegraphics[width=.35\linewidth]{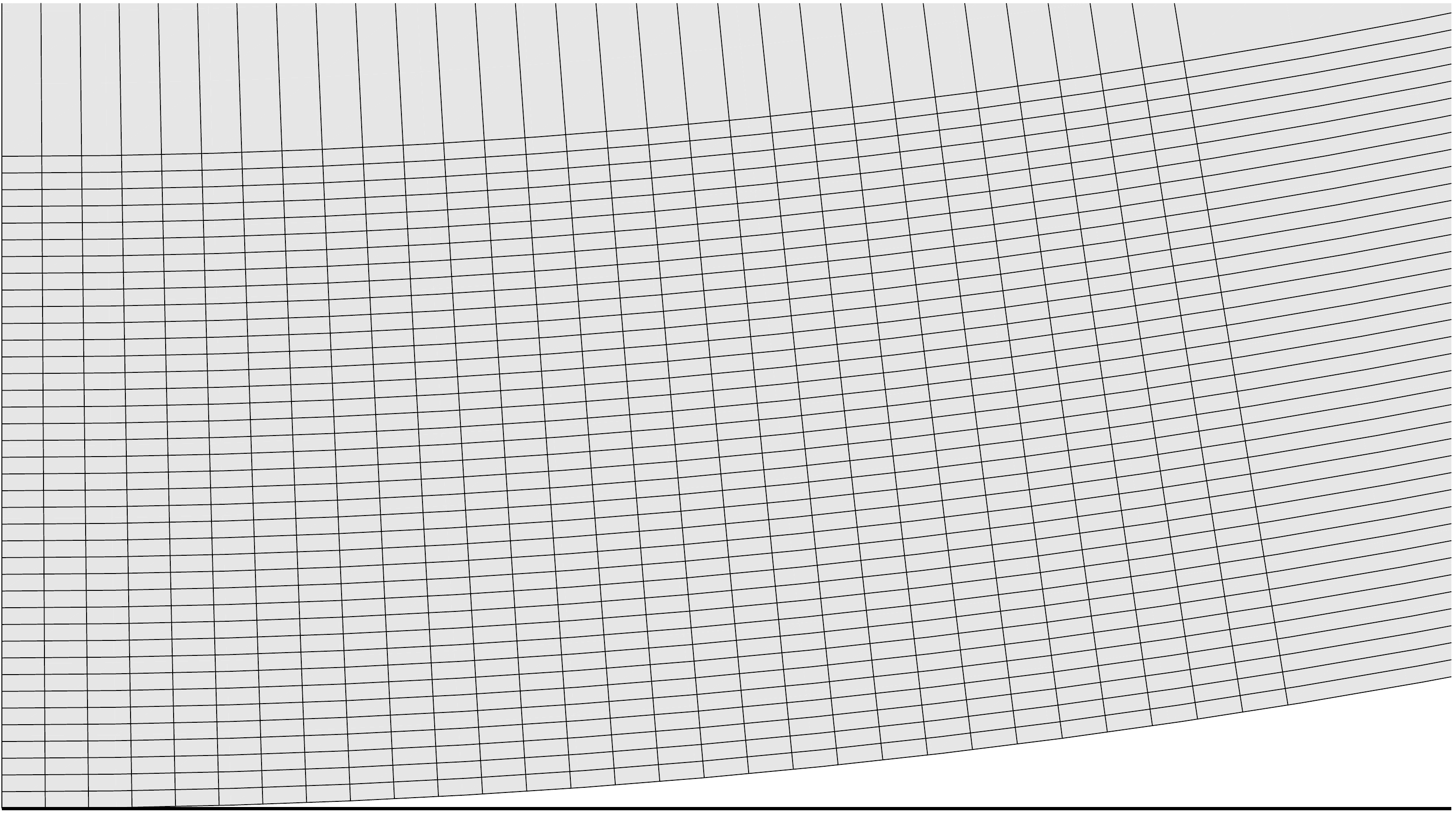}\label{fig:ref_config}} ~~~~~~~
	\subfloat[]{\includegraphics[width=.35\linewidth]{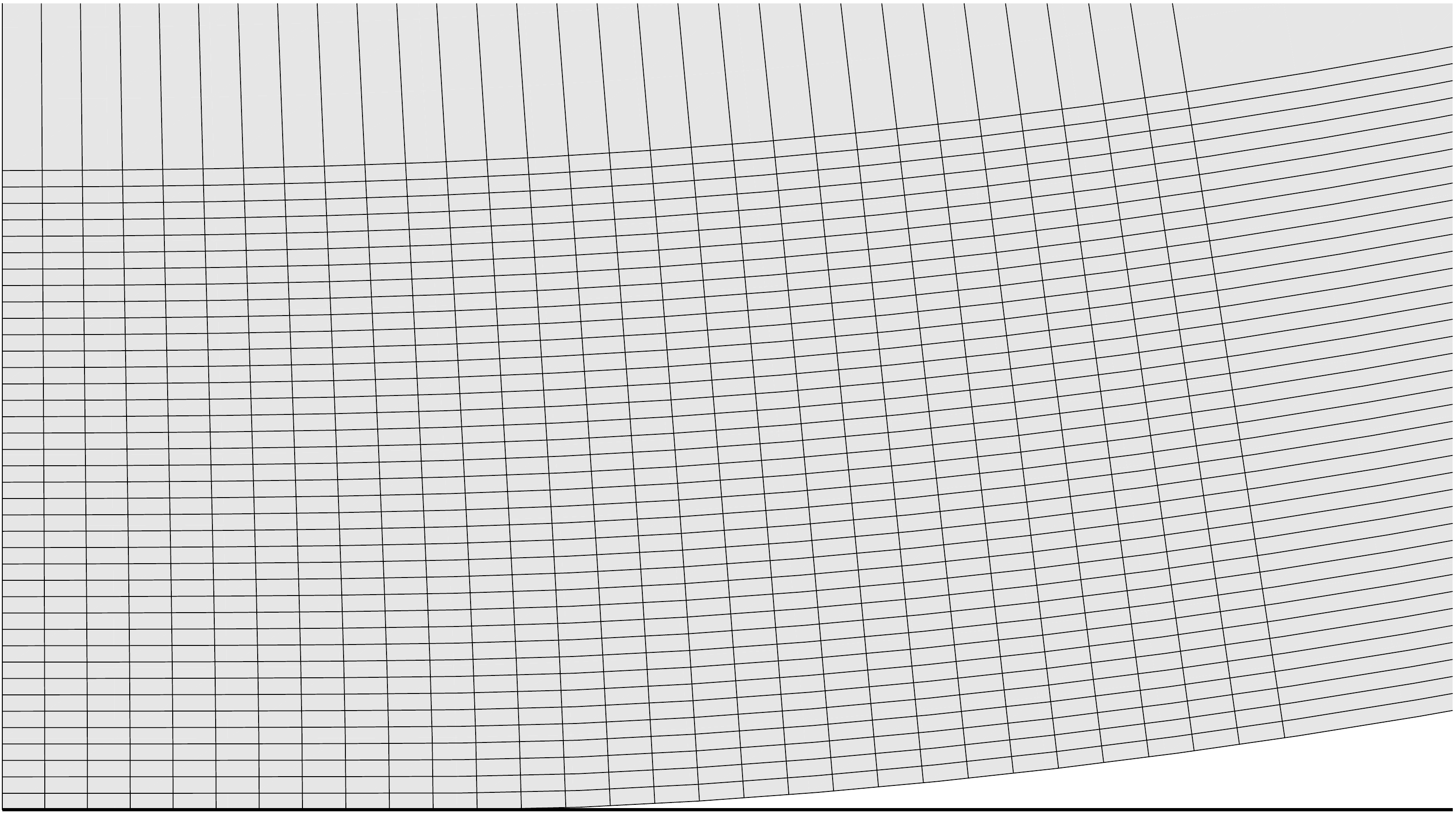}\label{fig:cur_config}}
	\caption{(a) The undeformed and (b) deformed configuration of the setup with N$ _2 $ using mesh m$ _3 $ after performing the mesh relocation.}\label{fig:configurations}
\end{figure}

The example considers the contact between a deformable infinitely long cylinder, having an outer radius $ R_o = 1 $, and a rigid surface under the plane strain setting~\cite{johnson1987}. The setup of the problem along with the boundary condition, which is taken from~\cite{Temizer2011}, is shown in  Fig.~\ref{fig:Hertz_setupa}. Due to symmetry, only a quarter of the geometry is utilized. The natural NURBS based description of geometry consists a non-zero internal radius $ R_i  = 0.1$. The numerical investigation reveals that it exerts no influence on the quality of the solution. For its modeling, a linearly elastic material with Young's modulus $ E = 1 $ and Poisson's ration $ \nu = 0.3 $ is used. Here, the top surface of the cylinder is subjected to prescribed vertical displacement $ \bar{v} $ and a penalty parameter $ \epsilon_{\text{N}} = 2\times 10^3 $ is chosen as a default value. Six uniformly refined meshes that are driven by $ 9n \times 48 $ number of elements along each parametric directions $(\textrm{where } n = 1,\, 2,\, 4,\, 8,\, 16\textrm{ and } 32)$ and are denoted by m$ _1 $ to m$_6 $ are used for this analysis. For numerical efficiency, approximately $ 80 \% $ of the total elements in each parametric direction are relocated in such a manner that they lie within the $ 10 \% $ surface length of the geometry, as shown in Fig.~\ref{fig:configurations}. It has been tested that $ 48 $ number of redistributed elements along the radial direction is adequate for the considered example. The undeformed and deformed configurations of the setup with N$ _2 $ using mesh m$ _3 $ are shown in Figs.~\ref{fig:ref_config} and~\ref{fig:cur_config}, respectively.

\begin{figure}[!t]
	\centering	
	\subfloat[]{\includegraphics[width=.34\linewidth]{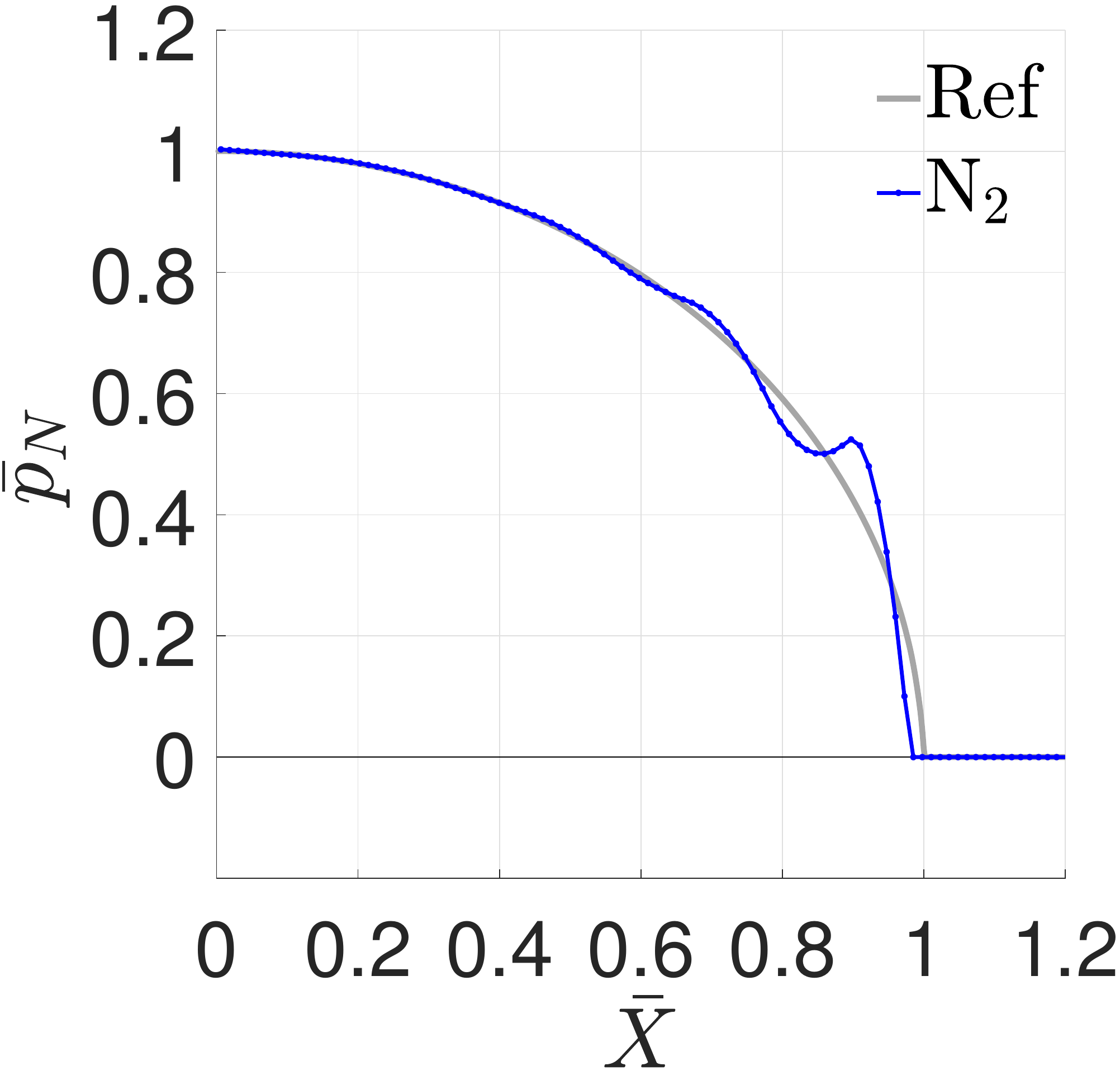}\label{fig:N2_m1}} 
	\subfloat[]{\includegraphics[width=.34\linewidth]{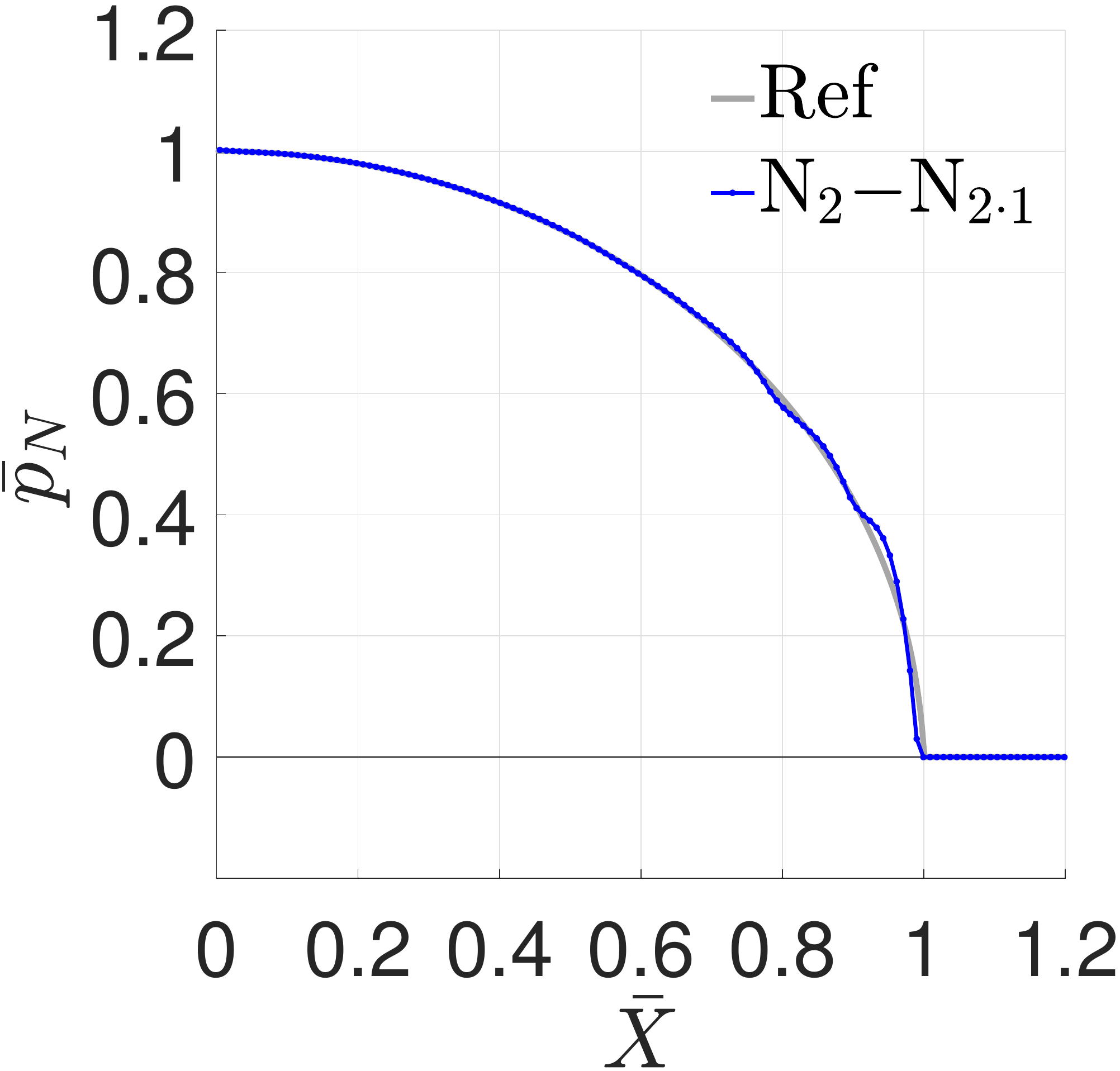}\label{fig:N21_m1}} 
	\subfloat[]{\includegraphics[width=.34\linewidth]{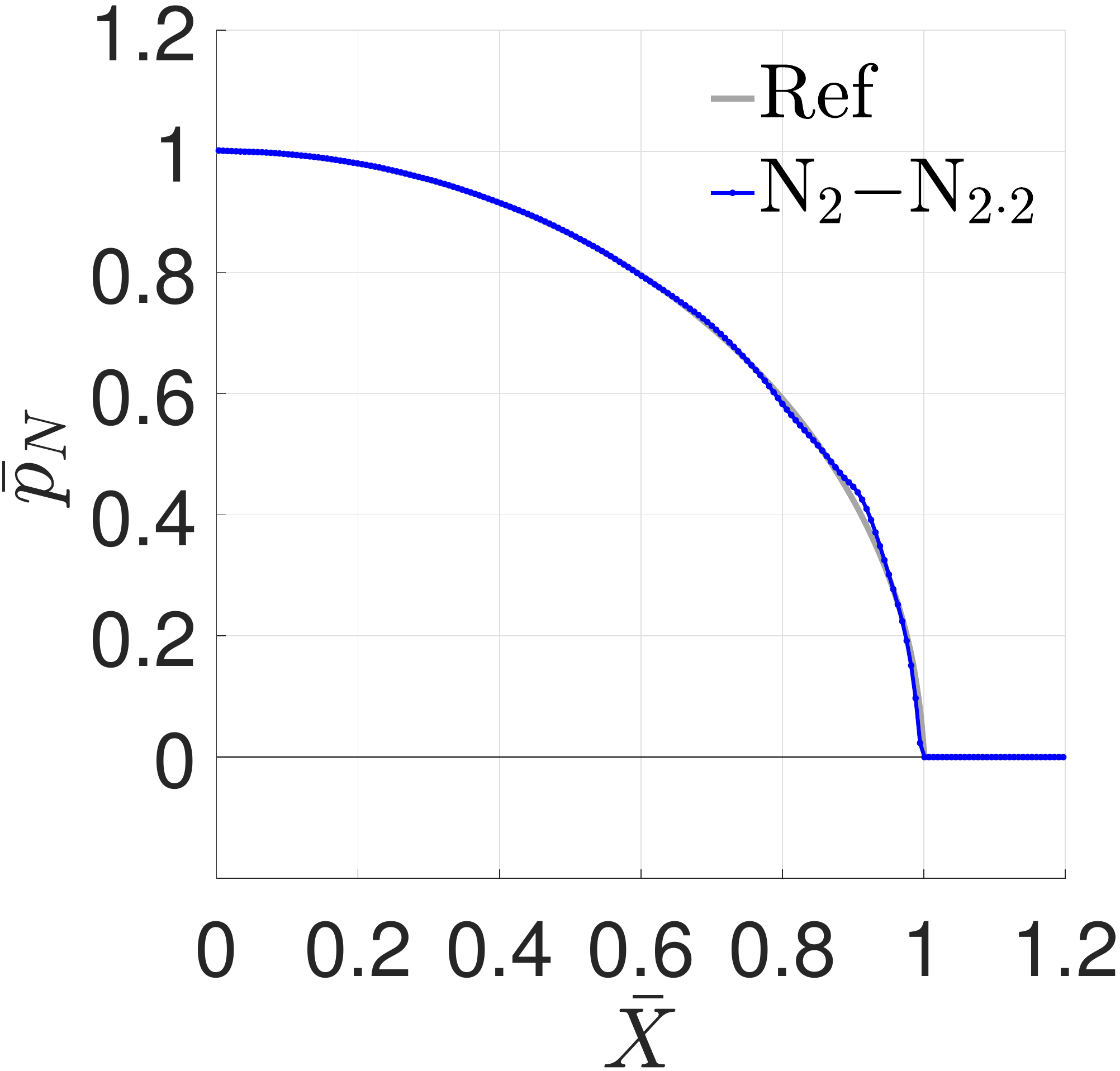}\label{fig:N22_m1}} \\
	\subfloat[]{\includegraphics[width=.34\linewidth]{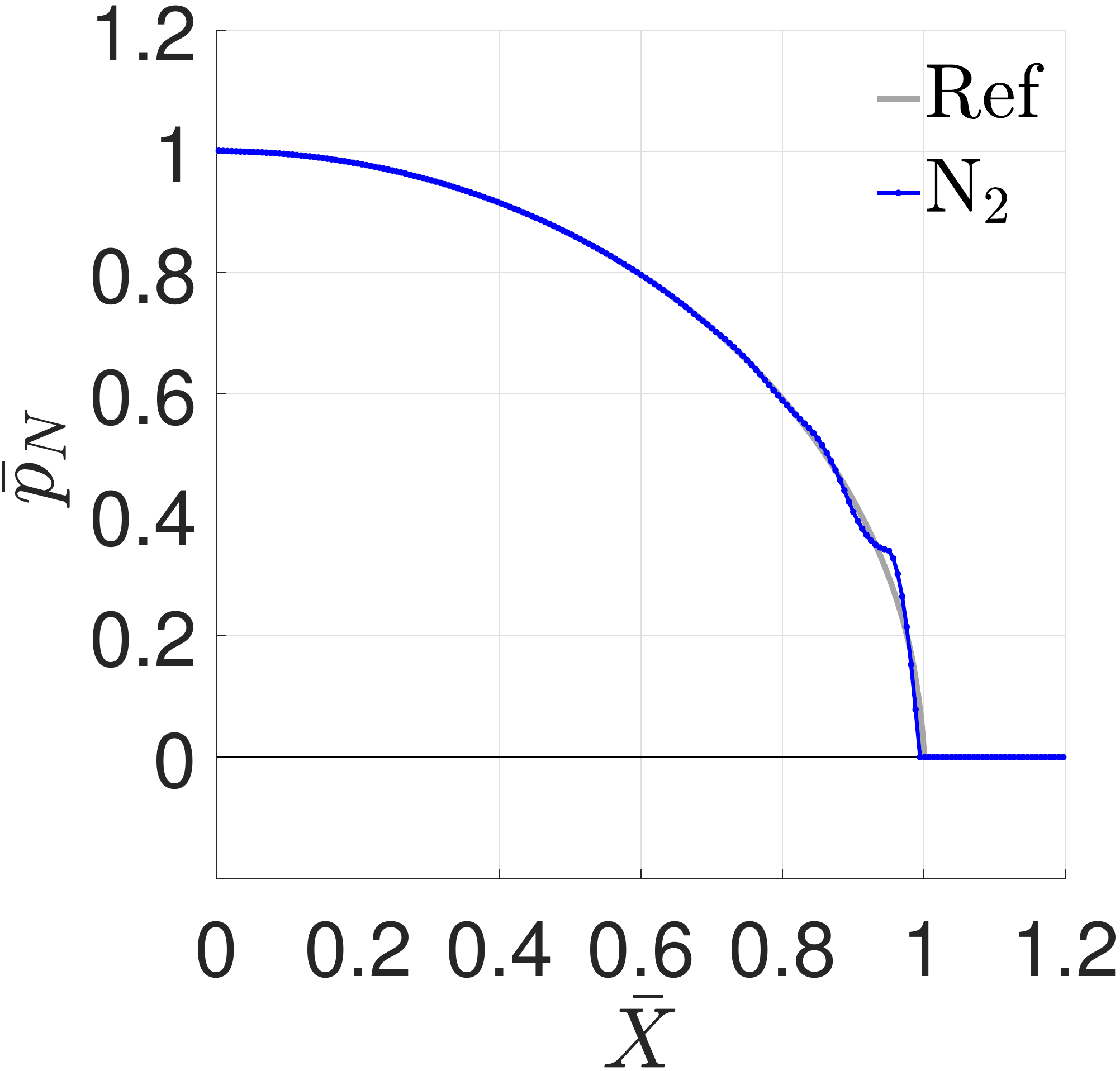}\label{fig:N2_m2}}  
	\subfloat[]{\includegraphics[width=.34\linewidth]{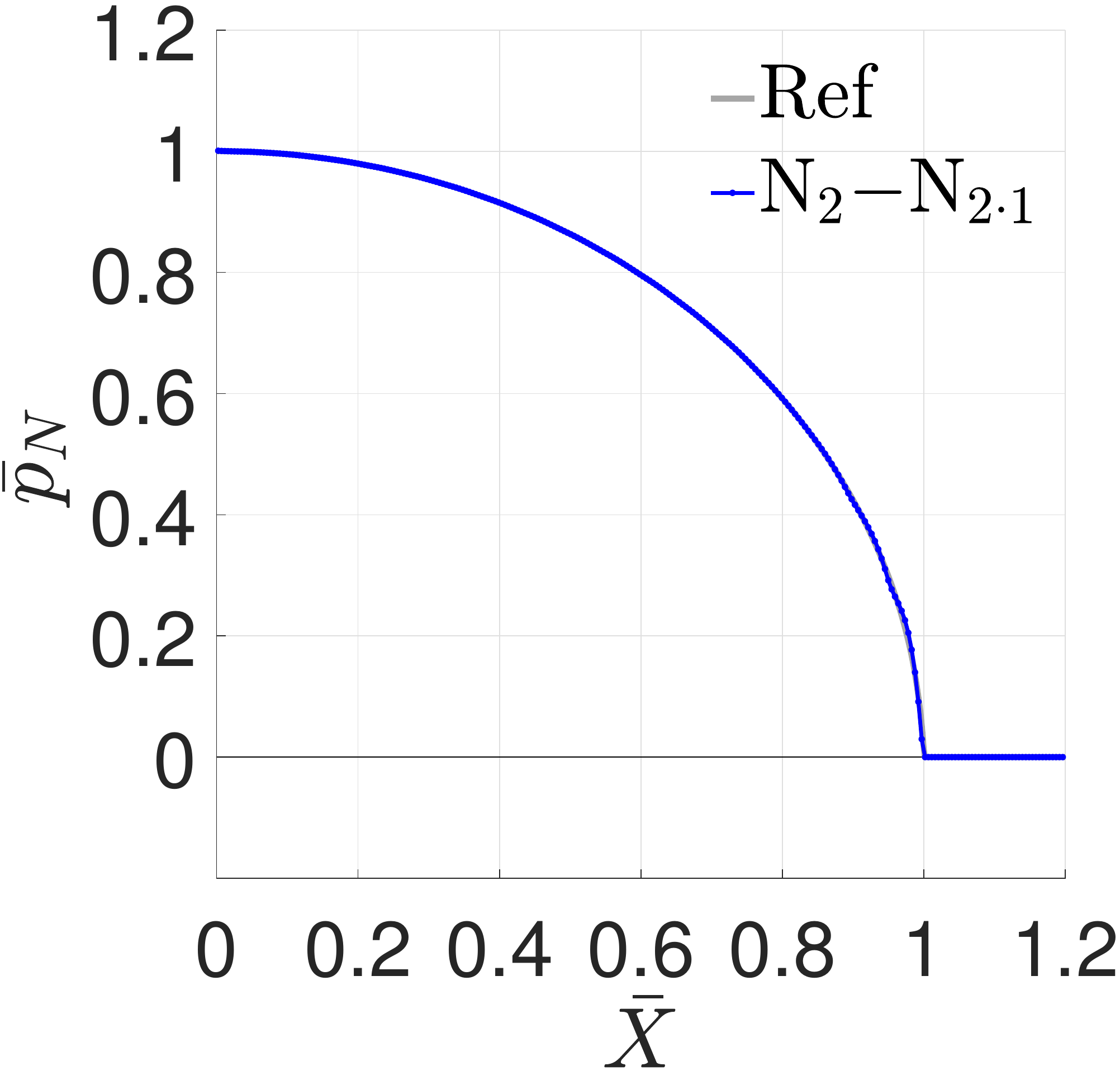}\label{fig:N21_m2}} 
	\subfloat[]{\includegraphics[width=.34\linewidth]{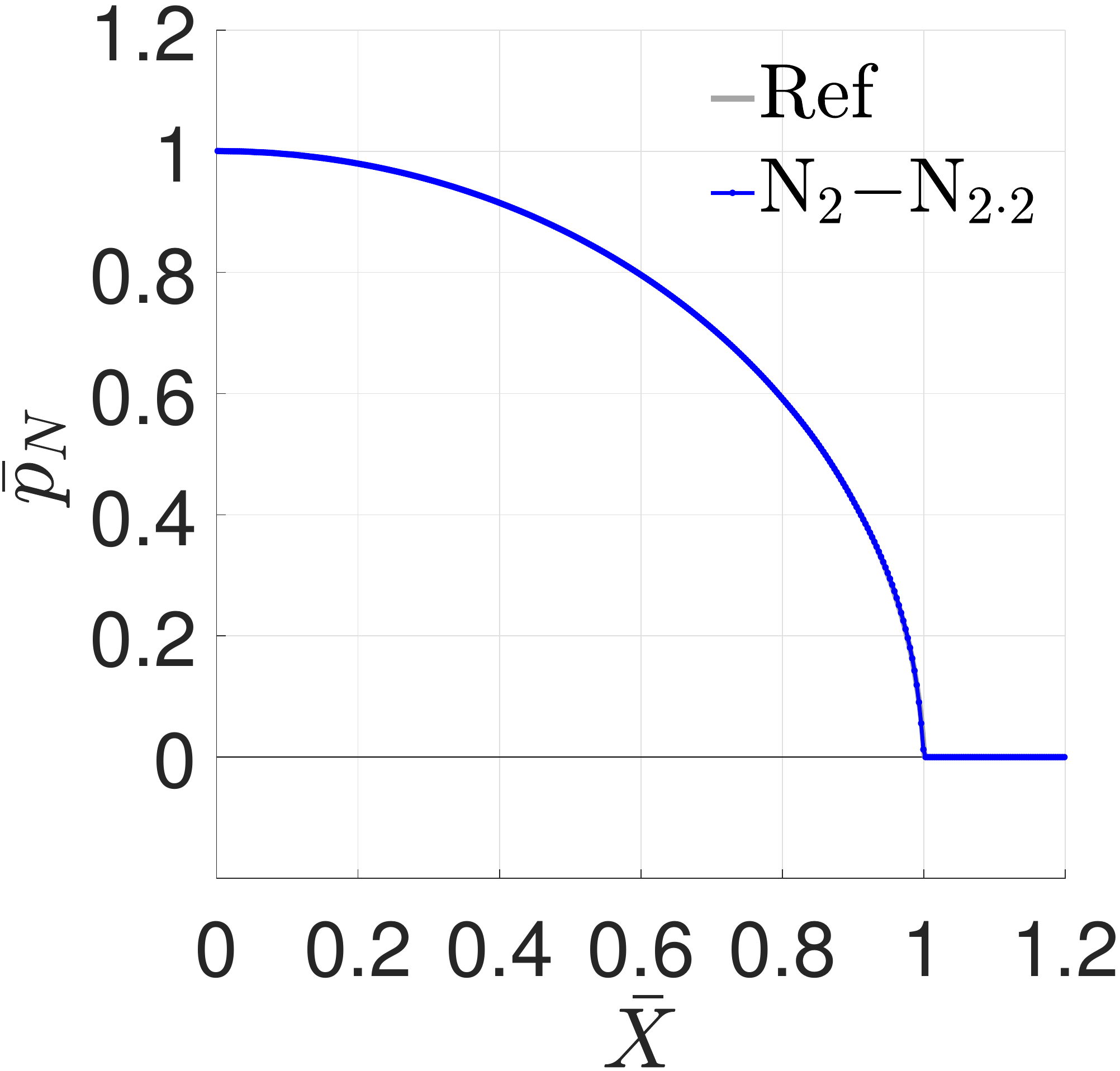}\label{fig:N22_m2}} \\
	\subfloat[]{\includegraphics[width=.34\linewidth]{Fig18f.pdf}\label{fig:N2_m4}}
	\subfloat[]{\includegraphics[width=.34\linewidth]{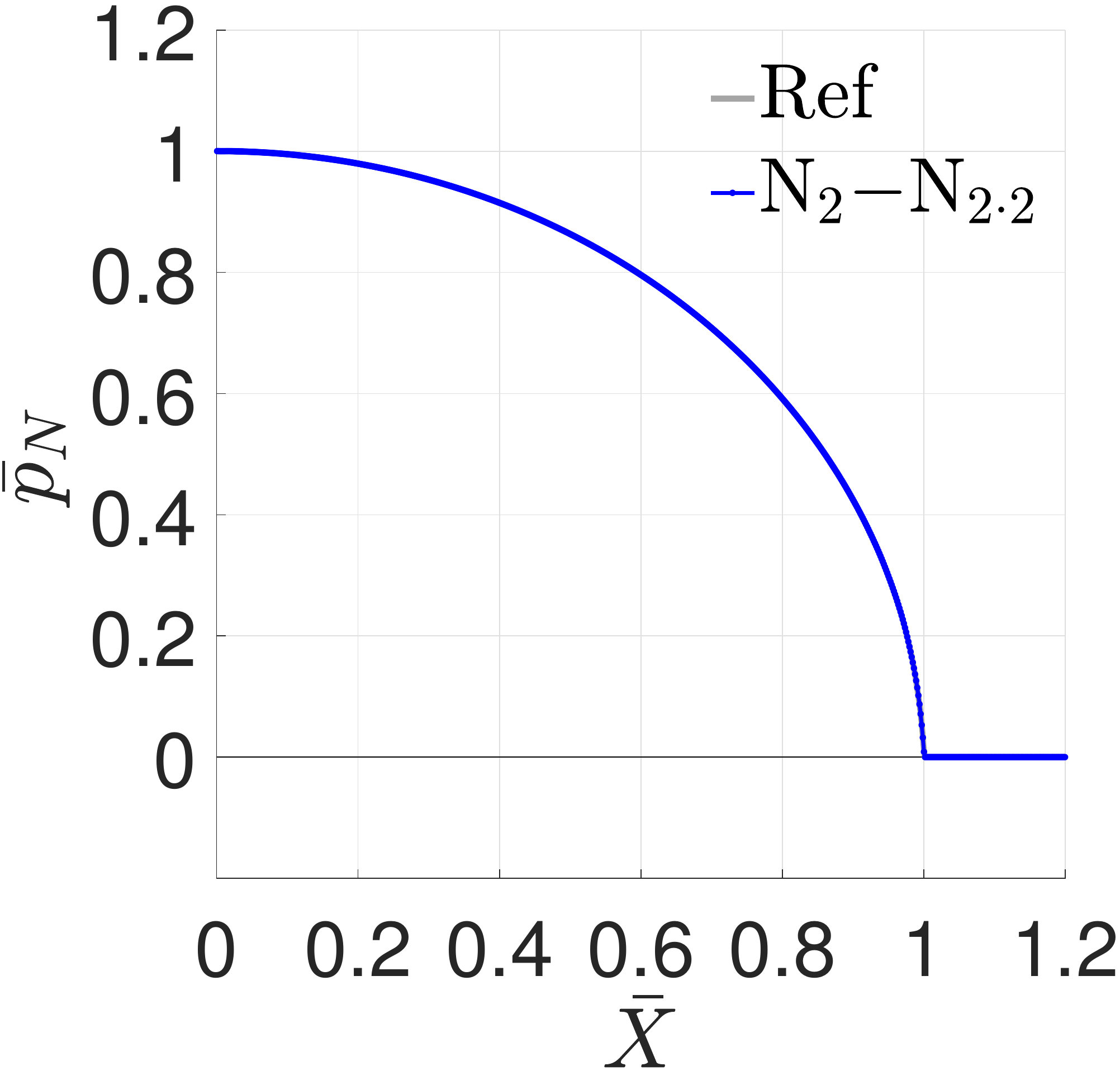}\label{fig:N21_m4}}
	\subfloat[]{\includegraphics[width=.34\linewidth]{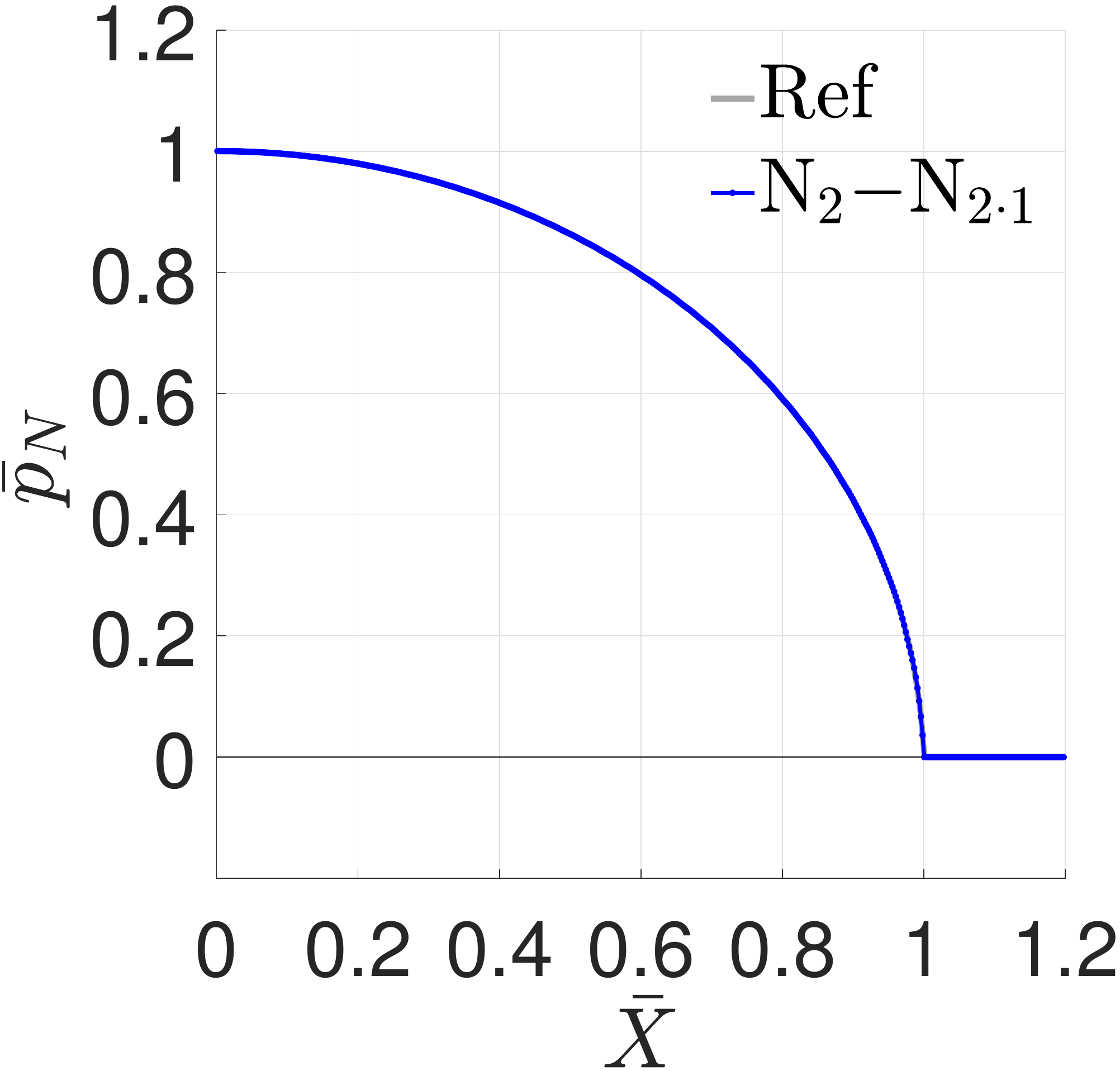}\label{fig:N22_m4}}
	\caption{Comparison of the variation of contact pressure distribution with standard and VO based NURBS discretizations at: (a)-(c) $ 36\times 48 $; (d)-(f) $ 72\times 48 $; (g)-(i) $ 144\times 48 $ mesh arrangements. The result with N$ _2 $ using finest mesh m$ _6 $ is used as a reference.} \label{fig:pressure_distribution}
\end{figure}

Figure\ref{fig:pressure_distribution} shows the comparison between the local solution quality of contact pressure distribution for both the standard N$ _2 $, and VO based N$ _{2\cdot 1} $ and N$ _{2\cdot 2} $ discretizations with that of reference contact pressure distribution for meshes m$ _3 $ to m$ _5 $. Here, the dimensionless normal contact pressure $ \bar{p}_N = p_N/p_{N_{ref}} $ is plotted versus the dimensionless contact coordinate $ \bar{X} = x/a_{ref} $, where $ p_N $ is the normal contact pressure evaluated at an active integration point and $ x $ is the distance of this points from the first point of contact. The result obtained with N$ _2 $ at the finest mesh level m$_6 $ is used as a reference. The corresponding maximum contact pressure $ p_{N_{ref}} $ and the contact area $ a_{ref} $ are used for the normalization purpose. It is noted that the results obtained with N$ _2 $ based discretization for different meshes are in agreement to those reported in literature~\cite{Temizer2011, DeLorenzis2011, Temizer2012, DeLorenzis2012, Kim2012, Dimitri2014Tspline, Matzen2016, SEITZ2016} for Hertzian-type contact problem. In case of VO based discretizations, the quadratic order of NURBS that are sufficient to represent the considered CAD geometry exactly are kept fixed for the bulk description. For contact layer, one and two steps of additional order elevation refinements that yields N$ _2- $N$ _{2\cdot1} $ and N$ _2- $N$ _{2\cdot 2} $ discretizations, respectively, are performed. The total number of DOFs present in the contact interface and in the remaining bulk region with both the standard and VO based NURBS discretizations for each mesh are listed in Table~\ref{table:DOFs}.

From Fig~\ref{fig:pressure_distribution}, it can be observed that for the same mesh resolution, VO based discretizations provide more accurate results than with the standard N$ _2 $ based discretization. With N$ _2 $, the oscillation of the contact pressure near the contact interface boundary is clearly visible at the coarse mesh m$_3 $, see Fig~\ref{fig:N21_m1}. As the mesh is refined, the magnitude and extent of the oscillation reduces, as can be seen from Figs.~\ref{fig:N2_m1}, \ref{fig:N2_m2}, and~\ref{fig:N2_m4}. With a very fine mesh m$ _5 $, although the contact pressure curve seems to match exactly with the reference solution, a slight kink is still present near the contact interface boundary. It is noted that in case of Hertzian-type contact, the quality of the solution is often affected by the sudden change of contact status from active to non-active within an element that lies across the edge of contact region. To this date, a number of solution approaches have been introduced to fix this issue, see e.g~\cite{Franke2010, Kim2012}. Moreover, T-splines based discretization that allows adaptive mesh refinement delivers higher accurate solution than standard NURBS for a fixed number of DOFs, see~\cite{Dimitri2014Tspline}. However, in the present work, we particularly focus on alleviating such an error with simpler NURBS based discretizations that exclude the remeshing strategies and is computationally inexpensive.

From Figs.~\ref{fig:N2_m1}-\ref{fig:N22_m1} it can be observed that with VO based discretizations, especially with N$ _2- $N$ _{2\cdot 2} $, superior quality result is obtained even at the coarse mesh m$ _3 $ as compared to N$ _2 $. It is attributed to the additional number of DOFs present across the contact interface with N$ _{2\cdot 2} $ over N$ _2 $. This, improves the resolution of edge region at the fixed mesh. At the intermediate mesh level m$ _4 $, the contact pressure distribution with N$ _2- $N$ _{2\cdot 2} $ seems to match very well with the reference solution, see Fig.~\ref{fig:N22_m2}. This is impressive as N$_2- $N$ _{2\cdot 2} $ takes only half of the elements to deliver the solution comparable to N$ _2 $ at mesh m$ _5 $, see Fig.~\ref{fig:N2_m4}. Further, from Figs.~\ref{fig:N21_m2}-\ref{fig:N22_m2} and~\ref{fig:N21_m4}-\ref{fig:N22_m4} it can be seen that the quality of results remain unchanged on further increasing the mesh resolution.

\begin{figure}[!b]
	\begin{center}		
		\subfloat{\includegraphics[width = 0.52\linewidth]{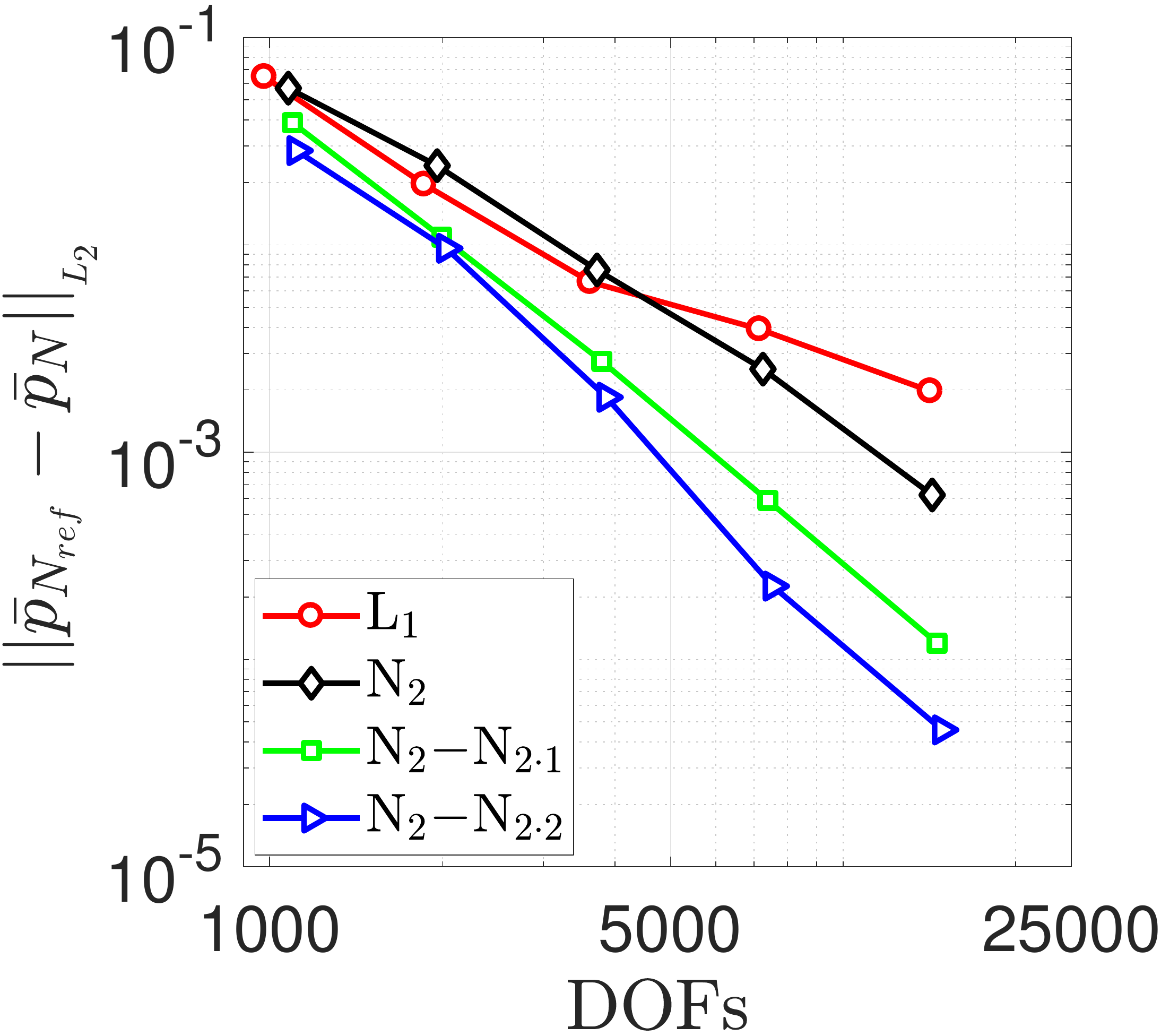}}
		\scriptsize{\caption {The $ L_2 -$norm based error in the normalized contact pressure distribution for both the discretization methods using meshes m$ _1 $ to m$ _5 $.}\label{fig:mesh_conv}}
	\end{center}
\end{figure}

Next, in Fig.~\ref{fig:mesh_conv}, we quantitatively show the comparison between the error in the normalized contact pressure distribution for both the standard and VO based discretizations for each mesh level. Moreover, the error decay over these meshes with the classical linear order of Lagrange discretization, i.e. L$ _1 $, is also included for comparison purpose. The four noded element mesh describing the circular geometry in an approximate form is generated through the conversion of NURBS elements for all mesh level, i.e. m$ _1 $ to m$ _5 $. For each mesh, the error is calculated using the $ L_2 $-norm 
\begin{equation}\label{eq:L2_norm}
{L}_2 = \sqrt{ \int_{\Gamma_c} \left[\bar{p}_{N_{ref}} - \bar{p}_{N}\right]^2\,\textrm{d}\Gamma}
\end{equation}
where $ \bar{p}_{N_{ref}} $ is the normalized reference normal contact pressure distribution that is computed with N$ _2 $ at mesh m$ _6 $. In the current work, the spatial rate of convergence is not monitored as the number of elements along the radial direction is fixed during the mesh refinement (driven by $ 9n \times 48 $). Here, we are rather interested in the difference between the accuracy of the result with VO and standard NURBS based discretizations on increasing the mesh resolution along the contact boundary layer of the geometry. It is highlighted that for the unilateral contact problems, the spatial convergence rate is generally limited by the reduced regularity of the solution under uniform mesh refinement, see e.g.~\cite{SEITZ2016} for the detailed description on the convergence rate with a different order of standard NURBS discretizations using the dual mortar based contact algorithm.%

Figure~\ref{fig:mesh_conv} clearly shows the advantage of using VO as compared to standard N$ _2 $ and L$ _1 $ based discretizations for the same mesh level. The error curve with N$ _2- $N$ _{2\cdot 2} $ lies below the N$ _2- $N$ _{2\cdot 1} $, N$ _2 $, and L$ _1 $ curves. This is attributed to additional number of DOFs present across contact interface with N$ _2- $N$ _{2\cdot 2} $, see Table~\ref{table:DOFs} for DOFs density data. A closer look reveals that accuracy in the result with N$ _2- $N$ _{2\cdot 2} $ at mesh m$ _3 $ is comparable to N$ _2 $ using mesh m$ _4 $. Further, at mesh m$ _4 $, N$ _2- $N$ _{2\cdot 2} $ provides a much more accurate result than with N$ _2 $ using a finer mesh m$ _5 $. As compared to L{$ _1 $}, a considerable gain in accuracy is obtained with N{$ _2- $N$ _{2\cdot 2} $} for the same mesh level. The result with N$ _2 -$N$ _{2\cdot 2} $ using intermediate mesh m$ _3 $ are of similar accuracy to that obtained with L$ _1 $ at a very fine mesh m$ _5 $.
\begin{table}[!t]
	\begin{center}
		\begin{tabular}{|c| c |c |c c c| c c c|}
			\hline
			\textbf{\centering{ Mesh}} &\textbf{Discretizations} & 
			\multicolumn{7}{|c|}{\textbf{DOFs}}   \\[1ex]
			\cline{3-9} 		
			&  & \textbf{Interface} & & \textbf{Bulk}& & &\textbf{Total} & \\  
			\hline
			& L$_1 $ &  20  & & 960 & & & 980 & \\			
			m$_1 $ & N$_2 $ &  22  & & 1056 & & & 1078 & \\
			& N$_2- $N$_{2\cdot 1} $ & 40 & & 1056 &  & & 1096 & \\
			& N$_2- $N$ _{2\cdot 2} $ & 58 & & 1056  &  & & 1114 & \\
			\hline
			& L$_1 $ &  38  & & 1824 & & & 1862 & \\			
			m$_2 $ & N$_2 $ &  40  & & 1920 &  & & 1960 & \\
			& N$_2- $N$_{2\cdot 1} $ &  76  & &1920&   &  & 1996 & \\ 		 
			& N$_2- $N$ _{2\cdot 2} $ & 112 & &1920&  &  & 2032 & \\
			\hline
			& L$_1 $ & 74 & & 3552 & & & 3626 & \\			
			m$_3 $ & N$_2 $ & 76 & &3648&  & & 3724 & \\
			& N$_2- $N$_{2\cdot 1} $ & 148   &  &3648&  &  & 3796 & \\ 		 
			& N$_2- $N$ _{2\cdot 2} $ & 220  & & 3648  & & &3868& \\
			\hline
			& L$_1 $ & 146 & & 7008 & & & 7154 & \\			
			m$_4 $ & N$_2 $ & 148 & &7104&  & & 7252 & \\
			& N$_2- $N$_{2\cdot 1} $ & 292 & &7104& & & 7396 & \\
			& N$_2- $N$ _{2\cdot 2} $ & 436 & &7104& & & 7540 &  \\
			\hline
			& L$_1 $ &  290  & & 13920 & & & 14210 & \\
			m$_5 $ & N$_2 $ &  292 & &14016&  &  & 14308 & \\
			& N$_2- $N$_{2\cdot 1} $ &  580  &  &14016&  &  & 14596 &  \\
			& N$_2- $N$ _{2\cdot 2} $ & 868 & &14016& &   & 14884 &   \\
			\hline				
			m$_6 $ & N$_2 $ & 580  & &27840&  & & 28420 &  \\
			\hline						
		\end{tabular} \caption{Degrees of freedom density data for both the discretization approaches with different mesh arrangements.} \label{table:DOFs}
	\end{center}
\end{table}

\begin{figure}[!t]
	\begin{center}	
		\subfloat{\includegraphics[width = 0.52\linewidth]{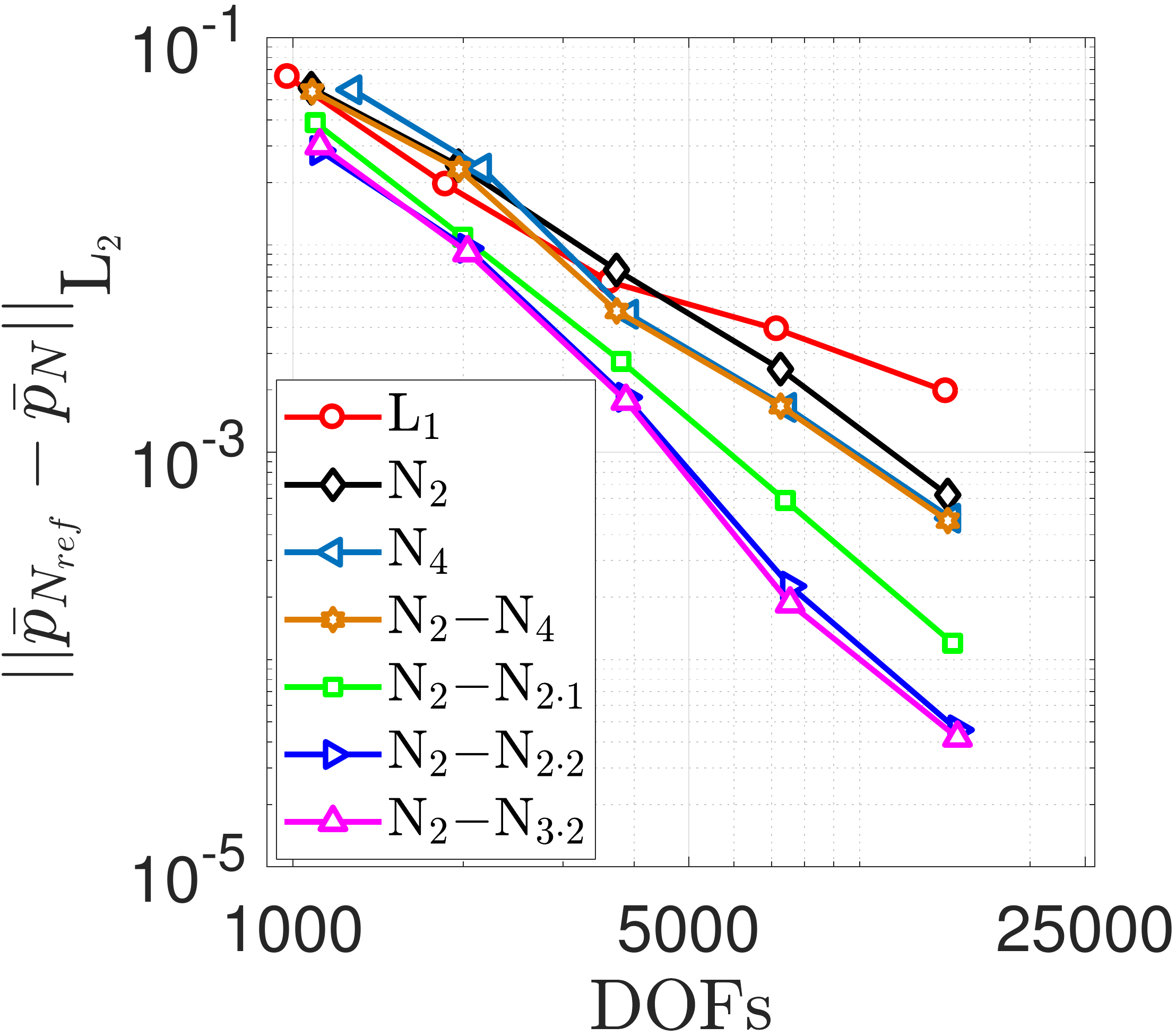}}
		\scriptsize{\caption {The $L_2 -$norm based error in the normalized contact pressure distribution for the  higher-continuous NURBS discretizations as compared to others using meshes m$ _1 $ to m$ _5 $.}\label{fig:mesh_conv_order}}
	\end{center}
\end{figure}

Besides, a comparison with the higher than $ C^1- $continuous based NURBS discretizations, i.e. with N$ _2- $N$ _4 $, its counterpart standard N$ _4 $, and N$ _2- $N$ _{3\cdot 2} $, is also carried out and the obtained results are shown in  Fig.~\ref{fig:mesh_conv_order}. It can be seen that although on increasing the continuity from $ C^0 $ (with L$_1  $) to $ C^1 $ (with N$ _2 $), a more accurate result can be obtained at a fixed mesh, i.e. at intermediate m$ _4 $ and fine  m$ _5 $ mesh level. But, with the more than $ C^1- $continuous NURBS discretization, i.e. with N$ _2- $N$ _4 $ or N$ _4 $, the absolute error values are nearly equivalent to those obtained with N$ _2 $ for the same mesh. This observation also holds for the higher-continuous version, i.e. 
N$ _2- $N$ _{3\cdot 2} $, of the N$ _2- $N$ _{2\cdot2}$ discretization. Consequently, it shows that NURBS functions smoother than $ C^1 $ are not advantageous to analyze the Hertzian contact problem. This behaviour is due to the inability of higher smooth NURBS functions to adequately capture the local changes in the distribution of the contact pressure result.

\begin{figure}[!t]
	\centering
	{\includegraphics[width = 0.48\linewidth]{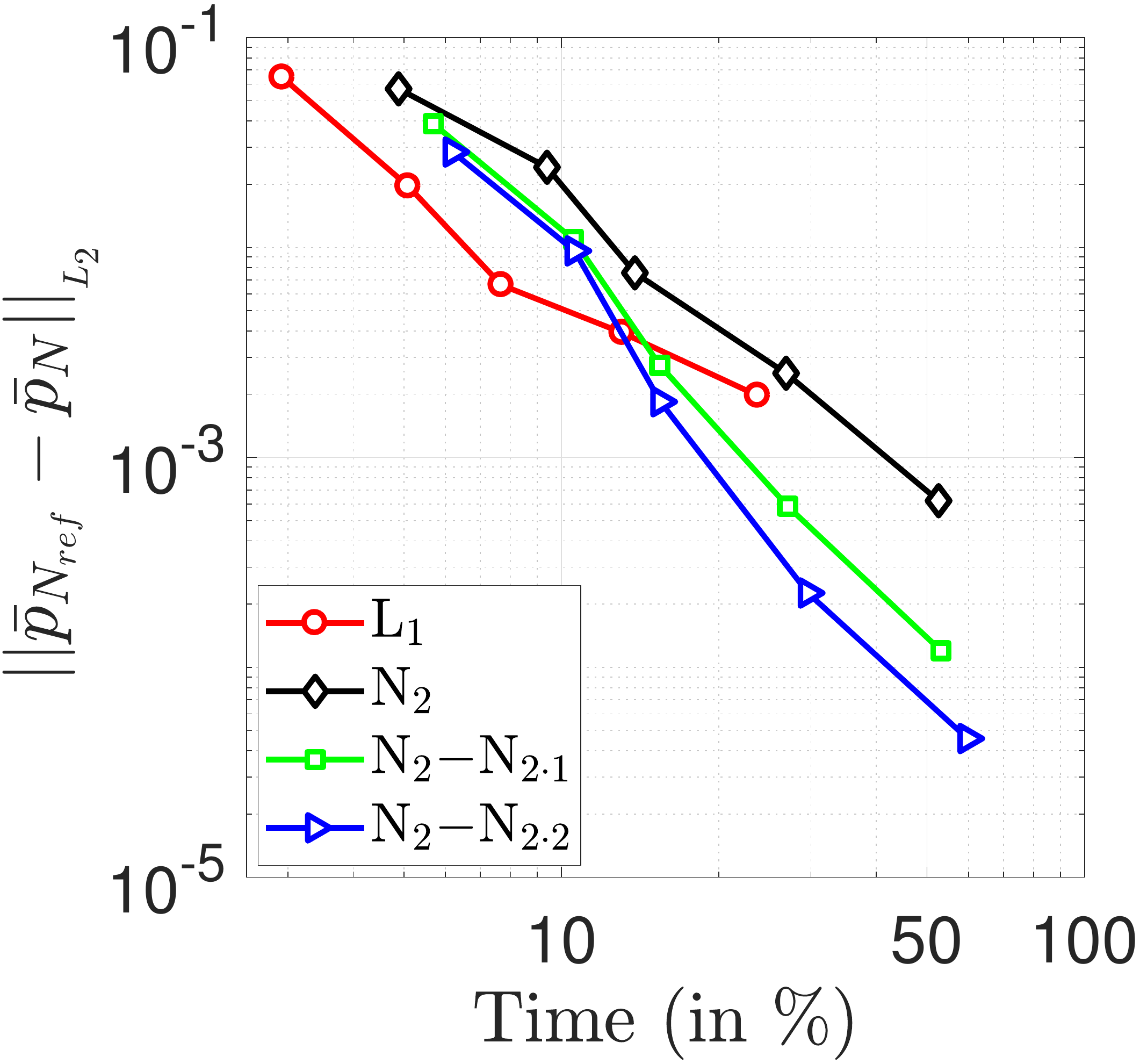}}
	\scriptsize{\caption {Error with total analysis time for different discretizations with meshes m$ _1 $ to m$ _5 $.} \label{fig:Time_Hertz}}		
\end{figure}

Finally, we compare the associated computational efforts in terms of the overall analysis time with both the  discretization approaches for similar accuracy.  Figure~\ref{fig:Time_Hertz} shows the plots for $ L_2 $-norm of the error in the contact pressure with the total analysis time with meshes m$ _1 $ to m$ _5 $ for L$ _1 $, N$ _2 $, N$ _2- $N$ _{2\cdot 1} $ and N$ _2- $N$ _{2\cdot 2} $ discretizations. Here, the time in $ \% $ is computed using Eq.~(\ref{eq:analysis_time}) and the maximum analysis time, that is used for normalizing the time values, is with N$ _2 $ at mesh m$ _6 $. From Fig.~\ref{fig:Time_Hertz} it is evident that to attain the accuracy similar to N$ _2 $, N$ _2- $N$ _{2\cdot 1} $ and N$ _2 -$N$ _{2\cdot 2} $ require much lower computational efforts. With a closer look, it turns out that N$ _2- $N$ _{2\cdot 2} $ at mesh m$ _3 $ and m$ _4 $ take approximately $ 42.48 \% $ and $ 43.74\% $ lower analysis time than that with N$ _2 $ at m$ _4 $ and m$_5 $,  respectively, to attain the similar accuracy. Subsequent comparison with L$ _1 $ discretization shows that to attain the comparable accuracy, N$ _2- $N$ _{2\cdot 2} $ at mesh m$ _1 $ and m$ _2 $ takes a slightly more analysis time as compared to L$ _1 $ at m$ _2 $ and m$ _3 $, respectively. Further, N$ _2- $N$ _{2\cdot 2} $ at mesh m$ _3 $ requires approximately $ 34.86\% $ lower computational efforts as compared to L$ _1 $ at m$ _5 $ for the similar accuracy level. The additional usage of a fine mesh, e.g. m$ _4 $ or m$ _5 $, with N$ _2- $N$ _{2\cdot 2} $ leads to a significant gain in the accuracy, which, however, comes at the expense of additional computational efforts as compared to L$ _1 $ at m$ _5 $.

In summary, this example shows the much higher accuracy and the computational efficiency of the VO based NURBS discretizations as compared to standard N$ _2 $ based discretization. For the fixed mesh, VO NURBS yields much more accurate results at a cost slightly more than with the standard NURBS based discretization. To deliver a similar accuracy, it takes much lower computational efforts as compared to N$ _2 $ discretization.

\subsection{Two elastic rings contact problem}
\begin{figure}[!t]
	\centering	
	{\includegraphics[width=.385\linewidth]{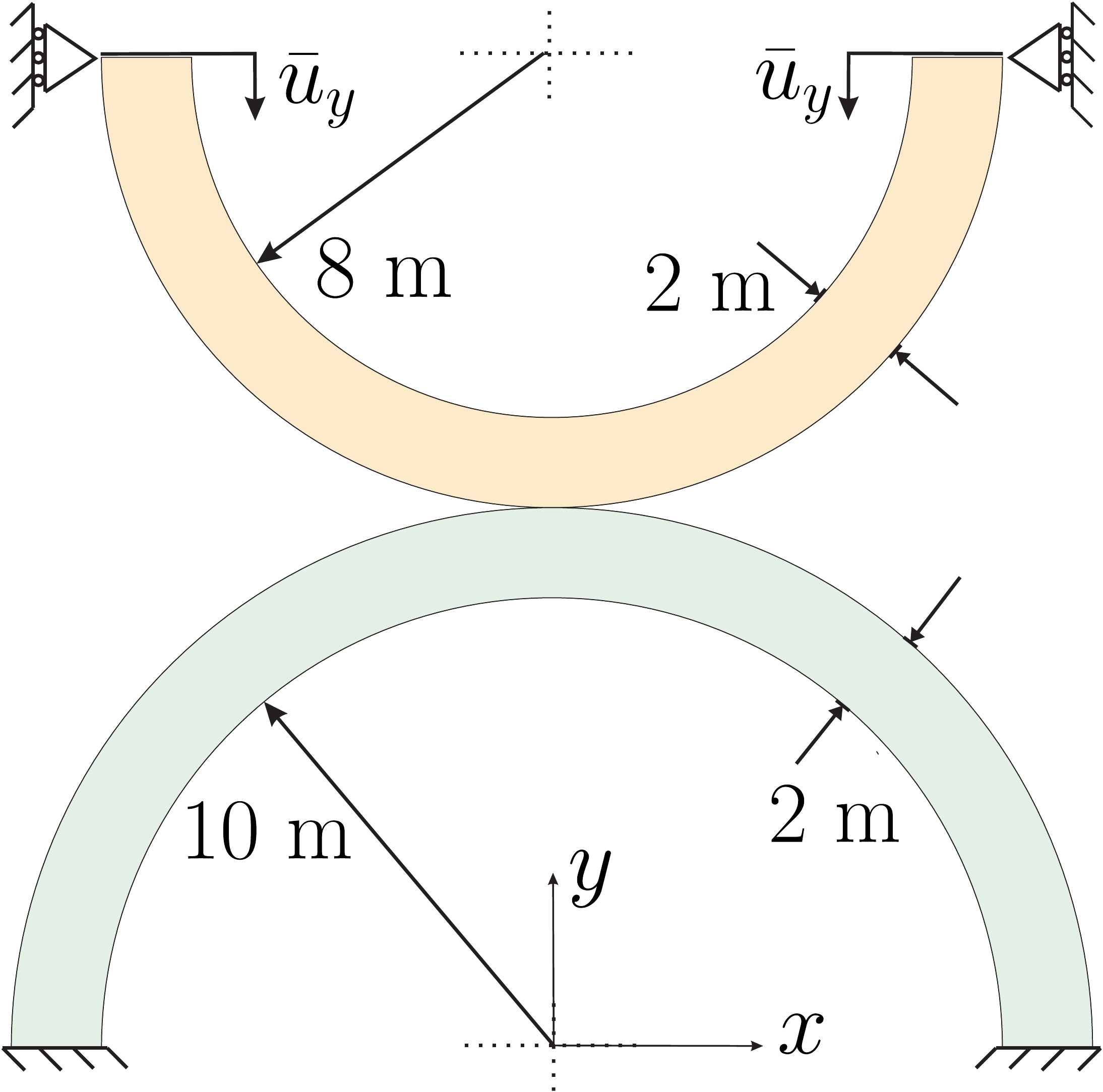}
		\caption{The setup of the two rings contact problem.} \label{fig:Setup_2Ring}}
\end{figure}
In the final example, the frictional contact between the two deformable rings undergoing the large deformation is simulated. The setup of the problem, taken from~\cite{DIAS2019}, is shown in Fig.~\ref{fig:Setup_2Ring}. The outer and inner radii of the upper and lower rings are $ (10, 8) $ and $ (12, 10) $ m, respectively. The bottom surface of the lower ring is held stationary by fixing the displacement in both directions. A Neo-Hookean hyper-elastic material model, according to Eq.~(\ref{key:Neo-Hookean}), is considered. The material parameters are: Young's modulus $ E_{\textrm{upper}} = 1\times 10^2 $ N$ / $m$ ^2 $ and $ E_{\textrm{lower}} = 3 \times 10^2 $ N$ / $m$ ^2 $, and Poisson's ratio $ \nu = 0.3 $ for both the rings. The plane strain assumption is used. Due to the symmetry with respect to the $ y- $axis, only half of the setup is utilized. The bottom surface of the upper ring is used as a slave surface, while the outer top surface of the lower one as the master surface. A vertical displacement $ \bar{u}_y = 4.0 $ m is applied in the downward direction to the top surface of the upper ring in $ 40 $ uniform load steps. The penalty parameters are taken as $ \epsilon_{\text{N}} = 100 $  and $ \epsilon_{\text{T}} = 10 $. The coefficient of friction $ \mu_f = 0.1 $ is used. The objective is to predict the normal and tangential reaction responses at the contact interface of the two rings using the VO based NURBS discretizations. The results obtained with N$ _2 $ based discretization are used for comparisons.
\begin{figure}[!ht]
	\centering
	\subfloat[]{\includegraphics[width=.25\linewidth]{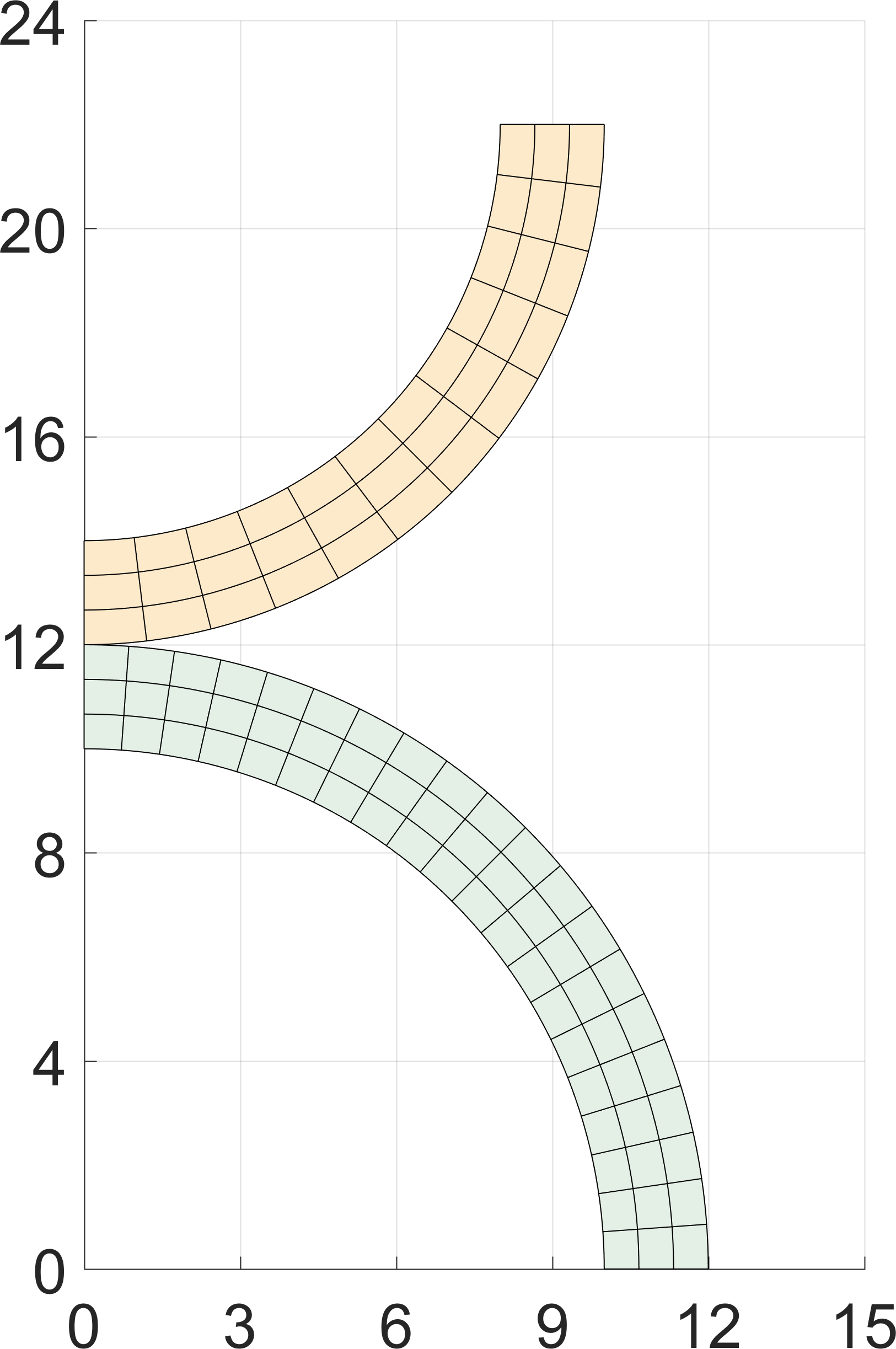}\label{fig:Def_t=0}} \\
	\subfloat[]{\includegraphics[width=.246\linewidth]{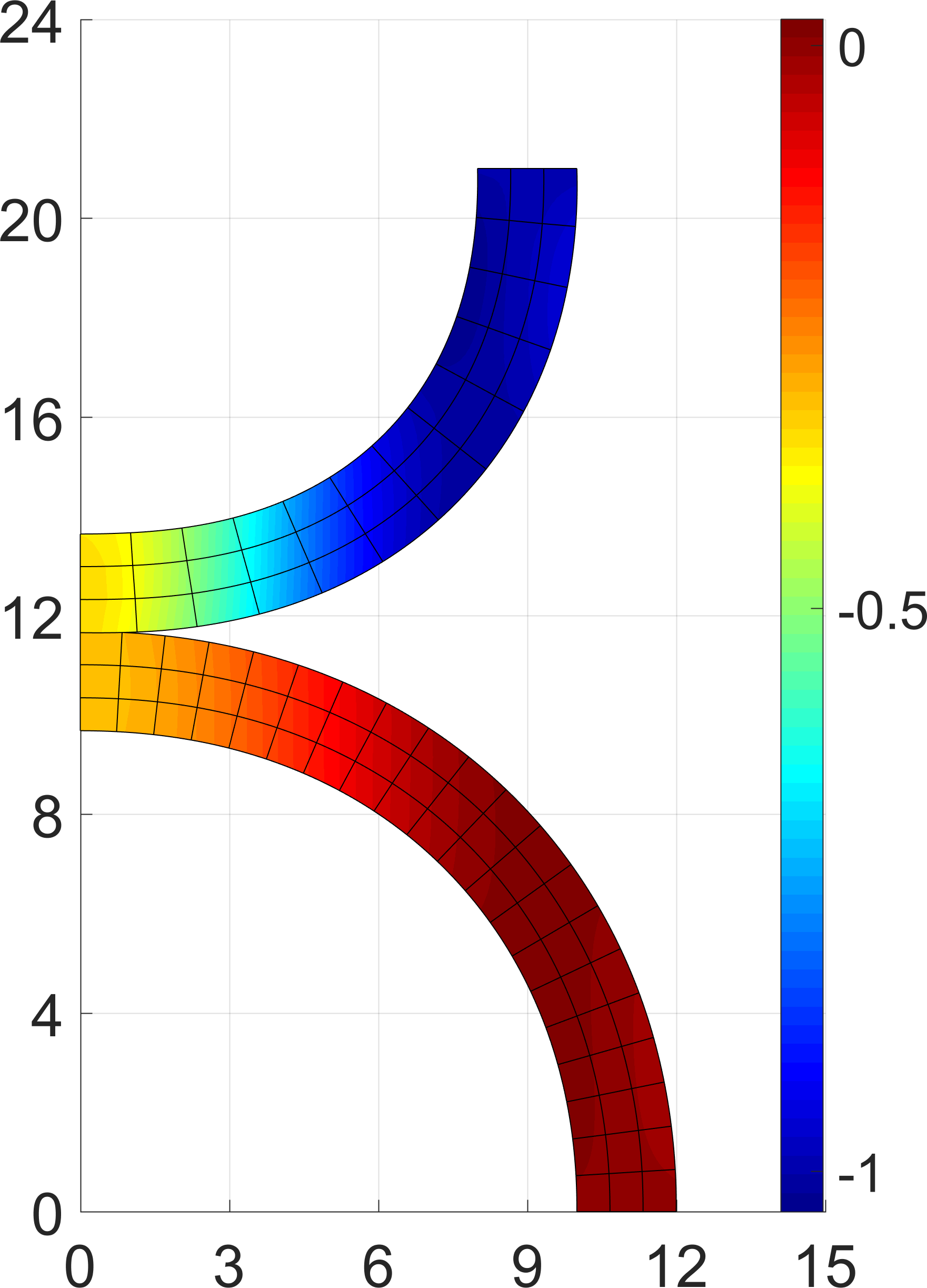}\label{fig:Def_t=10}}~ 
	\subfloat[]{\includegraphics[width=.246\linewidth]{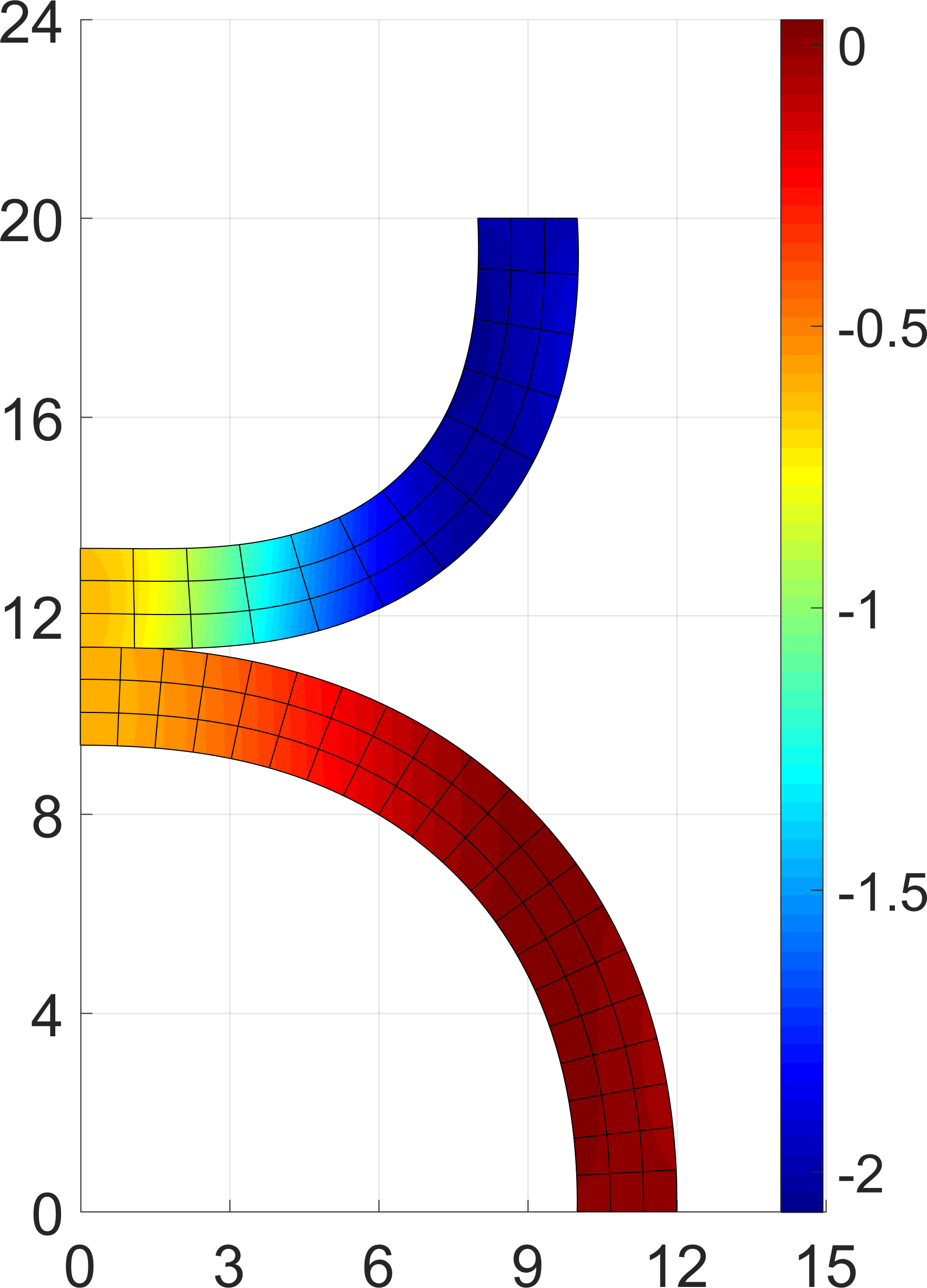}\label{fig:Def_t=20}}~
	\subfloat[]{\includegraphics[width=.246\linewidth]{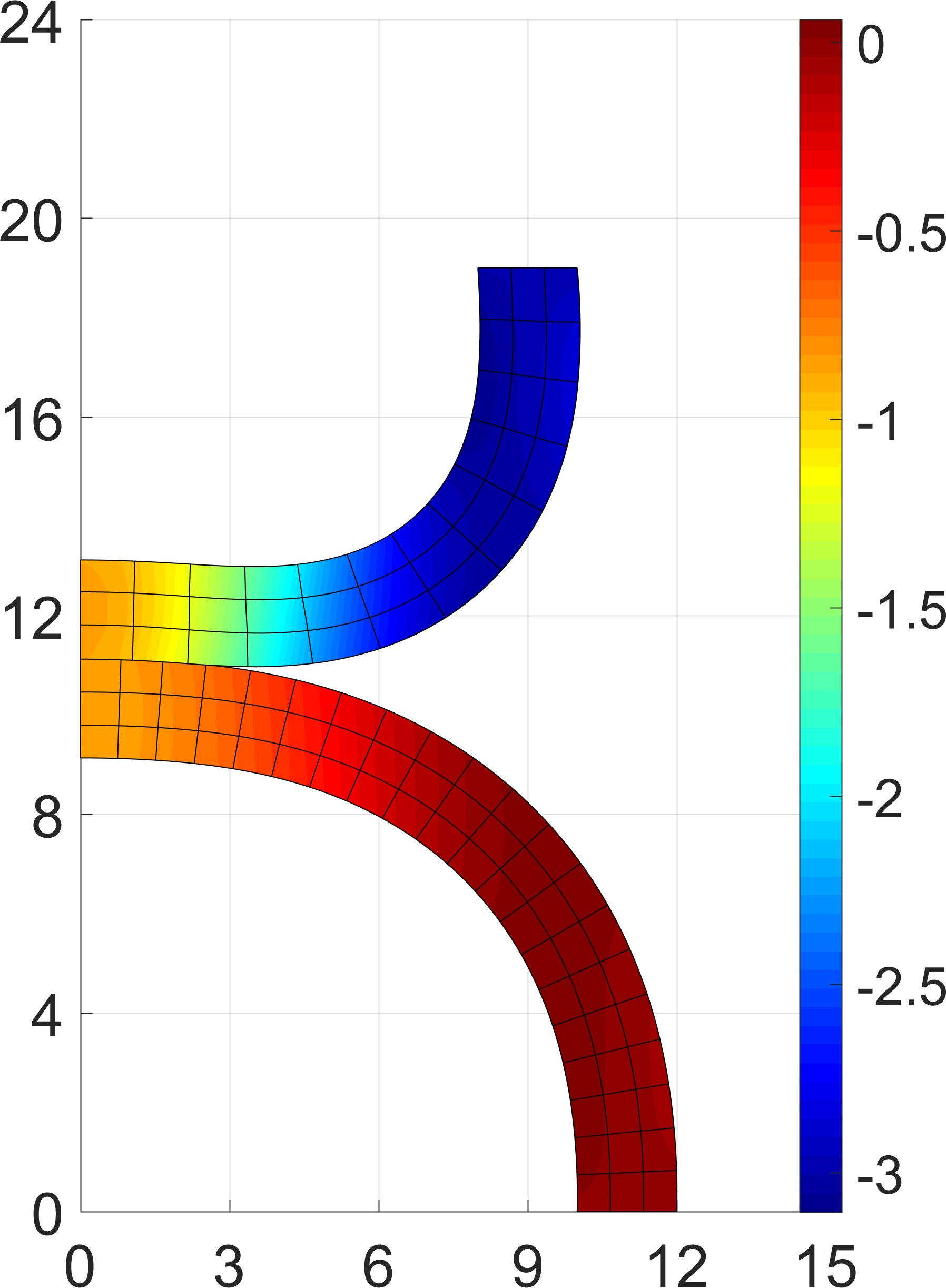}\label{fig:Def_t=30}}~
	\subfloat[]{\includegraphics[width=.246\linewidth]{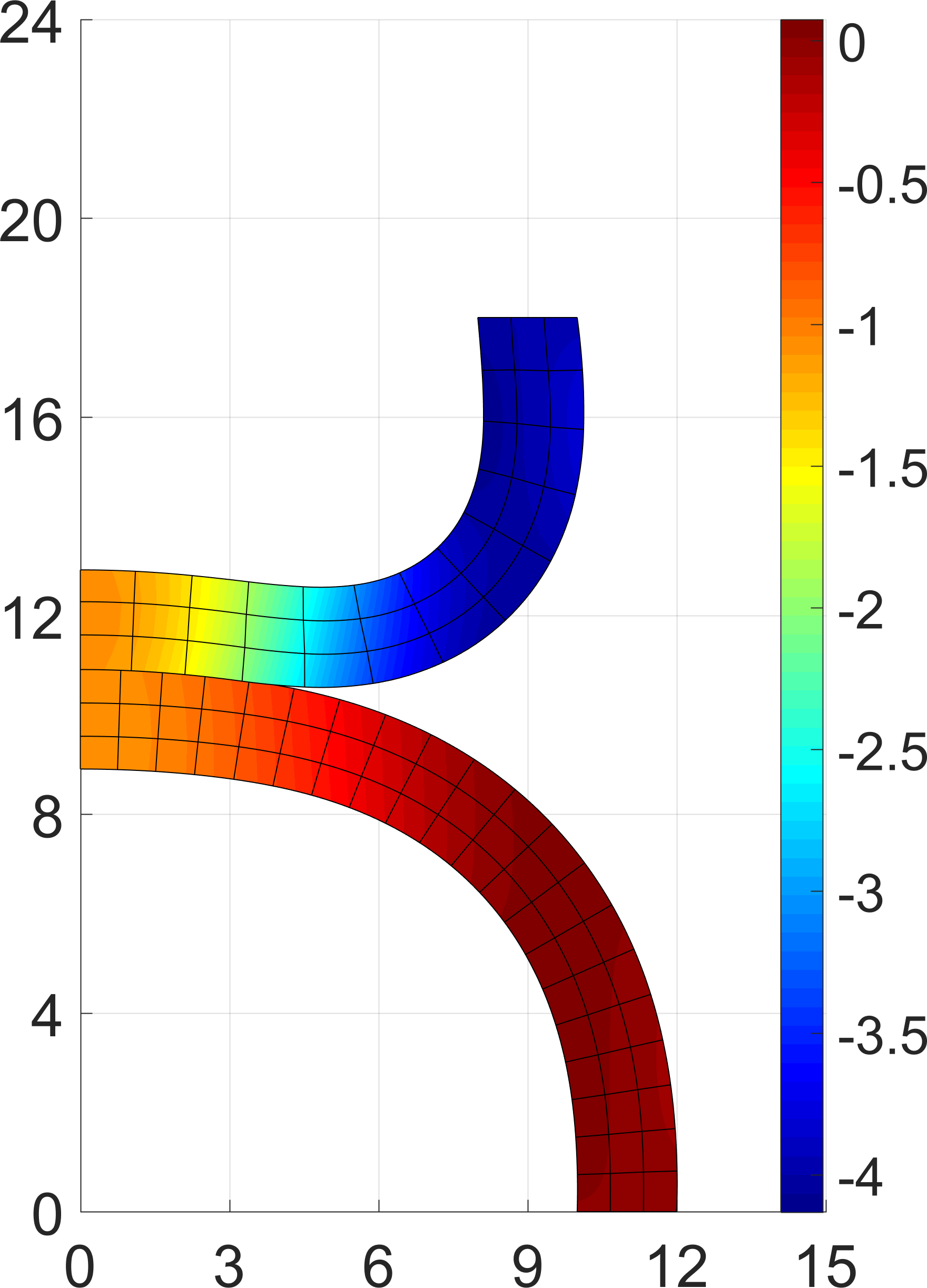}\label{fig:Def_t=40}}
	\caption{(a) The coarse mesh m$ _1 $ and the initial configuration of the setup. The deformation configuration with the vertical displacement contours at load steps (b) $ t =10 $; (c) $ t = 20 $; (d) $ t = 30 $; and  (e) $ t = 40 $.} \label{fig:Def_2Rings}
\end{figure}

A coarse mesh having $ 12 \times 3 $ number of elements in the angular and radial  directions of the upper geometry, and $ 20 \times 3 $ of the lower one is shown in Fig.~\ref{fig:Def_t=0}. Overall, three different meshes are used to carry out this simulation. These meshes have $ (12, 20), (24, 40) $, and $ (48, 80) $ number of elements across the contact boundary layer of the upper and lower geometries, and are denoted by m$ _1 $, m$ _2 $, and m$ _3 $, respectively. Numerical investigation reveals that N$ _2 $ based discretization using fine mesh m$ _3 $ delivers the converged normal and tangential contact pressure profiles for this example. Thus, the obtained result are used as the reference solutions. The deformed shapes of the setup with the vertical displacement field $ u_y $ at four different values of load steps with N$ _2 $ using mesh m$ _1 $ are shown in Figs.~\ref{fig:Def_t=10}-\ref{fig:Def_t=40}.

\begin{figure}[!t]
	\centering	
	\subfloat{\includegraphics[width=.4\linewidth]{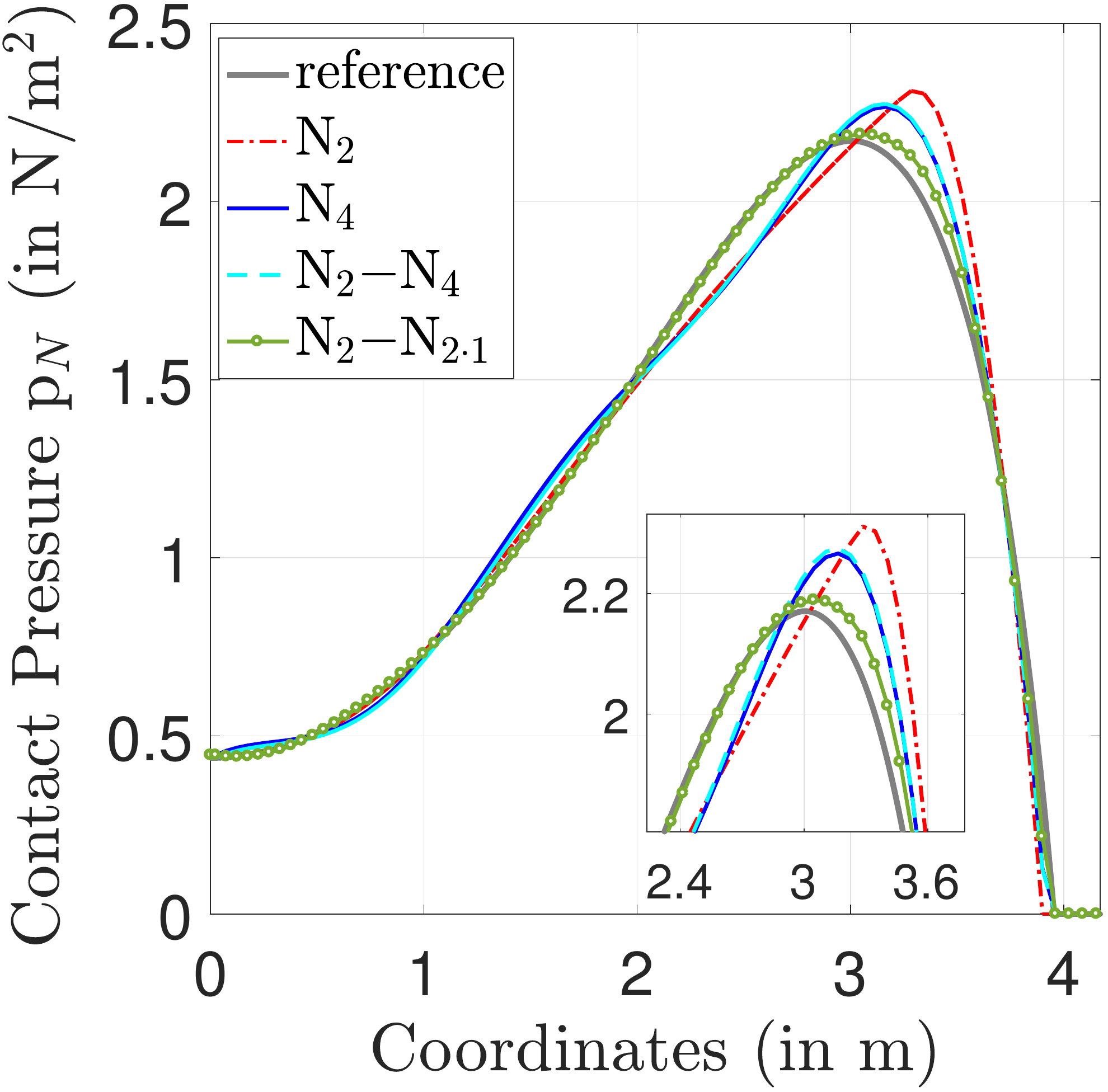}\label{fig:pN_m12}} ~~~~~~~
	\subfloat{\includegraphics[width=.4\linewidth]{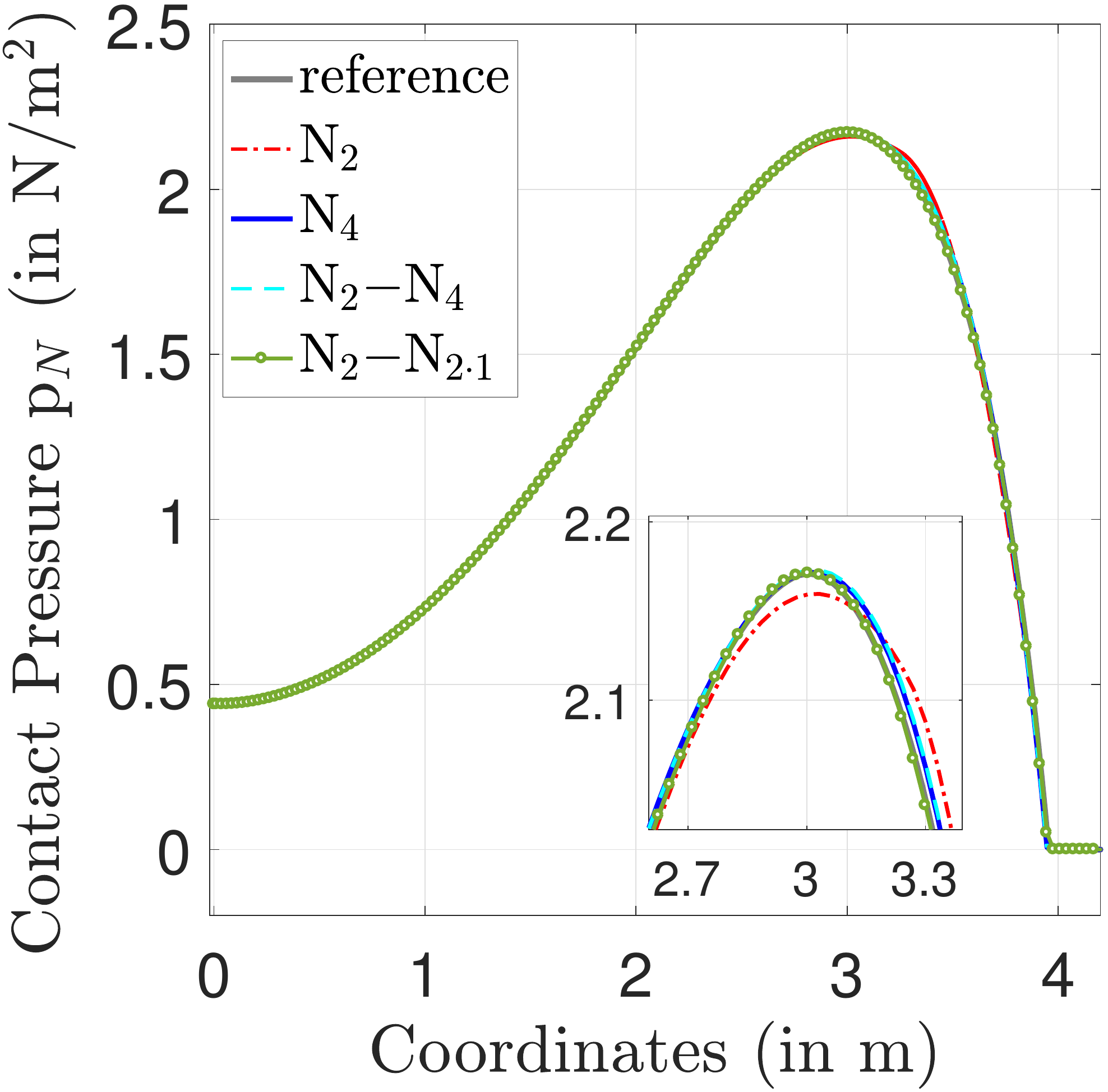}\label{fig:pN_m24}}\\ 
	\subfloat{\includegraphics[width=.4\linewidth]{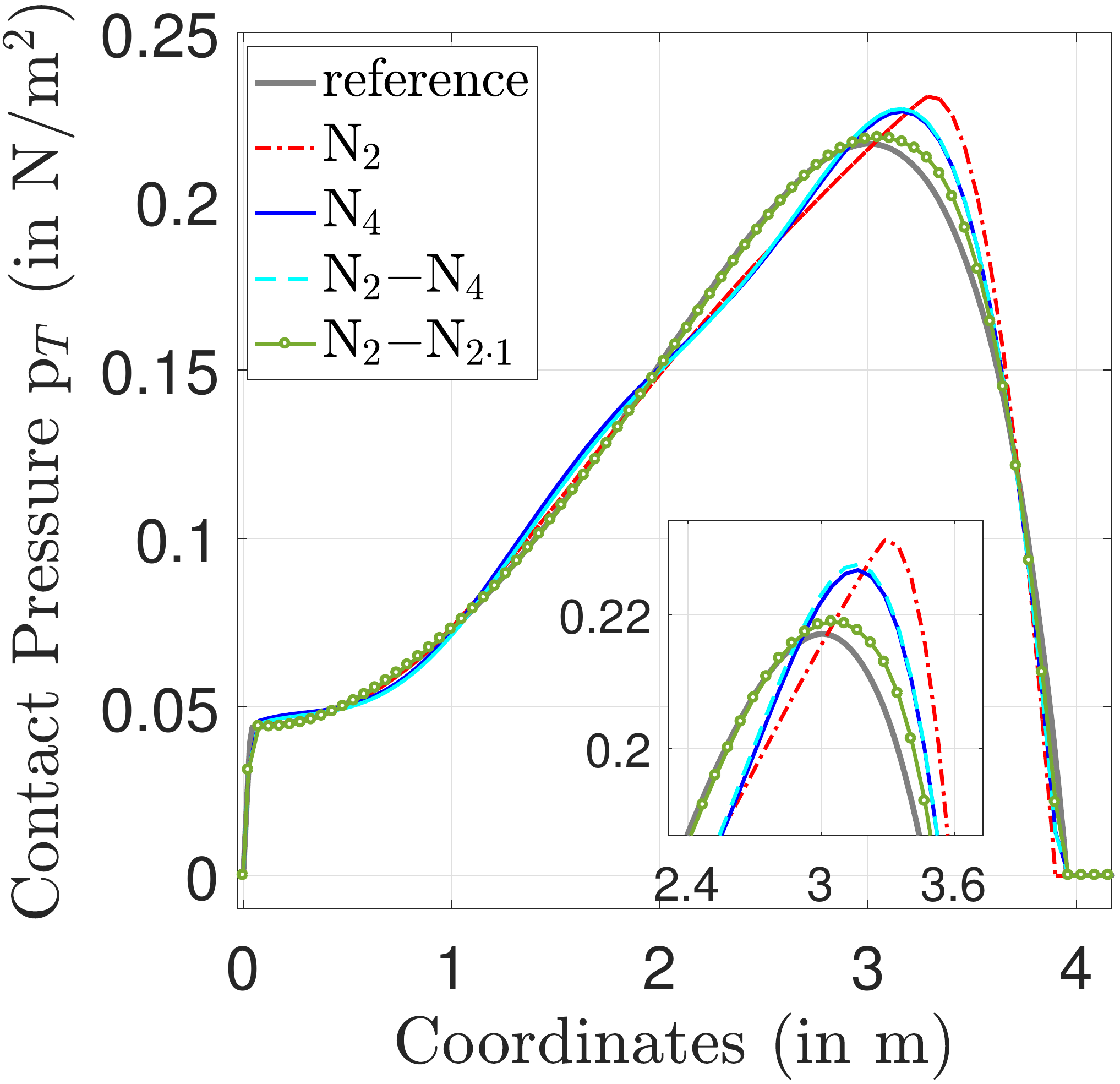}\label{fig:pT_m12}} ~~~~~~~
	\subfloat{\includegraphics[width=.4\linewidth]{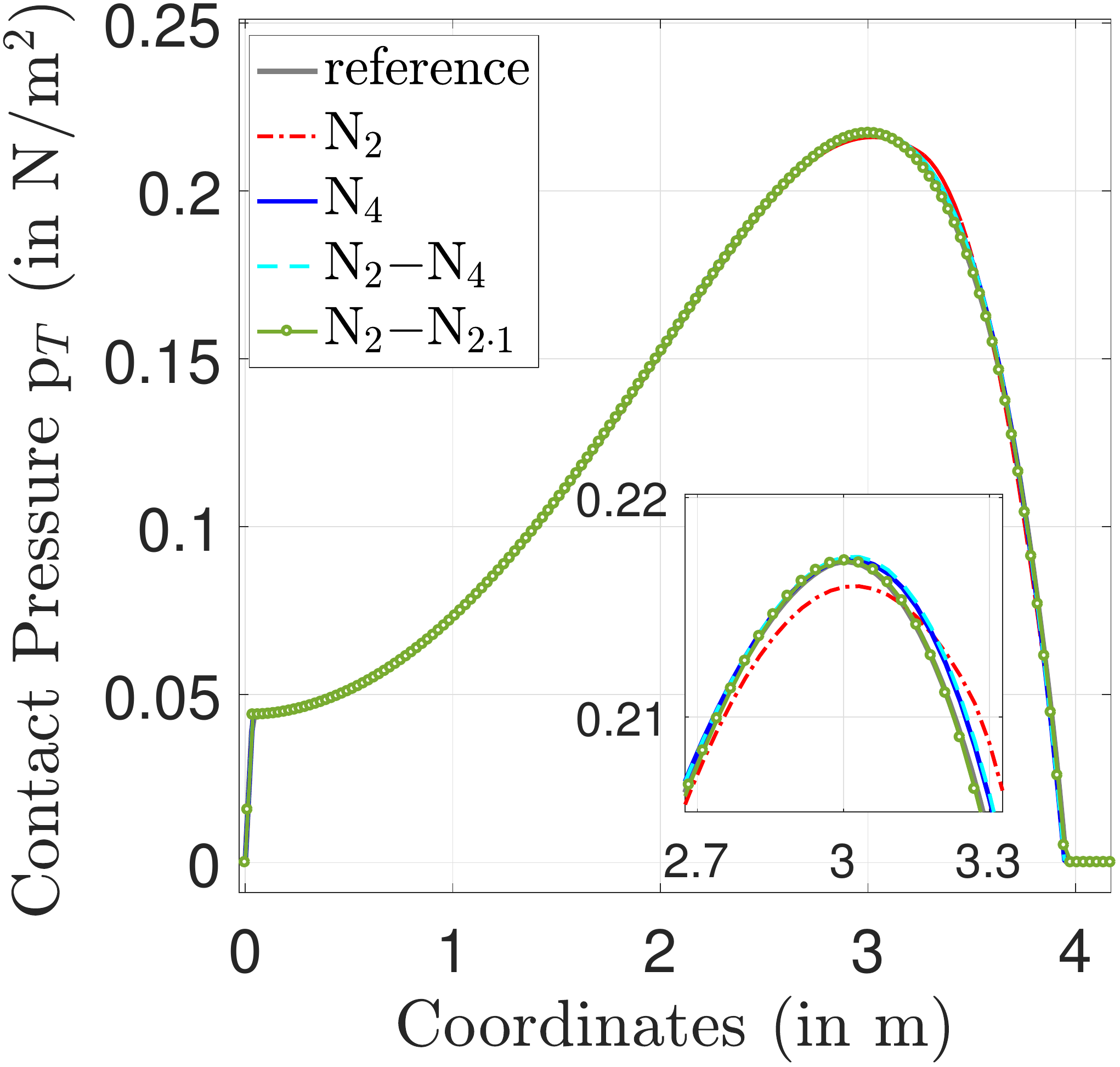}\label{fig:pT_m24}}
	\caption{Comparison of the contact pressure p\mbox{$ _N $} (top row) and p\mbox{$ _T $} (bottom row) for different discretization at meshes m\mbox{$ _1 $} (first column) and m\mbox{$ _2 $} (last column).} \label{fig:2Ring_plots}
\end{figure}

The normal contact pressure p$ _N $ and tangential contact pressure p$ _T $ with standard N$ _2 $ and N$ _4 $, and VO based N$ _2- $N$ _4 $ and N$ _2- $N$ _{2\cdot 1} $ discretizations using meshes m$ _1 $ and m$ _2 $ are shown in Fig.~\ref{fig:2Ring_plots}. The corresponding $ L_2 $-norm of the errors in p$ _N $ and p$ _T $, calculated using Eq.~(\ref{eq:L2_norm}), are summarized in Table~\ref{table:2Ring_conv}. The overall computational time taken by different discretizations using both the meshes are provided in Table~\ref{table:2Ring_Time}. From Fig.~\ref{fig:2Ring_plots} two observations are made. First, the contact pressure curves of p$ _N $ and p$ _T $ with higher-continuous VO based N$ _2- $N$ _4 $ are nearly identical to those obtained with standard N$ _4 $ discretization for meshes m$ _1 $ and m$ _2 $. However, only a marginal improvement in the accuracy of results, attributed to higher smooth NURBS functions, is obtained as compared to N$ _2 $ discretization. Second, concerning the distribution of contact pressure across the contact interface at a fixed mesh, N$ _2- $N$ _{2\cdot 1} $ is able to capture it more accurately than N$ _2 $ and N$ _2- $N$ _4 $ (or N$ _4 $) based discretizations. As can be seen from Tables~\ref{table:2Ring_conv} and~\ref{table:2Ring_Time}, 
it improves the accuracy of the normal contact pressure distribution by a factor of $ 3.4 $ and $ 9.9 $, and the tangential contact pressure distribution by $ 3.4 $ and $ 7.7 $ at a slightly more cost, i.e. $ 5.8 \%$ and $ 7.3 \% $, than with the N$ _2 $ based discretization with meshes m$ _1 $ and m$ _2 $. In summary, this example demonstrates the benefit of using VO based NURBS discretization as compared to standard NURBS discretization. It provides a major gain in the accuracy at a cost only slightly more than with the standard NURBS discretization.	

\begin{table}[!t]
	\begin{center}
		\begin{tabular}{c c c c c c}
			\hline
			\textbf{Discretization} &  \multicolumn{2}{c}{{$ || \textrm{p}_{N_{ref}} - \textrm{p}_N ||_2 $}} &  & \multicolumn{2}{c}{{$ || \textrm{p}_{T_{ref}} - \textrm{p}_T ||_2 $}}\\ [0.4ex]
			\cline{2-3}\cline{5-6}
			&  m$ _1 $ & m$ _2 $ & & m$ _1 $ & m$ _2 $ \\			
			\hline
			N$ _2 $ & $ 24.38 \times 10^{-2}$ & $ 34.32 \times 10^{-3}$ &     & $ 24.38 \times 10^{-3}$ & $ 34.33 \times 10^{-4} $ \\
			N$ _4 $ & $ 17.45 \times 10^{-2} $& $ 25.78 \times 10^{-3} $&     & $ 17.45 \times 10^{-3}$ & $ 25.99 \times 10^{-4} $ \\
			N$ _2- $N$ _4 $ & $ {17.18} \times 10^{-2} $ & $ 26.80 \times 10^{-3} $ &  & $ 17.18 \times 10^{-3}$ & $ 27.01 \times 10^{-4} $ \\
			N$ _2- $N$ _{2\cdot1} $ & $ 7.095 \times 10^{-2} $ & $ 3.433 \times 10^{-3} $ &    & $ 7.097 \times 10^{-3} $  &  $ 4.429 \times 10^{-4} $ \\
			\hline
		\end{tabular} \caption{$ L_2 $-norm of the error in p$ _N $ and p$ _T $ for different discretizations and meshes.} \label{table:2Ring_conv}
	\end{center}
\end{table}
\begin{table}[!t]
	\begin{center}
		\begin{tabular}{c c c c}
			\hline
			\textbf{Discretization} &  & \multicolumn{2}{c}{{ Time }}\\ [0.4ex]
			\cline{3-4}
			&  & m$ _1 $     & m$ _2 $ \\		
			\hline
			N$ _2 $              &  &  $ 1.000 $  &  $ 1.357 $ \\
			N$ _4 $              &  &  $ 1.451 $  &  $ 2.296 $ \\
			N$_2- $N$ _4 $       &  &  $ 1.216 $  & $ 1.724 $ \\
			N$ _2-$N$_{2\cdot1}$ &  &  $ 1.058 $  &  $ 1.475 $ \\
			\hline
		\end{tabular} \caption{Overall analysis time taken by various discretizations with different meshes. The time with N$ _2 $ using mesh m$ _1 $ is used for normalization.} \label{table:2Ring_Time}
	\end{center}
\end{table}

\section{Conclusion} \label{sec:conclusion}
In this work, a novel varying-order based NURBS discretization methodology is introduced for isogeometric contact analysis. With the proposed method, a much higher numerical accuracy is achieved in the contact computations even at a coarse mesh as compared to existing fixed-order based NURBS discretization. The improvement in the accuracy is obtained for both the small and large deformation contact problems, considering with or without friction. In other words, the proposed method exhibits a major gain in computational efficiency over the standard NURBS based discretizations for the similar accuracy level. The gain in the accuracy is attributed to the application of additional order-elevation based refinement strategy to the contact boundary layer of the NURBS discretized geometry. With this, a large number of additional DOFs are introduced across the contact interface than with the standard NURBS discretizations at a fixed mesh. The accuracy further improves on further elevating the interpolation-order of the contact boundary layer. As shown in the ironing example, the quality of the result improves marginally on additionally increasing the inter-element continuity of order-elevation  discretized contact layer. Refining the geometry with the proposed method is quite robust, as illustrated by the numerical examples. The simplicity of the current approach allows itself to be conveniently embedded into the existing isogeometric contact codes. The extension of the developed methodology to three-dimensional contact problems in the framework of isogeometric analysis is currently underway.

\section*{Acknowledgments}
The authors gratefully acknowledge the support from SERB, DST, under projects  SB/FTP/\newline ETA-0008/2014 and IMP/2019/000276. The authors would like to express their gratitude towards Prof. Roger A. Sauer, AICES RWTH Aachen University, and the anonymous reviewers for their helpful comments and suggestions. The first author also gratefully acknowledges the support from IIT Guwahati and Ministry of Human Resource and Development (MHRD), Government of India, for the financial assistantship.

\titleformat{\section}{\normalfont\Large\bfseries}{\appendixname~\thesection.}{1em}{}
\begin{appendices} 
	\numberwithin{equation}{section}
	\makeatletter 		
	\newcommand{\section@cntformat}{Appendix \thesection:\ }
	\makeatother
	\section{Contact tangent matrices}  \label{appendix:A}
	In this appendix, we provide the contact tangent matrices for the full-pass contact formulation used for the isogeometric analysis of the above numerical examples. In the following, the subscript $ n+1 $, denoting the contact quantities in the current configuration, is dropped for convenience. The presented expressions for the elemental contact tangent matrices $ \mathbf{k}_{\textrm{c}}^{kke} $ (where $ k= \{s,m\} $) are obtained from the linearization of the contact force vector $ \mathbf{f}_{\textrm{c}}^{ke} $. We refer to the works by Sauer and De Lorenzis~\cite{unbiased2013, unbiased2015} for the detailed description of the formulation, including the linearization procedure.
	
	\begingroup		
	\begin{align}
	\mathbf{k}_{\textrm{c}}^{\textrm{ss}e} &= -\int_{\Gamma^{\textrm{s}e}_c} \mathbf{N}^{\textrm{s}\textrm{T}} \, \frac{\partial \bm{t}^s}{\partial \mathbf{u}^{\textrm{s}e}} ~\textrm{d}\Gamma - \int_{\Gamma^{\textrm{s}}_c} \mathbf{N}^{\textrm{s}\textrm{T}} \, \bm{t}^{\textrm{s}} \otimes \bar{\bm{\uptau}}^1 \, \mathbf{N}^{\textrm{s}\textrm{T}}_{,1}  ~\textrm{d}\Gamma  \nonumber \\
	\mathbf{k}_{\textrm{c}}^{\textrm{sm}e} &=  \int_{\Gamma^{\textrm{s}e}_c} \mathbf{N}^{\textrm{s}\textrm{T}} \frac{\partial \bm{t}^s}{\partial \mathbf{u}^{\textrm{m}e}}  ~\textrm{d}\Gamma\\
	\mathbf{k}_{\textrm{c}}^{\textrm{ms}e} &=   \int_{\Gamma^{\textrm{s}e}_c} \mathbf{N}^{\textrm{m}\textrm{T}}_{,1}\, \bm{t}^{\textrm{s}}\otimes\frac{\partial \bar{\xi}^{\textrm{m}}}{\partial \mathbf{u}^{\textrm{s}e}}  ~\textrm{d}\Gamma + \int_{\Gamma^{\textrm{s}e}_c} \mathbf{N}^{\textrm{m}\textrm{T}} \, \frac{\partial \bm{t}^s}{\partial \mathbf{u}^{\textrm{s}e}} ~\textrm{d}\Gamma + \int_{\Gamma^{\textrm{s}}_c} \mathbf{N}^{\textrm{m}\textrm{T}} \, \bm{t}^{\textrm{s}} \otimes \bar{\bm{\uptau}}^1 \, \mathbf{N}^{\textrm{s}\textrm{T}}_{,1} ~\textrm{d}\Gamma \nonumber \\
	\mathbf{k}_{\textrm{c}}^{\textrm{mm}e} &=   \int_{\Gamma^{\textrm{s}}_c} \mathbf{N}^{\textrm{m}\textrm{T}}_{,1}\, \bm{t}^{\textrm{s}}\otimes\frac{\partial \bar{\xi}^{\textrm{m}}}{\partial \mathbf{u}^{\textrm{m}e}} 
	~\textrm{d}\Gamma  +  \int_{\Gamma^{\textrm{s}}_c}  \mathbf{N}^{\textrm{m}\textrm{T}} \, \frac{\partial \bm{t}^s}{\partial \mathbf{u}^{\textrm{m}e}} ~\textrm{d}\Gamma \nonumber
	\end{align}
	\endgroup
	
	where $ \mathbf{N}^{\textrm{s}} = \mathbf{N}(\xi^{\textrm{s}})  $ and $ \mathbf{N}^{\textrm{m}} = \mathbf{N}(\bar{\xi}^{\textrm{m}}) $. The  normal and tangential components of the contact traction vector $ \bm{t}^\textrm{s} := \bm{t}^\textrm{s}_{\textrm{N}} - \bm{t}^\textrm{s}_{\textrm{T}} $ are defined in Eqs.~(\ref{eq:normal_contact}) and~(\ref{eq:tangential_traction}) for frictionless and frictional contact, respectively. The derivatives of $ \bm{t}^{\textrm{s}} $ are given as 
	\begin{equation}\label{key}
	\frac{\partial \bm{t}^{\textrm{s}}}{\partial \mathbf{u}^{ke}} = \frac{\partial \bm{t}^{\textrm{s}}_{\textrm{N}}}{\partial \mathbf{u}^{ke}} - \frac{\partial \bm{t}^{\textrm{s}}_{\textrm{T}}}{\partial \mathbf{u}^{ke}}, ~~~~~~~~\textrm{where } k = \textrm{s, m}\,.
	\end{equation}
	The derivative of normal component of contact traction vector for each contact surface is given by
	\begin{equation}\label{key}
	\begin{aligned}
	\frac{\partial \bm{t}^{\textrm{s}}_{\textrm{N}}}{\partial \mathbf{u}^{se}} &= - \epsilon_{\text{N}}\left[  \bm{I} - c^{11} \bar{\bm{\uptau}}_1 \otimes \bar{\bm{\uptau}}_1 \right] \mathbf{N}^{\textrm{s}}\,,~~~\textrm{and} \\
	\frac{\partial \bm{t}^{\textrm{s}}_{\textrm{N}}}{\partial \mathbf{u}^{me}} &= \epsilon_{\text{N}}\left[  \bm{I} - c^{11} \bar{\bm{\uptau}}_1 \otimes \bar{\bm{\uptau}}_1 \right] \mathbf{N}^{\textrm{m}} + \epsilon_{\text{N}}\textrm{g}_{\textrm{N}}c^{11} \bar{\bm{\uptau}}_1 \otimes \bar{\bm{n}}\, \mathbf{N}^{\textrm{m}}_{,1}
	\end{aligned}
	\end{equation}
	with $ c^{11} = 1/[m_{11} - \textrm{g}_{\textrm{N}}	\bar{\bm{n}}\cdot \bar{\bm{x}}_{,11} ]$. The derivative of tangential component of contact traction vector for a stick step are given by
	\begingroup
	\allowdisplaybreaks
	\begin{align}
	\frac{\partial \bm{t}^{\textrm{s}}_{\textrm{T}}}{\partial \mathbf{u}^{se}} &= \epsilon_{\text{T}} \, c^{11} \bar{\bm{\uptau}}_1 \otimes \bar{\bm{\uptau}}_1\, \mathbf{N}^{\textrm{s}}\,,~~~\textrm{and} \\
	\frac{\partial \bm{t}^{\textrm{s}}_{\textrm{T}}}{\partial \mathbf{u}^{me}} &= -\epsilon_{\text{T}} \, c^{11} \bar{\bm{\uptau}}_1 \otimes \bar{\bm{\uptau}}_1\, \mathbf{N}^{m} + \epsilon_{\text{T}}\left[ \mathbf{N}^{m} - \mathbf{N}^{m}(\xi_{sl\, n}^{\textrm{m}}) \right] + \epsilon_{\text{T}} \, c^{11} \bar{\bm{\uptau}}_1 \otimes \bar{\bm{n}}\,\textrm{g}_{\textrm{N}} \mathbf{N}^{\textrm{m}}_{,1}\,\,, \nonumber
	\end{align}
	\endgroup
	and for a slip step are
	\begingroup
	\begin{align}
	\frac{\partial \bm{t}^{\textrm{s}}_{\textrm{T}}}{\partial \mathbf{u}^{se}} &= \left[ -\mu_f\, \epsilon_{\text{T}} \, \bm{n}_{\textrm{T}} \otimes \bar{\bm{n}} + \frac{\mu_f \textrm{t}_{\textrm{N}}}{|| \bm{t}^{\textrm{trial}}_{\textrm{T}}||} [\bm{I} - \bm{n}_{\textrm{T}}\otimes \bm{n}_{\textrm{T}}]\, \epsilon_{\text{T}} c^{11} \,\bar{\bm{\uptau}}_1\otimes \bar{\bm{\uptau}}_1  \right] \mathbf{N}^{\textrm{s}}
	\,,~~~\textrm{and} \nonumber \\
	\frac{\partial \bm{t}^{\textrm{s}}_{\textrm{T}}}{\partial \mathbf{u}^{me}} &=  \left[ \mu_f\, \epsilon_{\text{T}} \, \bm{n}_{\textrm{T}} \otimes \bar{\bm{n}} - \frac{\mu_f \textrm{t}_{\textrm{N}}}{|| \bm{t}^{\textrm{trial}}_{\textrm{T}}||} [\bm{I} - \bm{n}_{\textrm{T}}\otimes \bm{n}_{\textrm{T}}]\, \epsilon_{\text{T}} c^{11} \,\bar{\bm{\uptau}}_1\otimes \bar{\bm{\uptau}}_1  \right]\mathbf{N}^{\textrm{m}} \\ 
	&+ ~~ \epsilon_{\text{T}} \frac{\mu_f \textrm{t}_{\textrm{N}}}{|| \bm{t}^{\textrm{trial}}_{\textrm{T}}||} [\bm{I} - \bm{n}_{\textrm{T}}\otimes \bm{n}_{\textrm{T}}]\left[ \mathbf{N}^{m} - \mathbf{N}^{m}(\xi_{sl\, n}^{\textrm{m}}) \right] \nonumber \\
	&+ ~~ \epsilon_{\text{T}}\,\textrm{g}_{\textrm{N}}c^{11} \frac{\mu_f \textrm{t}_{\textrm{N}}}{|| \bm{t}^{\textrm{trial}}_{\textrm{T}}||} [\bm{I} - \bm{n}_{\textrm{T}}\otimes \bm{n}_{\textrm{T}}] \bar{\bm{\uptau}}_1 \otimes \bar{\bm{n}} \, \mathbf{N}^{m}_{,1}\,\,. \nonumber
	\end{align}
	\endgroup
	Moreover, the derivative of closest projection point are defined as
	\begin{equation}\label{key}
	\begin{aligned}
	\frac{\partial \bar{\xi}^m}{\partial \mathbf{u}^{se}} =  c^{11} \bar{\bm{\uptau}}_1 \mathbf{N}^{\textrm{s}}\,, ~~~\textrm{and}~~~~~
	\frac{\partial \bar{\xi}^m}{\partial \mathbf{u}^{me}} = -c^{11}\left( \bar{\bm{\uptau}}_1 \mathbf{N}^{\textrm{m}} - \textrm{g}_{\textrm{N}}\bar{\bm{n}}\mathbf{N}^{\textrm{m}}_{,1} \right)\,.
	\end{aligned}
	\end{equation}
\end{appendices}

\begingroup
\fontsize{10pt}{10pt}\selectfont
\printbibliography

@book{Konyukhov2013,
	author = {Konyukhov, A. and Schweizerhof, K.},	
	title = {Computational Contact Mechanics -– Geometrically Exact Theory for Arbitrary Shaped Bodies},
	volume = {67},
	year = {2013},
	PUBLISHER = {Springer, Heildelberg},
	isbn = {978-3-642-31531-2},
	doi = {10.1007/978-3-642-31531-2},
}

@book{Konyukhov2015,
	author = {Konyukhov, Alexander and Izi, Ridvan},
	title = {Introduction to Computational Contact Mechanics{:} A Geometrical Approach},
	year = {2015},	
	PUBLISHER = {Wiley},
	ADDRESS   = {Chichester},
	isbn = {978-1-118-77065-8},
}

@article{Konyukhov2009,
	author = "Konyukhov, A. and Schweizerhof, K.",
	title = "Incorporation of contact for high--order finite elements in covariant form",
	journal = "Computer Methods in Applied Mechanics and Engineering",
	volume = "198",
	number = "13",
	pages = "1213 -- 1223",
	year = "2009",
	issn = "0045-7825",
	doi = "https://doi.org/10.1016/j.cma.2008.04.023"
}

@article{Verhoosel2010,
	author = {Verhoosel, C. V. and Scott, M. A. and de Borst, R. and Hughes, T. J. R.},
	title = {An isogeometric approach to cohesive zone modeling},
	journal = {International Journal for Numerical Methods in Engineering},
	volume = {87},
	number = {1-‐5},
	pages = {336--360},
	year = {2011},
	doi = {10.1002/nme.3061}
}

@article{Singh2018b,
	author = "Singh, S. K. and Singh, I. V. and Bhardwaj, G. and Mishra, B. K.",
	title = "A {B}\'{e}zier extraction based {XIGA} approach for three--dimensional crack simulations",
	journal = "Advances in Engineering Software",
	volume = "125",
	pages = "55 -- 93",
	year = "2018",
	issn = "0965--9978",
	doi = "https://doi.org/10.1016/j.advengsoft.2018.08.014"
}

@article{Singh2019,
	author = "Singh, S. K. and Singh, I. V. and Mishra, B. K. and Bhardwaj, G.",
	title = "Analysis of cracked functionally graded material plates using {XIGA} based on generalized higher--order shear deformation theory",
	journal = "Composite Structures",
	volume = "225",
	pages = "111038",
	year = "2019",
	issn = "0263--8223",
	doi = "https://doi.org/10.1016/j.compstruct.2019.111038"
}

@article{Singh2018a,
	author = "Singh, S. K. and Singh, I. V. and Mishra, B. K. and Bhardwaj, G. and Singh, S. K.",
	title = "Analysis of cracked plate using higher--order shear deformation theory{:} {A}symptotic crack--tip fields and {XIGA} implementation",
	journal = "Computer Methods in Applied Mechanics and Engineering",
	volume = "336",
	pages = "594 -- 639",
	year = "2018",
	issn = "0045--7825",
	doi = "https://doi.org/10.1016/j.cma.2018.03.009"
}

@article{Bhardwaj2015,
	author  = "Bhardwaj, G. and Singh, I. V. and Mishra, B. K.",
	title   = "Stochastic fatigue crack growth simulation of interfacial crack in bi--layered {FGM}s using {XIGA}",
	journal = "Computer Methods in Applied Mechanics and Engineering",
	volume = "284",
	pages = "186 -- 229",
	year = "2015",	
	issn = "0045--7825",
	doi = "https://doi.org/10.1016/j.cma.2014.08.015"
}

@article{Bhardwaj2016,
	author = "Bhardwaj, G. and Singh, S. K. and Singh, I. V. and  Mishra, B. K. and Rabczuk, T.",
	title = "Fatigue crack growth analysis of an interfacial crack in heterogeneous materials using homogenized {XIGA}",
	journal = "Theoretical and Applied Fracture Mechanics",
	volume = "85",
	pages = "294 -- 319",
	year = "2016",
	issn = "0167--8442",
	doi = "https://doi.org/10.1016/j.tafmec.2016.04.004"
}

@article{Singh2020,
	author = "Singh, S. K. and Singh, I. V.",
	title = "Analysis of cracked functionally graded piezoelectric material using {XIGA}",
	journal = "Engineering Fracture Mechanics",
	volume = "230",
	pages = "107015",
	year = "2020",
	issn = "0013--7944",
	doi = "https://doi.org/10.1016/j.engfracmech.2020.107015"
}

@ARTICLE{hughes2005,
	author={Hughes, T. J. R. and Cottrell, J. A. and Bazilevs, Y.},
	title={Isogeometric analysis: {CAD}, finite elements, {NURBS}, exact geometry and mesh refinement},
	journal={Computer Methods in Applied Mechanics and Engineering},
	year={2005},
	volume={194},
	number={39--41},
	pages={4135--4195},
	doi = "10.1016/j.cma.2004.10.008",
}

@BOOK{Cottrell2009,
	AUTHOR    = {Cottrell, J. A. and Hughes, T. J. R. and Bazilevs, Y.},
	TITLE     = {Isogeometric Analysis: Toward Integration of {CAD} and {FEA}},
	PUBLISHER = {Wiley},
	ADDRESS   = {},
	YEAR      =  "2009"
}

@ARTICLE{Padmanabhan2001,
	author={Padmanabhan, V. and Laursen, T. A.},
	title={A framework for development of surface smoothing procedures in large deformation frictional contact analysis},
	journal={Finite Elements in Analysis and Design},
	year={2001},
	volume={37},
	number={3},
	pages={173--198},
	doi = "10.1016/S0168-874X(00)00029-9",
}

@ARTICLE{El-Abbasi2001,
	author={Elabbasi, N. and Meguid, S. A. and Czekanski, A.},
	title={On the modelling of smooth contact surfaces using cubic splines},
	journal={International Journal for Numerical Methods in Engineering},
	year={2001},
	volume={50},
	number={4},
	pages={953--967},
	doi = {10.1002/1097-0207(20010210)50:4<953::AID-NME64>3.0.CO;2-P}
}

@ARTICLE{Wriggers2001,
	author={Wriggers, P. and Krstulovic-Opara, L. and Korelc, J.},
	title={Smooth ${C}^{1}$-interpolations for two\textendash dimensional frictional contact problems},
	journal={International Journal for Numerical Methods in Engineering},
	year={2001},
	volume={51},
	number={12},
	pages={1469--1495},
	doi = {10.1002/nme.227}
}

@Article{Krstulovi2002,
	author="Krstulovi{\'{c}}-Opara, L. and Wriggers, P.
	and Korelc, J.",
	title="A ${C}^1-$continuous formulation for {3D} finite deformation frictional contact",
	journal="Computational Mechanics",
	year="2002",
	volume="29",
	number="1",
	pages="27--42",
	doi="10.1007/s00466-002-0317-z"
}

@ARTICLE{Al-Dojayli2002,
	author={Al-Dojayli, M. and Meguid, S. A.},
	title={Accurate modeling of contact using cubic splines},
	journal={Finite Elements in Analysis and Design},
	year={2002},
	volume={38},
	number={4},
	pages={337--352},
	doi = "10.1016/S0168-874X(01)00088-9"
}

@ARTICLE{Stadler2003,
	author={Stadler, M. and Holzapfel, G. A. and Korelc, J.},
	title={${C}^{n}$-continuous modelling of smooth contact surfaces using {NURBS} and application to {2D} problems},
	journal={International Journal for Numerical Methods in Engineering},
	year={2003},
	volume={57},
	number={15},
	pages={2177--2203},
	doi = {10.1002/nme.776},
}

@ARTICLE{Temizer2011,
	author={Temizer, \.{I}. and Wriggers,  P. and Hughes, T. J. R.},
	title={Contact treatment in {i}sogeometric analysis with {NURBS}},
	journal={Computer Methods in Applied Mechanics and Engineering},
	year={2011},
	volume={200},
	number={9--12},
	pages={1100--1112},
	doi = "10.1016/j.cma.2010.11.020"
}

@article{Temizer2012,
	author = "Temizer, \.{I}. and Wriggers, P. and Hughes, T. J. R. ",
	title = "Three\textendash dimensional mortar\textendash based frictional contact treatment in {i}sogeometric analysis with {NURBS} ",
	journal = "Computer Methods in Applied Mechanics and Engineering ",
	volume = "209\textendash 212",
	number = "",
	pages = "115--128",
	year = "2012",
	note = "",
	doi = "10.1016/j.cma.2011.10.014"
}

@ARTICLE{Lu2011,
	author={Lu, J.},
	title={Isogeometric contact analysis: {G}eometric basis and formulation for frictionless contact},
	journal={Computer Methods in Applied Mechanics and Engineering},
	year={2011},
	volume={200},
	number={5--8},
	pages={726--741},
	doi = "10.1016/j.cma.2010.10.001",
}

@article {DeLorenzis2011,
	author = {De Lorenzis, L. and Temizer, \.{I}. and Wriggers, P. and Zavarise, G.},
	title = {A large deformation frictional contact formulation using {NURBS}\textendash based isogeometric analysis},
	journal = {International Journal for Numerical Methods in Engineering},
	volume = {87},
	number = {13},
	publisher = {John Wiley & Sons, Ltd.},
	pages = {1278--1300},
	year = {2011},
	doi = {10.1002/nme.3159}
}

@ARTICLE{DeLorenzis2012,
	author={De Lorenzis, L. and Wriggers, P. and Zavarise, G.},
	title={A mortar formulation for {3D} large deformation contact using {NURBS}\textendash based isogeometric analysis and the augmented {L}agrangian method},
	journal={Computational Mechanics},
	year={2012},
	volume={49},
	number={1},
	pages={1--20},
	doi={10.1007/s00466-011-0623-4},
	note={},
}

@article{Papadopoulos1992,
	author = {Papadopoulos, P. and Taylor, R. L.},
	title = {A {m}ixed {f}ormulation for the {f}inite {e}lement {s}olution of {c}ontact {p}roblems},
	journal = {Computer Methods in Applied Mechanics and Engineering},
	volume = {94},
	number = {3},
	year = {1992},
	issn = {0045--7825},
	pages = {373--389},
	doi = {10.1016/0045-7825(92)90061-N},
	address = {Lausanne, Switzerland}
}

@article{ALART1991,
	author = "Alart, P. and Curnier, A.",
	title = "A mixed formulation for frictional contact problems prone to {N}ewton like solution methods",
	journal = "Computer Methods in Applied Mechanics and Engineering",
	volume = "92",
	number = "3",
	pages = "353--375",
	year = "1991",
	doi = "10.1016/0045-7825(91)90022-X"
}

@article{Matzen2013,
	author = "Matzen, M. E. and Cichosz, T. and Bischoff,  M. ",
	title = "A point to segment contact formulation for isogeometric, {NURBS} based finite elements ",
	journal = "Computer Methods in Applied Mechanics and Engineering ",
	volume = "255",
	number = "",
	pages = "27--39",
	year = "2013",
	note = "",
	doi = "10.1016/j.cma.2012.11.011"
}

@article{Matzen2016,
	author = "Matzen, M. E. and Bischoff,  M.",
	title = "A weighted point\textendash based formulation for isogeometric contact ",
	journal = "Computer Methods in Applied Mechanics and Engineering ",
	volume = "308",
	number = "",
	pages = "73--95",
	year = "2016",
	note = "",
	issn = "",
	doi = "10.1016/j.cma.2016.04.010"
}

@ARTICLE{Kim2012,
	author={Kim, J. Y. and Youn, S. K.},
	title={Isogeometric contact analysis using mortar method},
	journal={International Journal for Numerical Methods in Engineering},
	year={2012},
	volume={89},
	number={12},
	pages={1559--1581},	
	doi = {10.1002/nme.3300}
}

@article{Temzier2013Multiscale,
	author = {Temizer, \.{I}. },
	title = {Multiscale thermomechanical contact: {C}omputational homogenization with isogeometric analysis},
	journal = {International Journal for Numerical Methods in Engineering},
	year   = {2013},
	volume = {97},
	number = {8},
	pages = {582--607},
	doi = {10.1002/nme.4604}
}

@article{Dittmann2014,
	author = " Dittmann, M. and Franke, M. and Temizer, \.{I}. and Hesch, C.",
	title = "Isogeometric Analysis and thermomechanical Mortar contact problems",
	journal = "Computer Methods in Applied Mechanics and Engineering",
	volume = "274",
	pages = "192--212",
	year = "2014",
	issn = "",
	doi = "10.1016/j.cma.2014.02.012"
}

@article{BRIVADIS2015,
	author = "Brivadis, E. and Buffa, A. and Wohlmuth, B. and Wunderlich, L.",
	title = "Isogeometric mortar methods",
	journal = "Computer Methods in Applied Mechanics and Engineering",
	volume = "284",
	pages = "292--319",
	year = "2015",
	issn = "0045--7825",
	doi = "10.1016/j.cma.2014.09.012"
}

@article{SEITZ2016,
	author = "Seitz, A. and Farah, P. and Kremheller, J. and Wohlmuth, B. I. and  Wall, W. A. and Popp, A. ",
	title = "Isogeometric dual mortar methods for computational contact mechanics",
	journal = "Computer Methods in Applied Mechanics and Engineering",
	volume = "301",
	pages = "259--280",
	year = "2016",
	doi = "10.1016/j.cma.2015.12.018"
}

@Article{Duong2018,
	author="Duong, T. X. and De Lorenzis, L. and Sauer, R. A.",
	title="A segmentation-free isogeometric extended mortar contact method",journal="Computational Mechanics",
	year="2019",
	volume="63",
	number="2",
	pages="383--407",	
	doi="10.1007/s00466-018-1599-0"
}

@article{KRUSE2015,
	author = "Kruse, R. and Nguyen-Thanh, N.  and De Lorenzis, L. and Hughes, T. J. R.",
	title = "Isogeometric collocation for large deformation elasticity and frictional contact problems",
	journal = "Computer Methods in Applied Mechanics and Engineering",
	volume = "296",
	pages = "73--112",
	year = "2015",
	issn = "0045--7825",
	doi = "10.1016/j.cma.2015.07.022"
}

@article{DELORENZIS2015,
	author = "De Lorenzis, L. and Evans, J. A. and Hughes, T. J. R. and Reali, A.",
	title = "Isogeometric collocation: {N}eumann boundary conditions and contact",
	journal = "Computer Methods in Applied Mechanics and Engineering",
	volume = "284",
	pages = "21--54",
	year = "2015",
	issn = "0045--7825",
	doi = "10.1016/j.cma.2014.06.037"
}

@Article{Weeger2018, 
	author="Weeger, O. and Narayanan, B. and Dunn, M. L.",
	title="Isogeometric collocation for nonlinear dynamic analysis of {C}osserat rods with frictional contact",
	journal="Nonlinear Dynamics",
	year="2018",
	volume="91",
	number="2",
	pages="1213--1227",
	doi="10.1007/s11071-017-3940-0"
}

@ARTICLE{DeLorenzis2014,
	author={De Lorenzis, L. and Wriggers, P. and Hughes, T. J. R.},
	title={Isogeometric contact: {A} review},
	journal={GAMM Mitteilungen},
	year={2014},
	volume={37},
	number={1},
	pages={85--123},
	doi = {10.1002/gamm.201410005}
}

@article{BAZILEVS2010Tspline,
	author = "Bazilevs, Y. and Calo, V. M. and Cottrell, J. A.  and Evans, J. A. and Hughes, T. J. R. and Lipton, S. and  Scott, M. A. and Sederberg, T. W.",
	title = "Isogeometric analysis using {T-}splines",
	journal = "Computer Methods in Applied Mechanics and Engineering",
	volume = "199",
	number = "5",
	pages = "229--263",
	year = "2010",
	doi = "10.1016/j.cma.2009.02.036"
}

@article{DORFEL2010,
	author = "D$\ddot{\text{o}}$rfel, M. R. and J$\ddot{\text{u}}$ttler, B. and Simeon, B.",
	title = "Adaptive isogeometric analysis by local h{-}refinement with {T-}splines",
	journal = "Computer Methods in Applied Mechanics and Engineering",
	volume = "199",
	number = "5",
	pages = "264--275",
	year = "2010",
	doi = "10.1016/j.cma.2008.07.012"
}

@article{Scott2011,
	author = {Scott, M. A. and Borden, M. J. and Verhoosel, C. V. and Sederberg, T. W. and Hughes, T. J. R.},
	title = {Isogeometric finite element data structures based on {B}\'{e}zier extraction of {T-}splines},
	journal = {International Journal for Numerical Methods in Engineering},
	volume = {88},
	number = {2},
	pages = {126--156},
	year = {2011},
	doi = {10.1002/nme.3167}
}

@article{Dimitri2014Tspline,
	author = "Dimitri, R. and De Lorenzis, L. and Scott, M. A. and Wriggers, P. and Taylor, R. L. and Zavarise, G.",
	title = "Isogeometric large deformation frictionless contact using {T-}splines",
	journal = "Computer Methods in Applied Mechanics and Engineering",
	volume = "269",
	pages = "394--414",
	year = "2014",	
	doi = "10.1016/j.cma.2013.11.002"
}

@article{Dimitri2014,
	author = {Dimitri, R. and De Lorenzis, L. and Wriggers, P. and Zavarise, G.},
	title = {{NURBS}\textendash and {T}\textendash spline-based isogeometric {c}ohesive {z}one {m}odeling of {i}nterface {d}ebonding},
	journal = {Computational Mechanics},
	volume = {54},
	number = {2},
	year = {2014},
	pages = {369--388},
	numpages = {20},
	publisher = {Springer-Verlag New York, Inc.},
	address = {Secaucus, NJ, USA},
	doi = {10.1007/s00466-014-0991-7}
}

@Article{Dimitri2017,
	author="Dimitri, R. and Zavarise, G.",
	title="Isogeometric treatment of frictional contact and mixed mode debonding problems",
	journal="Computational Mechanics",
	year="2017",
	volume="60",
	number="2",
	pages="315--332",
	doi="10.1007/s00466-017-1410-7"
}

@ARTICLE{Temizer20161,
	author={Temizer, \.{I}. and Hesch, C.},
	title={Hierarchical {NURBS} in frictionless contact},
	journal={Computer Methods in Applied Mechanics and Engineering},
	year={2016},
	volume={299},
	pages={161--186},
	doi={10.1016/j.cma.2015.11.006}
}

@article{Temizer2014,
	author = "Temizer, \.{I}. and Abdalla, M. M. and Gurdal, Z.",
	title = "An interior point method for isogeometric contact",
	journal = "Computer Methods in Applied Mechanics and Engineering",
	volume = "276",
	pages = "589--611",
	year = "2014",
	issn = "0045-7825",
	doi = "https://doi.org/10.1016/j.cma.2014.03.018"
}

@article{Temizer20162,
	author = "Hesch, C. and Franke, M. and Dittmann, M. and Temizer, \.{I}.",
	title = "Hierarchical {NURBS} and a higher--order phase--field approach to fracture for finite-deformation contact problems",
	journal = "Computer Methods in Applied Mechanics and Engineering",
	volume = "301",
	pages = "242--258",
	year = "2016",
	issn = "0045--7825",
	doi = "10.1016/j.cma.2015.12.011"
}

@Article{Sauer2017,
	author = {{Zimmermann}, C. and {Sauer}, R. A.},
	title = "{Adaptive local surface refinement based on {LR} {NURBS} and its application to contact}",
	journal="Computational Mechanics",
	year="2017",
	volume="60",
	number="6",
	pages="1011--1031",
	doi="10.1007/s00466-017-1455-7"
}

@article{Corbett2014,
	author = "Corbett, C. J. and Sauer, R. A. ",
	title = "{NURBS}\textendash enriched contact finite elements ",
	journal = "Computer Methods in Applied Mechanics and Engineering ",
	volume = "275",
	number = "",
	pages = "55--75",
	year = "2014",
	note = "",
	doi = "10.1016/j.cma.2014.02.019"
}

@article{Corbett2015,
	author = "Corbett, C. J. and Sauer, R. A.",
	title = "Three\textendash dimensional {i}sogeometrically enriched finite elements for frictional contact and mixed\textendash mode debonding ",
	journal = "Computer Methods in Applied Mechanics and Engineering ",
	volume = "284",
	number = "",
	pages = "781--806",
	year = "2015",
	doi = "10.1016/j.cma.2014.10.025"
}

@article{RASOOL2016182,
	author = "Rasool, R. and Corbett, C. J. and Sauer, R. A.",
	title = "A strategy to interface isogeometric analysis with Lagrangian finite elements-—Application to incompressible flow problems",
	journal = "Computers \& Fluids",
	volume = "127",
	pages = "182--193",
	year = "2016",
	doi = "10.1016/j.compfluid.2015.12.016"
}

@article{MALEKIJEBELI2018,
	author = "Maleki-Jebeli, S.  and Mosavi-Mashhadi, M. and Baghani, M.",
	title = "A large deformation hybrid isogeometric--finite element method applied to cohesive interface contact/debonding",
	journal = "Computer Methods in Applied Mechanics and Engineering",
	volume = "330",
	pages = "395--414",
	year = "2018",
	doi = "10.1016/j.cma.2017.10.017"
}

@Article{Otto2018,
	author="Otto, P. and De Lorenzis, L. and Unger, J. F.",
	title="Coupling a {NURBS} contact interface with a higher order finite element discretization for contact problems using the mortar method",
	journal="Computational Mechanics",
	year="2019",
	volume="63",
	number="6",
	pages="1203--1222",
	doi="10.1007/s00466-018-1645-y"
}

@Article{DIAS2019,
	author="Dias, A. P. C. and Proenca, S. P. B. and Bittencourt, M. L.",
	title="High{-}order mortar{-}based contact element using {NURBS} for the mapping of contact curved surfaces",
	journal="Computational Mechanics",
	year="2019",
	volume="64",
	number="1",
	pages="85--112",
	doi="10.1007/s00466-018-1658-6"
}

@article {Roger2011_enrichment,
	author = {Sauer, R. A.},
	title = {Enriched contact finite elements for stable peeling computations},
	journal = {International Journal for Numerical Methods in Engineering},
	volume = {87},
	number = {6},
	publisher = {John Wiley & Sons, Ltd.},
	year = {2011},
	doi = {10.1002/nme.3126}
}

@article{Sauer2013,
	author = {Sauer, R. A.},
	title = {Local finite element enrichment strategies for 2{D} contact computations and a corresponding post--processing scheme},
	journal = {Computational Mechanics},
	volume = {52},
	number = {2},
	year = {2013},
	pages = {301--319},
	numpages = {19},
	doi = {10.1007/s00466-012-0813-8}
}

@Article{Fischer2005,
	author="Fischer, K. A. and Wriggers, P.",
	title="Frictionless 2{D} Contact formulations for finite deformations based on the mortar method",
	journal="Computational Mechanics",
	year="2005",
	volume="36",
	number="3",
	pages="226--244",
	doi="10.1007/s00466-005-0660-y"
}

@article{PUSO2004,
	author = "Puso, M. A. and Laursen, T. A.",
	title = "A mortar segment{-}to{-}segment contact method for large deformation solid mechanics",
	journal = "Computer Methods in Applied Mechanics and Engineering",
	volume = "193",
	number = "6",
	pages = "601--629",
	year = "2004",
	doi = "10.1016/j.cma.2003.10.010"
}

@book{laursen2003,
	title={Computational {C}ontact and {I}mpact {M}echanics: {F}undamentals of {M}odeling {I}nterfacial {P}henomena in {N}onlinear {F}inite {E}lement {A}nalysis},
	author={Laursen, T. A.},
	isbn={9783540429067},
	lccn={2002511027},
	year={2002},
	publisher={Springer, Berlin}
}

@BOOK{wriggers2006,
	AUTHOR    = {Wriggers, P.},
	TITLE     = {Computational {C}ontact {M}echanics (2nd edn.)},
	PUBLISHER = {Springer, Berlin},
	ADDRESS   = {},
	YEAR      =  "2006"
}

@article{unbiased2013,
	author = "Sauer, R. A. and De Lorenzis, L.",
	title = "A computational contact formulation based on surface potentials",
	journal = "Computer Methods in Applied Mechanics and Engineering",
	volume = "253",
	pages = "369--395",
	year = "2013",
	doi = "10.1016/j.cma.2012.09.002"
}

@article{unbiased2015,
	author = "Sauer, R. A. and De~Lorenzis, L.",
	title = "An unbiased computational contact formulation for {3D} friction",
	journal = "International Journal for Numerical Methods in Engineering",
	volume = {101},
	number = {4},
	pages = {251--280},
	year = {2015},
	doi = {10.1002/nme.4794}
}

@BOOK{nurbsbook,
	AUTHOR    = {Piegl, L. and Tiller, W.},
	TITLE     = {The {NURBS} book (Monographs in Visual Communication)},
	SERIES    = {},
	PUBLISHER = {Springer, Berlin Heidelberg},
	ADDRESS   = {},
	YEAR      =  "2012"
}

@article {Bezier_extraction_NURBS,
	author = {Borden, M. J. and Scott, M. A. and Evans, J. A. and Hughes, T. J. R.},
	title = {Isogeometric finite element data structures based on {B}\'{e}zier extraction of {NURBS}},
	journal = {International Journal for Numerical Methods in Engineering},
	volume = {87},
	number = {1--5},
	publisher = {John Wiley \& Sons, Ltd.},
	pages = {15--47},
	year = {2011},
	doi = {10.1002/nme.2968}
}

@ARTICLE{Hughes2010,
	author={Hughes, T. J. R. and Reali, A. and Sangalli, G.},
	title={Efficient quadrature for {NURBS}\textendash based isogeometric analysis},
	journal={Computer Methods in Applied Mechanics and Engineering},
	year={2010},
	volume={199},
	number={5--8},
	pages={301--313},
	doi = "10.1016/j.cma.2008.12.004"
}

@article{AURICCHIO2012,
	author = "Auricchio, F. and Calabr\`{o}, F. and Hughes, T. J .R. and Reali, A. and Sangalli, G.",
	title = "A simple algorithm for obtaining nearly optimal quadrature rules for NURBS{-}based isogeometric analysis",
	journal = "Computer Methods in Applied Mechanics and Engineering",
	volume = "249--252",
	pages = "15--27",
	year = "2012",
	doi = "10.1016/j.cma.2012.04.014"
}

@article{FAHRENDORF2018,
	author = "Fahrendorf, F. and De Lorenzis, L. and Gomez, H.",
	title = "Reduced integration at superconvergent points in isogeometric analysis",
	journal = "Computer Methods in Applied Mechanics and Engineering",
	volume = "328",
	pages = "390--410",
	year = "2018",
	doi = "10.1016/j.cma.2017.08.028"
}

@Article{Agrawal2018,
	author="Agrawal, V. and Gautam, S. S.",
	title="{IGA}: {A} {S}implified {I}ntroduction and {I}mplementation {D}etails for {F}inite {E}lement {U}sers",
	journal="Journal of The Institution of Engineers (India): Series C.~~",
	year="2019",
	volume="100",
	number="3",
	pages="561--585",
	doi="10.1007/s40032-018-0462-6"
}

@BOOK{johnson1987,
	AUTHOR    = {Johnson, K. L.},
	TITLE     = {{C}ontact {M}echanics},
	PUBLISHER = {Cambridge University Press},
	ADDRESS   = {Cambridge},
	YEAR      =  "1987"
}

@Article{Franke2010,
	author="Franke, D. and D{\"u}ster, A. and N{\"u}bel, V. and Rank, E.",
	title="A comparison of the $h-,~ p-,~ hp-,$ and $rp-$ version of the {FEM} for the solution of the {2D} {H}ertzian contact problem",
	journal="Computational Mechanics",
	year="2010",
	volume="45",
	number="5",
	pages="513--522",
	doi="10.1007/s00466-009-0464-6"
}

@book{bonet1997,
	author={Bonet, J. and Wood, R. D.},  
	title={Nonlinear Continuum Mechanics for Finite Element Analysis},
	isbn={9780521572729},
	lccn={97011366},
	year={1997},
	publisher={Cambridge University Press}
}

@book{Kikuchi1988,
	author="Kikuchi, N. and Oden, J. T.",
	title="Contact {P}roblems in {E}lasticity: {A} {S}tudy of {V}ariational {I}nequalities and {F}inite {E}lement {M}ethods",
	isbn={},
	lccn={},
	year={1998},
	publisher={Philadelphia: SIAM}
}

@article{Otto2019,
	author = {Otto, P. and De Lorenzis, L. and Unger, J. F.},
	title   = {Explicit dynamics in impact simulation using a NURBS contact interface},
	journal = {International Journal for Numerical Methods in Engineering},	
	volume = {},
	number = {},
	pages = {},
	year  = {2019},
	doi = {10.1002/nme.6264}
}

@article{Popp2010,
	author = {Popp, A. and Gitterle, M. and Gee, M. W. and Wall, W. A.},
	title = {A dual mortar approach for {3D} finite deformation contact with consistent linearization},
	journal = {International Journal for Numerical Methods in Engineering},
	volume = {83},
	number = {11},
	pages = {1428-1465},
	year = {2010},
	doi = {10.1002/nme.2866}
}

@article{Popp2013,
	author = "Popp, A. and Seitz, A. and Gee, M. W. and Wall, W. A.",
	title = "Improved robustness and consistency of 3D contact algorithms based on a dual mortar approach",
	journal = "Computer Methods in Applied Mechanics and Engineering",
	volume = "264",
	pages = "67 - 80",
	year = "2013",
	issn = "0045-7825",
	doi = "https://doi.org/10.1016/j.cma.2013.05.008"
}

@book{Popp2018,
	author="Popp, A.", 
	editor="Popp, A. and Wriggers, P.",
	title = "State-of-the-Art Computational Methods for Finite Deformation Contact Modeling of Solids and  Structures",
	booktitle = "Contact Modeling for Solids and Particles",
	publisher="Springer International Publishing",
	address="",
	pages="1--86",
	year="2018",
	isbn="978-3-319-90155-8",
	doi="10.1007/978-3-319-90155-8_1"
}

@article{Zavarise2009,
	author = {Zavarise, G. and De Lorenzis, L.},
	title = {A modified node{-}to{-}segment algorithm passing the contact patch test},
	journal = {International Journal for Numerical Methods in Engineering},
	year = {2009},
	volume = {79},
	number = {4},
	pages = {379-416},
	doi = {10.1002/nme.2559}
}
\endgroup
\end{document}